\documentclass[twoside,11pt]{article}

\usepackage{jmlr2e}

\usepackage{lipsum}
\usepackage{amsfonts}
\usepackage{graphicx}
\usepackage{epstopdf}
\usepackage{mathtools}
\usepackage{empheq}
\usepackage{amsfonts}
\usepackage{amsgen,amsmath,amstext,amsbsy,amsopn,amssymb,subfigure,stmaryrd,mathabx}
\usepackage{comment}
\usepackage{enumitem}
\usepackage{booktabs}
\usepackage{blkarray}
\usepackage{algorithm}
\usepackage{algpseudocode}
\usepackage{url}
\usepackage{longtable}
\usepackage{multirow}
\usepackage{xcolor}

\ifpdf
\DeclareGraphicsExtensions{.eps,.pdf,.png,.jpg}
\else
\DeclareGraphicsExtensions{.eps}
\fi

\ifpdf
\DeclareGraphicsExtensions{.eps,.pdf,.png,.jpg}
\else
\DeclareGraphicsExtensions{.eps}
\fi

\newtheorem{Theorem}{Theorem}

\newtheorem{Lemma}{Lemma}
\newtheorem{Remark}{Remark}
\newtheorem{Corollary}{Corollary}

\DeclareMathAlphabet\mathbfcal{OMS}{cmsy}{b}{n}

\newcommand{\be}{\begin{equation}}
\newcommand{\ee}{\end{equation}}
\newcommand{\bea}{\begin{eqnarray}}
\newcommand{\eea}{\end{eqnarray}}
\newcommand{\beas}{\begin{eqnarray*}}
	\newcommand{\eeas}{\end{eqnarray*}}

\newcommand{\bbR}{\mathbb{R}}

\newcommand{\cS}{\mathcal{S}}

\newcommand{\cM}{\mathcal{M}}

\newcommand{\0}{{\mathbf{0}}}

\renewcommand{\u}{{\mathbf{u}}}
\renewcommand{\v}{{\mathbf{v}}}

\renewcommand{\r}{{\mathbf{r}}}

\newcommand{\X}{{\mathbf{X}}}
\newcommand{\B}{{\mathbf{B}}}
\newcommand{\J}{{\mathbf{J}}}
\newcommand{\E}{{\mathbf{E}}}

\newcommand{\D}{{\mathbf{D}}}

\newcommand{\I}{{\mathbf{I}}}
\newcommand{\M}{{\mathbf{M}}}

\newcommand{\R}{{\mathbf{R}}}

\newcommand{\U}{{\mathbf{U}}}
\newcommand{\V}{{\mathbf{V}}}
\newcommand{\W}{{\mathbf{W}}}

\newcommand{\A}{{\mathbf{A}}}

\newcommand{\Z}{{\mathbf{Z}}}
\newcommand{\T}{{\mathbf{T}}}

\newcommand{\bPi}{{\mathbf{\Pi}}}

\newcommand{\bcY}{{\mathbfcal{Y}}}
\newcommand{\bcB}{{\mathbfcal{B}}}
\newcommand{\bcS}{{\mathbfcal{S}}}

\newcommand{\bcZ}{{\mathbfcal{Z}}}
\newcommand{\bcT}{{\mathbfcal{T}}}

\newcommand{\bcI}{{\mathbfcal{I}}}

\newcommand{\Proj}{P}

\newcommand{\bSigma}{\boldsymbol{\Sigma}}
\newcommand{\bLambda}{\boldsymbol{\Lambda}}

\newcommand{\rank}{{\rm rank}}

\newcommand{\tr}{{\rm tr}}

\newcommand{\SVD}{{\rm SVD}}

\newcommand{\tHS}{{\rm HS}}
\newcommand{\HOSVD}{{\rm T-HOSVD}}

\newcommand{\argmin}{\mathop{\rm arg\min}}

\newcommand{\bbP}{\mathbb{P}}
\newcommand{\bbE}{\mathbb{E}}

\makeatletter
\newcommand*{\rom}[1]{\expandafter\@slowromancap\romannumeral #1@}
\makeatother 
\allowdisplaybreaks

\usepackage{lastpage}
\jmlrheading{22}{2021}{1-\pageref{LastPage}}{8/20; Revised
1/21}{6/21}{20-919}{Yuetian Luo, Garvesh Raskutti, Ming Yuan, and Anru R. Zhang}
\ShortHeadings{title}{Luo, Raskutti, Yuan, and Zhang}

\ShortHeadings{Tensor Perturbation Bounds for Orthogonal Iteration}{Y. Luo, G. Raskutti, M. Yuan, and A. R. Zhang}
\firstpageno{1}

\begin{document}

\title{A Sharp Blockwise Tensor Perturbation Bound for Orthogonal Iteration}

\author{\name Yuetian Luo \email yluo86@wisc.edu\\
       \addr Department of Statistics\\
       University of Wisconsin-Madison\\
       Madison, WI 53706, USA
       \AND
       \name Garvesh Raskutti \email raskutti@stat.wisc.edu \\
       \addr Department of Statistics\\
       University of Wisconsin-Madison\\
       Madison, WI 53706, USA
   	   \AND
        \name Ming Yuan \email ming.yuan@columbia.edu \\
		\addr Department of Statistics\\
		Columbia University\\
		New York, NY 10027, USA       
		\AND
        \name Anru R. Zhang \email anruzhang@stat.wisc.edu \\
		\addr Department of Statistics\\
		University of Wisconsin-Madison\\
		Madison, WI 53706, USA
}

\editor{Animashree Anandkumar}

\maketitle

\begin{abstract}%
	In this paper, we develop novel perturbation bounds for the higher-order orthogonal iteration (HOOI) \citep{de2000best}. Under mild regularity conditions, we establish blockwise tensor perturbation bounds for HOOI with guarantees for both tensor reconstruction in Hilbert-Schmidt norm $\|\widehat{\bcT} - \bcT \|_{\tHS}$ and mode-$k$ singular subspace estimation in Schatten-$q$ norm $\| \sin \Theta (\widehat{\U}_k, \U_k) \|_q$ for any $q \geq 1$. We show the upper bounds of mode-$k$ singular subspace estimation are unilateral and converge linearly to a quantity characterized by blockwise errors of the perturbation and signal strength. For the tensor reconstruction error bound, we express the bound through a simple quantity $\xi$, which depends only on perturbation and the multilinear rank of the underlying signal. Rate matching deterministic lower bound for tensor reconstruction, which demonstrates the optimality of HOOI, is also provided. Furthermore, we prove that one-step HOOI (i.e., HOOI with only a single iteration) is also optimal in terms of tensor reconstruction and can be used to lower the computational cost. The perturbation results are also extended to the case that only partial modes of $\bcT$ have low-rank structure. We support our theoretical results by extensive numerical studies. Finally, we apply the novel perturbation bounds of HOOI on two applications, tensor denoising and tensor co-clustering, from machine learning and statistics, which demonstrates the superiority of the new perturbation results.
\end{abstract}


\begin{sloppypar}
	\section{Introduction}\label{sec:intro}
	
	The past decades have seen a large body of work on tensors or multiway arrays \citep{kolda2009tensor, sidiropoulos2017tensor, cichocki2015tensor,kroonenberg2008applied}. Tensors arise in numerous applications involving multiway data (e.g., brain imaging \citep{zhou2013tensor}, hyperspectral imaging \citep{li2010tensor}, recommender system design \citep{bi2018multilayer}). In addition, various methods have been developed and applied to fundamental tensor problems such as tensor completion \citep{yuan2016tensor,xia2017statistically,yuan2017incoherent,zhang2019cross}, tensor regression \citep{zhang2020islet,zhou2013tensor,raskutti2019convex,chen2019non}, tensor PCA/SVD \citep{zhang2018tensor,richard2014statistical,liu2017characterizing}, generalized tensor learning \citep{han2020optimal}. In many other problems where the observations are not necessarily tensors, tensor parameters arise, such as topic and latent variable models \citep{anandkumar2014tensor}, additive index models \citep{balasubramanian2018tensor}, high-order interaction pursuit \citep{hao2020sparse}. We refer readers to recent survey papers \citep{kolda2009tensor, sidiropoulos2017tensor, cichocki2015tensor}.
	
	Among these methods, tensor decomposition is one of the most important and a flurry of research have been devoted to it. Tensor decomposition plays a similar role to singular value decomposition (SVD) or eigendecomposition for matrices which is of fundamental importance throughout a wide range of fields including computer science, applied mathematics, machine leaning, statistics, signal processing, etc. For matrices, truncated SVD achieves the best low-rank approximation in terms of any unitarily invariant norm by the well-known Eckart-Young-Mirsky theorem \citep{eckart1936approximation,mirsky1960symmetric,golub1987generalization} and more importantly it is computationally efficient. Despite the well-established theory for low-rank decomposition of matrices, tensors present unique challenges. First there are several notions of low-rankness in tensors, moreover it has been shown that computing various best low-rank approximations of a given tensor is NP hard in general \citep{hillar2013most}. 
	
	Fortunately, many computationally efficient algorithms have been proposed to approximate the best low-rank tensor decomposition \citep{kroonenberg1980principal,de2000best,elden2009newton,ishteva2011best,ishteva2009differential,savas2010quasi,de2000multilinear,vannieuwenhoven2012new}. One popular choice is the higher-order orthogonal iteration (HOOI) proposed in \cite{de2000best}. 
	HOOI is based on alternating least-squares. It can be seen as a ``spectral" algorithm for tensors, generalizations of the 2D-PCA \citep{sheehan2007higher}, and the power iteration refinement of HOSVD \citep{de2000multilinear} and sequential HOSVD \citep{vannieuwenhoven2012new}. Convergence properties of HOOI have been studied in \cite{zhang2001rank,wang2014global,uschmajew2015new,xu2018convergence,zhang2018tensor}. 
	
	In addition to computing low-rank approximations of matrices and tensors, there is the more nuanced question of computing low-rank approximations under noise perturbation and determining how the perturbation impacts the quality of the decomposition. 
	For matrices, perturbation theory is well studied and a number of results has been established \citep{bhatia2013matrix,stewart1990matrix}. However, perturbation theory for tensors is still in its infancy. It is difficult to extend results for matrices to tensors due to the complexity of tensor algebra and the fact many well established theories and concepts in matrices such as SVD or eigendecomposition do not exist or are not easy to use for tensors. There are several attempts in the literature and most of them require the noise tensor to be random \citep{richard2014statistical,zhang2018tensor,liu2017characterizing,zhang2019optimal,zheng2015interpolating,xia2019sup}. In this paper, we move one step forward in this direction and provide the first general perturbation bounds of HOOI for tensors under the signal-plus-noise model
	\begin{equation}\label{eq: perturbation model}
	\widetilde{\bcT} = \bcT + \bcZ
	\end{equation}
	without putting any structural assumption on the perturbation noise $\bcZ$. Like the classical Wedin's perturbation theory for matrices, we provide perturbation guarantees of estimated mode-$k$ singular subspace in tensors. In addition, we also provide the perturbation bounds for tensor reconstruction. By providing the deterministic rate matching lower bound, we can further show HOOI with good initialization is optimal for tensor reconstruction.
	
	\subsection{Problem Statement} \label{sec: problem statement}
	Formally, this paper considers the tensor perturbation model \eqref{eq: perturbation model}, where $\bcT$ is the low-rank order-$d$ signal tensor and $\bcZ$ is the perturbation tensor with the same dimension as $\bcT$. Two popular choices of low-rankness in tensors are canonical polyadic (CP) low-rank and multilinear/Tucker low-rank and each of them has their respective benefits (see the surveys \cite{kolda2009tensor,sidiropoulos2017tensor,cichocki2015tensor,grasedyck2013literature}). The CP low-rank decomposition which approximates the original tensor by a sum of rank-1 outer products gives a compact and unique (under certain conditions) representation and multilinear/Tucker low-rank decomposition generally finds a better fit for the data by estimating the subspace of each mode. Since any CP low-rank tensor can be written as a multilinear low-rank tensor with a diagonal core tensor, we focus on the latter (Tucker rank) in this paper. Specifically, we assume $\bcT$ admits the following low multilinear rank (Tucker) decomposition:
	\begin{equation} \label{eq: decomposition of T}
	\bcT = \bcS \times_{\Omega_1} \U_1 \times \cdots \times_{\Omega_m} \U_m.
	\end{equation} 
	Here, $\bcS$ is an order-$d$ core tensor; $\{\U_i\}_{i=1}^m$ are group-$i$ singular subspaces; $\{\Omega_i \}_{i=1}^m$ are symmetric index groups which will be introduced next; and ``$\times_{\Omega_i}$" is the tensor matrix product along modes in $\Omega_i$. The formal definition of Tucker decomposition and tensor matrix product will be given in equations \eqref{eq:tucker-decomposition}, \eqref{eq: tensor-matrix product} and \eqref{eq: tensor matrix prod short notation} in Section \ref{sec: notation}.
	
	The symmetric index groups $\{\Omega_i \}_{i=1}^m$ represent the specific symmetricity structure of $\bcT$ and satisfy $\Omega_i \subseteq [d], \Omega_i \neq \emptyset$, $\Omega_i \bigcap \Omega_j = \emptyset$ for $i \neq j$, and $\bigcup_{i=1}^m \Omega_i = [d]$, where $[d] := \{1,\ldots,d \}$. It means by fixing indices outside the group and any permutation of indices within the group does not change the value of tensor corresponding to those indices. For example, if $\Omega_1 = \{1,2,\ldots, k \}$, then fixing coordinates in $\{\Omega_i \}_{i=2}^m$ i.e., coordinates $k+1, \ldots, d$ and for any permutation $\delta$ of $[k]$, we have $\bcT_{[i_{\delta(1)},\ldots, i_{\delta(k)},i_{k+1}, \ldots, i_d]} = \bcT_{[i_1, \ldots, i_d]}$. In addition to $\bcT$, we also assume $\bcS$ and $\bcZ$ have the same symmetric structure characterized by $\{\Omega_i \}_{i=1}^m$. 
	The symmetric index groups have two extreme cases:
	\begin{itemize}
		\item Asymmetric: $\Omega_i = \{ i\}$ for $i = 1, \ldots, d$,
		\item Supersymmetric: $\Omega_1 = \{1, 2, \ldots, d \}$. 
	\end{itemize}
	To simplify the notation, if mode $j \in \Omega_i$, then we denote $j' := i$ as the group index of mode $j$. For symmetric group $i$, the dimension of $\bcT$ in this group is denoted as $p_i$ and its mode rank is denoted as $r_i$. The rigorous definition of mode rank is postponed to Section \ref{sec: notation}. Also throughout the paper, for $i \in [m]$, let $\bar{i} = \min\{j: j \in \Omega_i\}$ be the smallest index in $\Omega_i$, and denote $ \widecheck{\Omega}_i := \Omega_i \setminus \{\bar{i}\}, \underline{\Omega}_i :=\bigcup_{j=1}^{i-1} \Omega_j, \overline{\Omega}_i := \bigcup_{j=i+1}^d \Omega_j$ with $\underline{\Omega}_1 = \overline{\Omega}_d = \emptyset$. Finally, we define the signal strength of $\bcT$ as $\lambda = \min_{i=1, \ldots,m} \sigma_{r_i} (\mathcal{M}_{\bar{i}}(\bcT))$, the smallest singular values of matricization of $\bcT$ in modes $\{\bar{i}\}_{i =1}^m$. Here for any matrix $\D$, $\sigma_i(\D)$ denotes the $i$th largest singular values of $\D$ and $\mathcal{M}_{\bar{i}}(\bcT)$ represents the matricization of tensor $\bcT$ along mode $\bar{i}$ and its formal definition will be given in \eqref{eq: tensor-matricization} in Section \ref{sec: notation}.
	
	In summary, for the dimensions of the perturbation model \eqref{eq: perturbation model}, we have $\widetilde{\bcT}, \bcT, \bcZ \in \bbR^{p_{1'} \times \cdots \times p_{d'}}$ with symmetric index groups $\{\Omega_i \}_{i=1}^m$. The HOOI algorithm we study is provided in Algorithm \ref{alg: HOOI}. It is worth noting the original HOOI algorithm in \cite{de2000best} mainly focuses on asymmetric tensor decomposition and we generalize it to accommodate arbitrary symmetric structures of $\bcT$ characterized by $\{\Omega_i \}_{i=1}^m$. In addition, in the literature \citep{de2000best,kolda2009tensor}, HOOI often refers to the overall procedure including both the initialization of $\widetilde\U_i^{(0)}$ by HOSVD and the orthogonal iteration updates as detailed in Algorithm \ref{alg: HOOI}. We shall point out that this paper studies the orthogonal iteration with any initializers satisfying some mild conditions. Our results accommodate different types of initialization scheme and can be applied to a wide range of scenarios (see Remark \ref{rem: initialization}). We conclude the algorithm by remarking that the matricization mode we choose in group $i$ to perform SVD in \eqref{eq: SVD in algorithm} does not matter due to symmetry. For simplicity, we choose the smallest index in group $i$.
	
	\begin{algorithm}[!h] \caption{Higher-Order Orthogonal Iteration for Tensor Decomposition} \label{alg: HOOI}
		\textbf{Input:} $\widetilde{\bcT}\in \bbR^{p_{1'} \times \cdots \times p_{d'}}$, symmetric index groups $\{\Omega_i \}_{i=1}^m$, initialization $\{\widetilde{\U}_i^{(0)}\}_{i=1}^m$ with $\widetilde{\U}_i^{(0)} \in \mathbb{O}_{p_i, r_i}$ having orthonormal columns, maximum number of iterations $t_{\max}$.\\
		\textbf{Output:} $ \{\widehat{\U}_i\}_{i=1}^m, \widehat{\bcT}$. 
		\begin{algorithmic}[1]
			\State For $t=1, \ldots, t_{\max}$, do
			\begin{enumerate}[label=(\alph*)]
				\item For $i = 1$ to $m$, update
				\begin{equation} \label{eq: SVD in algorithm}
				\begin{split}
				\widetilde{\U}_i^{(t+1)} =\SVD_{r_i} &  \big(\cM_{\bar{i}}( \widetilde{\bcT} \times_{\Omega_1} ( \widetilde{\U}_{1}^{(t+1)})^\top \times \cdots\\
				& \times_{\Omega_{i-1}} (\widetilde{\U}_{i-1}^{(t+1)})^\top \times_{\widecheck{\Omega}_i} (\widetilde{\U}_{i}^{(t)})^\top \times \cdots \times_{\Omega_{d}} (\widetilde{\U}_{d}^{(t)})^\top )  \big),
				\end{split}
				\end{equation}
				where for any matrix $\D$, $\SVD_r(\D)$ computes the subspace composed of the leading $r$ left singular vectors of $\D$.
			\end{enumerate}
			\State Let $\widehat{\U}_i = \widetilde{\U}_i^{(t_{\max})}$ for $i = 1, \ldots, m$ and compute 
			\begin{equation*}
			\widehat{\bcT} = \widetilde{\bcT}\times_{\Omega_1} P_{\widehat{\U}_1}\times \cdots \times_{\Omega_m} P_{\widehat{\U}_m}.
			\end{equation*}
		\end{algorithmic}
	\end{algorithm}
	
	\subsection{Our Contributions} \label{sec: contributions}
	Under the perturbation model \eqref{eq: perturbation model} and the HOOI Algorithm \ref{alg: HOOI}, we make the following major contributions to the tensor perturbation theory of low-rank tensor decomposition based on HOOI.
	\begin{itemize}
		\item We provide the first sharp blockwise perturbation bounds of HOOI for tensors with guarantees for both the estimated mode-$k$ singular subspace and the tensor reconstruction induced by low-rank approximation of the noise corrupted tensor. 
		Specifically, in Theorem \ref{th:HOOI general}, we show that under suitable initialization $\{\widetilde{\U}^{(0)}_i \}_{i=1}^m$ in Algorithm \ref{alg: HOOI}, the upper bound of $\max_i \left\|\sin\Theta\left(\widehat{\U}_i, \U_i\right)\right\|_q$ converges linearly with respect to the iteration number to a fixed quantity characterized by $\bcZ$ and the signal strength. 
		In addition, a practically useful tensor reconstruction error bound in Hilbert-Schmidt norm is provided. Surprisingly, we found the upper bound of tensor reconstruction is free of the ``condition number" of the signal tensor and can be expressed by a unified simple quantity $\xi$ characterized only by the noise tensor $\bcZ$ and the underlying multilinear rank of $\bcT$. $\xi$ is closely related to the Gaussian width \citep{gordon1988milman}, a common measure for the complexity of a given set, and its formal definition and more explanations will be given in Sections \ref{sec: error terms of Z} and \ref{sec: applications}. 
		
		\item In addition, we also generalize the main results to the case that only a subset of modes of $\bcT$ have low-rank structure. 
		
		\item Furthermore, we provide a deterministic minimax lower bound for the tensor reconstruction error under perturbation model \eqref{eq: perturbation model} in Theorem \ref{th: pertur lower bound}. The lower bound matches the perturbation upper bound in Theorem \ref{th:HOOI general} when the tensor order $d$ is fixed, which demonstrates the optimality of HOOI for tensor reconstruction.
		\item In addition, by combining Theorems \ref{th: pertur lower bound} and \ref{th:HOOI general}, we prove that the tensor reconstruction error rate of HOOI with only one iteration is also optimal and further iterations improve the coefficient in front of the error rate $\xi$. This suggests that in some applications where running full HOOI is too expensive and prohibitive compared to truncated HOSVD \citep{de2000multilinear} or sequentially truncated HOSVD \citep{vannieuwenhoven2012new}, we can just run HOOI for one iteration to obtain an optimal (up to constant) reconstruction. Details are provided in Remark \ref{rem: one-step optimality of HOOI} and numerical comparison is given in Section \ref{sec: numerical study: comparison}.
		\item In addition, we apply the new perturbation bounds of HOOI in two modern applications, tensor denoising and tensor co-clustering, from machine learning and statistics. Based on our perturbation results, we can easily recover the results of tensor denoising in current literature with a much shorter proof and provide the first guarantee of HOOI on tensor cocluster recovery that improves the state-of-the-art results.
		\item Finally, we perform extensive numerical studies to support our perturbation bounds and do a comparison with other existing low-rank tensor decomposition algorithms. 
	\end{itemize}
	
	\subsection{Related Literature} \label{sec: literature review}
	In this section, we give a brief overview of the literature on Tucker decomposition of tensors, matrix/tensor perturbation theory related to this article. Tensor decomposition has become one of the most important topics in the literature on tensors \citep{kolda2009tensor, sidiropoulos2017tensor, cichocki2015tensor}. An analogy of matrix SVD for Tucker decomposition of tensors, which is today commonly called high-order singular value decomposition (HOSVD), was first proposed in \cite{tucker1966some} and then popularized by \cite{de2000multilinear}. However, unlike the truncated matrix SVD, truncated HOSVD (T-HOSVD) can provide a reasonable but not necessarily optimal low-multilinear-rank approximation for a given tensor. It has been shown that computing the best low multilinear rank approximation is NP hard in general \citep{hillar2013most}. On the other hand, various computationally efficient algorithms have been proposed to obtain better approximations than HOSVD \citep{de2000best,elden2009newton,ishteva2011best,ishteva2009differential,savas2010quasi}. Among them, the higher-order orthogonal iteration (HOOI) in \cite{de2000best} has become a popular choice in literature. HOOI has been applied to problems including but not limited to tensor PCA/tensor SVD \citep{richard2014statistical,zhang2018tensor}, tensor completion \citep{xia2017statistically}, tensor regression \citep{zhang2020islet}, hypergraph community recovery \citep{ke2019community,jing2020community}, independent component analysis \citep{de2004dimensionality}, tensor clustering \citep{luo2020tensor}. Many variants of HOOI such as sparse high-order singular value decomposition (STAT-SVD) for tensors \citep{zhang2019optimal}, regularized HOOI \citep{ke2019community,jing2020community}, generalized higher-order orthogonal iteration (gHOI) \citep{liu2014generalized} have been proposed. Nowadays, HOOI has become a prevalent choice to compute the low-rank tensor approximation in many applications and been coded in common tensor software such as Matlab ``Tensor Toolbox" \citep{bader2012matlab}, ``Tensorlab" \citep{sorber2014tensorlab} and R ``rTensor" package \citep{li2018rtensor}. Moreover, it has been regarded as the gold standard guideline for comparison when developing even faster randomized or memory-efficient algorithms for low-rank tensor approximation \citep{sun2020low,malik2018low,tsourakakis2010mach,kolda2008scalable}. 
	
	Perturbation theory is a long-existing field in mathematics. In particular, the perturbation theory on matrices has attracted much attention. The original work in matrix dated back to Weyl, Davis-Kahan and Wedin \citep{davis1970rotation,wedin1972perturbation,weyl1912asymptotische,stewart1998perturbation} (see \cite{bhatia2013matrix,stewart1990matrix} for an overview of classical perturbation results and historical development) and recently it has been further developed in \cite{yu2014useful,cai2018rate,cape2019two}. In addition, various generalizations and extensions have been made in different settings including random perturbation \citep{vu2011singular,shabalin2013reconstruction,o2018random,wang2015singular,benaych2011eigenvalues,abbe2020entrywise,koltchinskii2016perturbation,benaych2012singular,cape2019signal,chen2018asymmetry}, structured perturbation \citep{fan2018eigenvector,stewart2006perturbation} and many others \citep{eldridge2017unperturbed}. Also the perturbation theory for matrices has been widely applied to a number of applications such as community detection \citep{rohe2011spectral,chaudhuri2012spectral,chin2015stochastic,sussman2012consistent,cape2019two}, covariance matrix estimation \citep{fan2018eigenvector,cape2019two}, matrix denoising \citep{cai2018rate}, matrix completion \citep{cai2016structured}, etc. 
	
	Studying perturbation theory for low-rank tensor approximation is much harder than matrices and there are only a few attempts in literature. 
	Perturbation results for the best low multilinear rank approximation have been developed in \cite{de2004first,elden2011perturbation}. These results are promising but less practically useful due to the computational hardness of the best multilinear rank approximation \citep{hillar2013most}. Hence more attention are given to the perturbation results for polynomial-time low-rank tensor approximation algorithms. A line of work has been done when the perturbation is random \citep{richard2014statistical,zhang2018tensor,hopkins2015tensor,liu2017characterizing,zhang2019optimal,zheng2015interpolating,xia2019sup} and various perturbation results for robust tensor decomposition with sparse noise have been studied in \cite{goldfarb2014robust,gu2014robust,shah2015sparse,anandkumar2016tensor}. Perturbation bounds for orthogonally decomposable tensors have been studied in \cite{mu2017greedy,auddy2020perturbation}. In addition, \cite{anandkumar2014tensor,anandkumar2014guaranteed} provided perturbation guarantees for power iteration algorithm for symmetric orthogonal and non-orthogonal CP low-rank decomposition. However, we are not aware of any perturbation result for polynomial time algorithms under partial symmetric multilinear low-rank setting. In this paper, we make the first attempt in this direction and provide the first perturbation bounds of HOOI for tensors with guarantees for both singular subspaces and tensor reconstruction in the general setting. It is worth mentioning that the reconstruction error bound of this paper is often significantly better than the simple estimator and truncated HOSVD without power iteration. This is fundamentally different from the perturbation results in \cite{anandkumar2014tensor}. See further discussions in Remarks \ref{rem: compare of HOSVD and HOOI} and \ref{rem: compare with Anandkumar}. 

	We end this section by remarking that in most situations there is a trade-off about the quality of low-rank tensor decomposition and computational cost of the algorithm. For example, computing truncated HOSVD and sequentially truncated HOSVD (ST-HOSVD) \citep{vannieuwenhoven2012new} may be much faster than iterative algorithms such as HOOI, (quasi-)Newton-Grassmann method \citep{elden2009newton,savas2010quasi}, geometric Newton method \citep{ishteva2009differential} and Riemannian trust region scheme \citep{ishteva2011best} in the large scale settings. On the other hand, these iterative algorithms achieve higher accuracy. In the perturbation model \eqref{eq: perturbation model}, we show HOOI could achieve optimal tensor reconstruction error, which is not true for HOSVD and ST-HOSVD in general. 
	
	\subsection{Organization of the Paper} \label{sec: organization}
	The remainder of the article is organized as follows. In Section \ref{sec: notation}, after a brief introduction of notation and preliminaries, we define various blockwise errors of $\bcZ$ as the key quantities in our perturbation bounds. We illustrate our main perturbation theorem in asymmetric order-$3$ case in Section \ref{sec: main results d = 3 case} and at the end of the same section, we provide the deterministic lower bound for tensor reconstruction. In Section \ref{sec:main}, we provide the perturbation bounds of HOOI applying on a corrupted general partial symmetric order-$d$ tensor. 
	In Section \ref{sec: HOOI in partial low rank}, we discuss the tensor perturbation bounds when the target tensor has the low-rank structure only along a subset of modes. In Section \ref{sec: applications}, we apply our perturbation bounds to two applications, tensor denoising and tensor co-clustering. In Section \ref{sec: numerical study}, we corroborate our theoretical results by extensive numerical studies. Conclusion and discussions are provided in Section \ref{sec: conclusion and discussion}. Due to space constraints, all technical proofs are postponed to the Appendix. 
	
	\section{Notation and Preliminaries} \label{sec: notation}
	
	The following notation will be used throughout this article. For any non-negative integer $n$, let $[n] = \{1, \ldots, n\}$. Lowercase letters (e.g., $a, b$), lowercase boldface letters (e.g., $\u, \v$), uppercase boldface letters (e.g., $\U, \V$), and boldface calligraphic letters (e.g., $\bcT, \bcZ$) are used to denote scalars, vectors, matrices, and order-3-or-higher tensors, respectively. For any two series of numbers, say $\{a_n\}$ and $\{b_n\}$, denote $a \asymp b$ or $a = O(b)$ if $ca_n \leq b_n\leq Ca_n$ or $a_n \leq Cb_n$ for some uniform constants $c, C>0$. For any matrix $\D \in \mathbb{R}^{m\times n}$ with singular value decomposition $\sum_{i=1}^{m \land n} \sigma_i(\D)\u_i \v_i^\top$ ($\sigma_1(\D) \geq \cdots \geq \sigma_{m \vee n}(\D)$), let $\D_{\max(r)}= \sum_{i=1}^{r} \sigma_i(\D)\u_i \v_i^\top$ be the leading rank-$r$ SVD approximation of $\D$ and $\D_{\max(-r)} = \sum_{i=r+1}^{m \land n} \sigma_i(\D)\u_i \v_i^\top$ be its complement. We also denote $\SVD_r(\D) := [\u_1 ~ \cdots \u_r]$ as the subspace composed of the leading $r$ left singular vectors of $\D$. The Schatten-$q$ norm of matrix $\D$ for $q \geq 1$ is defined as $\|\D\|_q := \left(\sum_{i=1}^{m \land n} \sigma^q_i(\D)\right)^{1/q}$. Frobenius norm $\|\cdot\|_F$ and spectral norm $\|\cdot\|$ of a matrix are special cases of Schatten-$q$ norm with $q = 2$ and $q = \infty$. In addition, $\I_r$ represents the $r$-by-$r$ identity matrix. Let $\mathbb{O}_{p, r} = \{\U: \U^\top \U=\I_r\}$ be the set of all $p$-by-$r$ matrices with orthonormal columns. For any $\U\in \mathbb{O}_{p, r}$, $P_{\U} = \U\U^\top$ represents the projection matrix onto the column span of $\U$; we also use $\U_\perp\in \mathbb{O}_{p, p-r}$ to represent the orthonormal complement of $\U$. We use bracket subscripts to denote sub-matrices. For example, $\D_{[i_1,i_2]}$ is the entry of $\D$ on the $i_1$-th row and $i_2$-th column; $\D_{[(r+1):m, :]}$ contains the $(r+1)$-th to the $m$-th rows of $\D$. For any matrices $\U\in \mathbb{R}^{p_1\times p_2}$ and $\V\in \mathbb{R}^{m_1\times m_2}$, let
	$$\U\otimes \V = \begin{bmatrix}
	\U_{[1,1]}\cdot \V & \cdots & \U_{[1, p_2]}\cdot \V\\
	\vdots & & \vdots\\
	\U_{[p_1,1]}\cdot \V & \cdots & \U_{[p_1, p_2]}\cdot \V\\
	\end{bmatrix}\in \mathbb{R}^{(p_1m_1)\times(p_2m_2)}$$ 
	be the Kronecker product of $\U$ and $\V$. 
	
	For any order-$d$ tensor $\bcT \in \mathbb{R}^{p_1\times \cdots\times p_d}$, let $\mathcal{M}_k(\cdot)$ be the matricization operation that unfolds or flattens the order-$d$ tensor $\bcT \in\mathbb{R}^{p_1\times \cdots \times p_d}$ along mode $k$ into the matrix $\mathcal{M}_k(\bcT)\in \mathbb{R}^{p_k\times p_{-k} }$ and here $p_{-k} := \prod_{j\neq k}p_j$. Specifically, the mode-$k$ matricization of $\bcT$ is formally defined as 
	\begin{equation} \label{eq: tensor-matricization}
	\mathcal{M}_k(\bcT) \in \mathbb{R}^{p_k\times p_{-k}}, \quad \left(\mathcal{M}_k(\bcT)\right)_{\left[i_k, j\right]} =\bcT_{[i_1,\ldots, i_d]}, \quad j = 1 + \sum_{\substack{l=1\\l\neq k}}^d\left\{(i_l-1)\prod_{\substack{m=1\\m\neq k}}^{l-1}p_m\right\}
	\end{equation}
	for any $1\leq i_l \leq p_l, l=1,\ldots, d$. Also see \cite[Section 2.4]{kolda2009tensor} for more discussion on tensor matricizations. Given two tensors $\bcT_1, \bcT_2 \in \mathbb{R}^{p_1\times \cdots\times p_d}$, define their inner product as $\langle \bcT_1, \bcT_2 \rangle = \sum_{i_1,\ldots, i_d} \langle \bcT_{1[i_1,\ldots, i_d]}, \bcT_{2[i_1,\ldots, i_d]} \rangle $. The Hilbert-Schmidt norm of $\bcT$ is defined as
	$\|\bcT\|_{\tHS} = \left( \langle \bcT, \bcT \rangle \right)^{1/2}.$ The multilinear rank of a tensor $\bcT$, $\rank(\bcT)$, is defined as a $d$-tuple $(r_1, \ldots, r_d)$, where $r_k = \text{rank}(\mathcal{M}_k(\bcT))$ is the mode-$k$ rank. For any multilinear rank-$(r_1, \ldots, r_d)$ tensor $\bcT$, it has Tucker decomposition \citep{tucker1966some}:
	\begin{equation}\label{eq:tucker-decomposition}
	\bcT = \llbracket \bcS; \U_1, \ldots, \U_d\rrbracket := \bcS\times_1 \U_1 \times \cdots \times_d \U_d,
	\end{equation}
	where $\bcS \in \bbR^{r_1 \times \cdots \times r_d}$ is the core tensor and $\U_k = \SVD_{r_k}(\cM_k(\bcT))$ is the mode-$k$ singular subspace. Here, the mode-$k$ product of $\bcT \in \mathbb{R}^{p_1 \times \cdots \times p_d}$ with a matrix $\U\in \mathbb{R}^{p_k\times r_k}$ is denoted by $\bcT \times_k \U^\top$ and is of size $p_1 \times \cdots \times p_{k-1}\times r_k \times p_{k+1}\times \cdots \times p_d$, such that 
	\begin{equation}\label{eq: tensor-matrix product}
		(\bcT \times_k \U^\top)_{[i_1, \ldots, i_{k-1}, j, i_{k+1}, \ldots, i_d]} = \sum_{i_k=1}^{p_k} \bcT_{[i_1, i_2, \ldots, i_d]} \U_{[i_k, j]}.
	\end{equation}
	Given $S = \{i_1, \ldots, i_d \}$, it is convenient to denote the product of $\bcT$ along the modes indexed by $S$ with the same matrix $\U$ and with different matrices $\{\U_i\}$ respectively as
	\begin{equation} \label{eq: tensor matrix prod short notation}
	\begin{split}
	\bcT \times_{ S}\U := \bcT \times_{i_1} \U \times \cdots \times_{i_d} \U,\quad \bcT \times_{i \in S}\U_i := \bcT \times_{i_1} \U_{i_1} \times \cdots \times_{i_d} \U_{i_d}.
	\end{split}
	\end{equation} 
	Given $\bcT = \llbracket \bcS; \U_1, \ldots, \U_d\rrbracket$, the following relationship of tensor matricization is used frequently in the proof:
	\begin{equation}\label{eq: matricization relationship}
	\mathcal{M}_k \left( \bcS\times_1 \U_1 \times_2 \cdots \times_d \U_d \right) = \U_k \mathcal{M}_k(\bcS) (\U_d^\top \otimes \cdots \otimes \U_{k+1}^\top \otimes \U_{k-1}^\top \otimes \cdots \otimes \U_{1}^\top ).
	\end{equation}
	We refer the readers to \cite{kolda2001orthogonal, kolda2009tensor} for a more comprehensive survey on tensor algebra. 
	
	Finally, we use $\sin \Theta$ distance to measure the difference between two $p$-by-$r$ column orthogonal matrices $\widehat{\U}$ and $\U$. Suppose the singular values of $\widehat{\U}^\top \U$ are $\sigma_1 \geq \sigma_2 \geq \ldots \geq \sigma_r \geq 0$. Then $\Theta(\widehat{\U}, \U)$ is defined as $$\Theta(\widehat{\U}, \U) = {\rm diag}\left( \cos^{-1} (\sigma_1), \cos^{-1} (\sigma_2), \ldots, \cos^{-1} (\sigma_r) \right).$$ Common properties of $\sin \Theta$ distance can be found in  \cite[Lemma 1]{cai2018rate} and \cite[Lemma 6]{luo2020schatten}.
	
	\subsection{Blockwise Errors of $\bcZ$} \label{sec: error terms of Z}
	In this subsection, we introduce key quantities appearing in the perturbation bounds that characterize the blockwise errors of $\bcZ$. For simplicity, we consider order-3 tensors and $\Omega_1 = \{1\}, \Omega_2 = \{2\}, \Omega_3 = \{3\}$ perturbation setting for illustration. 
	
	Define the blockwise errors of $\bcZ$ that characterize the tensor perturbation:
	\begin{equation} \label{eq: blockwise error d = 3}
	\begin{split}
	& \tau_1 = \max_{k\in [3]} \tau_{1k}, \quad \tau_{1k} = \left\|  \left(\mathcal{M}_k(\bcZ \times_{k+1}\U_{k+1}^\top \times_{k+2}\U_{k+2}^\top ) \right)_{\max(r_k)}\right\|_q, \quad k=1,2,3;\\
	& \tau_2 = \max_{k\in [3]}\Big\{\max_{\substack{\V\in \mathbb{R}^{(p_{k+1}-r_{k+1})\times r_{k+1}}\\ \|\V\|_q\leq 1}}\left\|\left(\mathcal{M}_k(\bcZ \times_{k+1}(\U_{k+1\perp}\V)^\top \times_{k+2}\U_{k+2}^\top ) \right)_{\max(r_k)}\right\|_q, \\
	&\quad\quad\quad\quad  \max_{\substack{\V\in \mathbb{R}^{(p_{k+2}-r_{k+2})\times r_{k+2}}\\ \|\V\|_q\leq 1}}\left\|\left(\mathcal{M}_k(\bcZ \times_{k+1}\U_{k+1}^\top \times_{k+2}(\U_{k+2\perp}\V)^\top ) \right)_{\max(r_k)}\right\|_q \Big\};\\
	& \tau_3 = \max_{k\in [3]} \max_{\substack{\V\in \mathbb{R}^{(p_{k+1}-r_{k+1})\times r_{k+1}} \\\V'\in \mathbb{R}^{(p_{k+2}-r_{k+2})\times r_{k+2}}\\ \|\V\|_q\leq 1, \|\V'\|_q\leq 1}} \left\|\left(\mathcal{M}_k(\bcZ \times_{k+1}(\U_{k+1\perp}\V)^\top \times_{k+2}(\U_{k+2\perp}\V')^\top ) \right)_{\max(r_k)}\right\|_q.
	\end{split}
	\end{equation}
	Here all mode indices $(\cdot)_k$ of an order-3 tensor are in the sense of modulo-3, e.g., $r_1=r_4$, $p_2 = p_5$.
	
\begin{figure}
		\centering
		\subfigure[$\tau_{11}$]{ \includegraphics[width = 0.3\textwidth, height=1.5in]{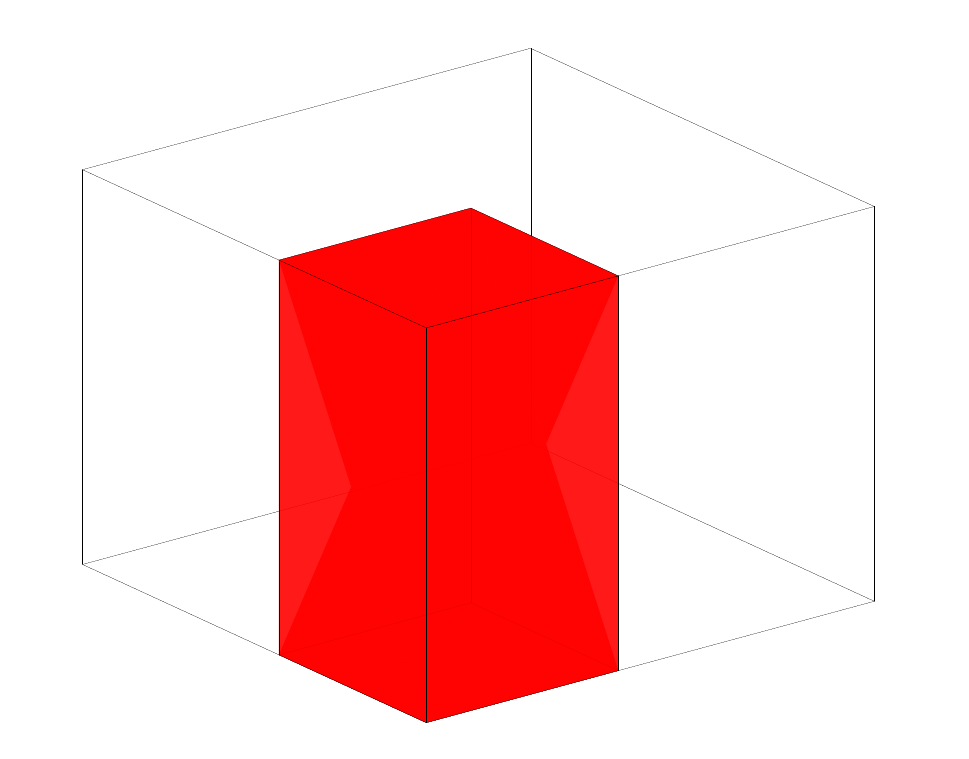} }
		\subfigure[$\tau_{12}$]{ \includegraphics[width = 0.3\textwidth, height=1.5in]{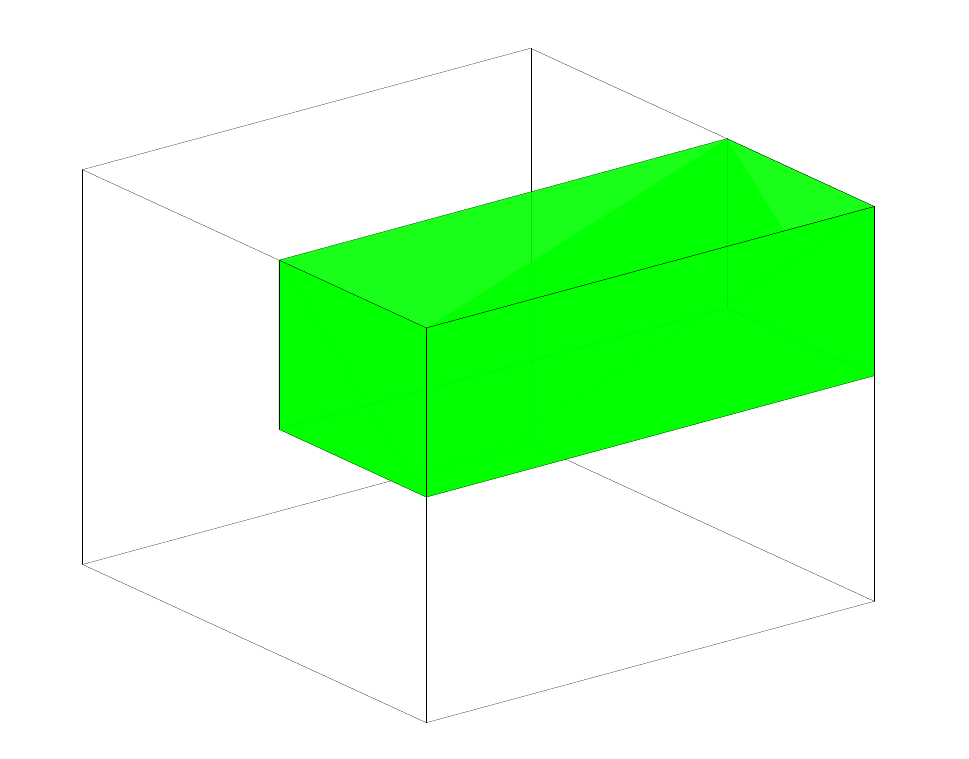} }
		\subfigure[$\tau_{13}$]{ \includegraphics[width = 0.3\textwidth, height=1.5in]{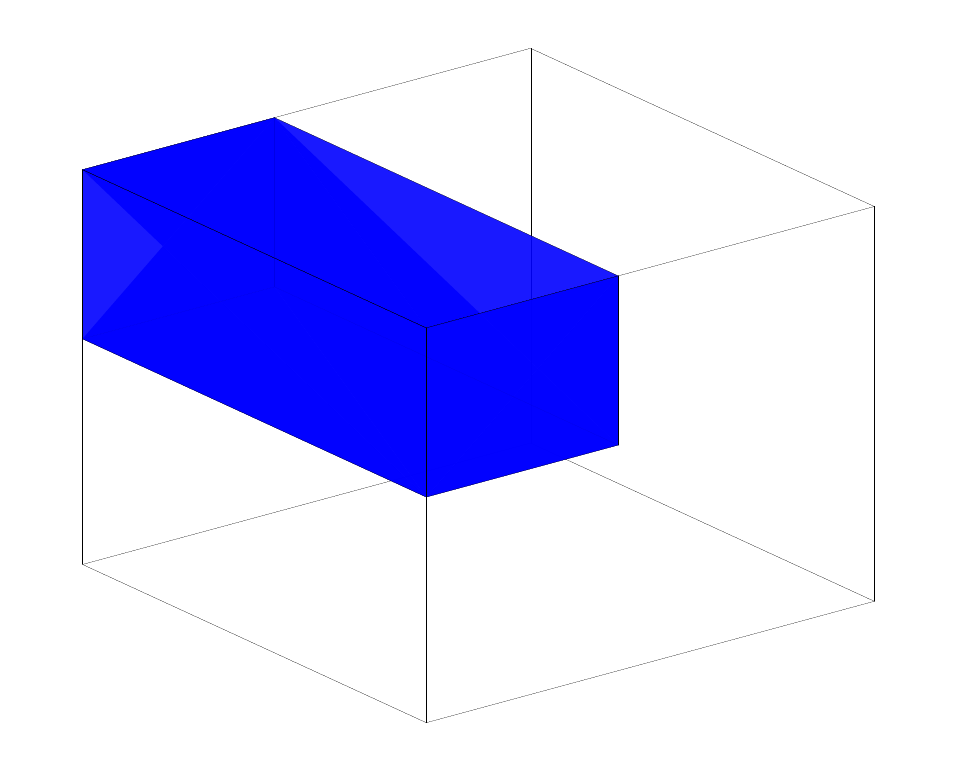} }
		\caption{Illustration of $\tau_{11},\tau_{12},\tau_{13}$. Here, we assume $\U_k^\top = [\I_{r_k} ~ \boldsymbol{0}_{r_k \times (p_k-r_k)}]$, $k=1,2,3$, for a better visualization. The red, green, blue blocks represent the corresponding blockwise errors $\tau_{11}, \tau_{12}, \tau_{13}$ in $\bcZ$.}\label{fig: illustration of tau}
	\end{figure}

	$\tau_1, \tau_2, \tau_3$ represent the maximum of blockwise errors of $\bcZ$ in the projection spaces expanded by $\U_1, \U_2, \U_3$ and their complements. For example, in Figure \ref{fig: illustration of tau} we illustrate the blockwise errors characterized by $\tau_{11}, \tau_{12}, \tau_{13}$. $\tau_2, \tau_3$ characterize blockwise errors of $\bcZ$ in a similar fashion but with more complicated projections.
	
	These blockwise errors of $\bcZ$ are in fact a generalization of error terms in matrix perturbation theory. In the matrix setting $\widetilde{\T} = \T + \Z$, let $\widehat{\U}, \widehat{\V}$ and $\U, \V$ be leading left and right singular vectors of $\widetilde{\T}$ and $\T$, respectively. Then by Wedin's perturbation theory, the upper bounds of $\left\|\sin \Theta(\widehat{\U}, \U)\right\|$ and $\left\|\sin \Theta(\widehat{\V}, \V)  \right\|$ involve $\|\Z \widehat{\V}\|$ and $\|\widehat{\U}^\top \Z \|$, which are also blockwise errors of $\Z$.
	
	Next, we introduce a simple quantity $\xi$ that characterizes the error bound for tensor reconstruction. In this order-3 asymmetric setting, $\xi$ is defined as
	\begin{equation} \label{eq: xi}
	\xi := \underset{\|\bcY \|_{\tHS} \leq 1,  \rank(\bcY) \leq (r_1, r_2, r_3)}\sup \langle \bcZ, \bcY \rangle.
	\end{equation}
	In the following Lemma \ref{lem: character of xi}, we give another equivalent characterization of $\xi$.
	\begin{Lemma}[Equivalent Characterizations of $\xi$]\label{lem: character of xi}
		\begin{equation*}
		\begin{split}
		\xi &:= \underset{\|\bcY \|_{\tHS} \leq 1,  \rank(\bcY) \leq (r_1, r_2, r_3)}\sup \langle \bcZ, \bcY \rangle \\
		& = \sup_{\U_i \in \mathbb{O}_{p_i, r_i}, 1 \leq i \leq 3} \| \bcZ \times_1 \U_1^\top \times \U_2^\top \times_d \U_3^\top \|_{\tHS}.
		\end{split}
		\end{equation*}
	\end{Lemma}
	By Lemma \ref{lem: character of xi}, $\xi$ measures the maximum magnitude of the projection of $\bcZ$ onto low rank subspaces in Hilbert-Schmidt norm. Another interpretation of $\xi$ is from Gaussian width \citep{gordon1988milman}, which will be discussed in Section \ref{sec: applications}.
	
	Although the exact values of $\tau_j$ and $\xi$ may be hard to compute in general, it is often practical to provide probabilistic bounds when we impose distributional assumptions on $\bcZ$. For example, if $\bcZ$ is a random tensor with i.i.d. standard normal entries and consider $r_1 = r_2 = r_3 = r, p_1 = p_2 = p_3 = p$, $q = \infty$, then by random matrix theory \citep{vershynin2010introduction}, we can show that with high probability $\tau_{1i} \asymp (\sqrt{p} + r)$, $\tau_2 \asymp  \sqrt{p r },\tau_3 \asymp  \sqrt{p r}$ and $\xi \asymp (r^{\frac{3}{2}} + \sqrt{pr} )$. More details about the upper bounds for $\tau_j$ and $\xi$ can be found in tensor denoising and tensor co-clustering applications in Section \ref{sec: applications}.
	
	\section{Illustration of Perturbation Bounds for $d = 3$ Asymmetric Case} \label{sec: main results d = 3 case}
	In this section, we present our main results in the $d = 3$ asymmetric case to better illustrate the main ideas in this paper. The specialized HOOI algorithm for the $d = 3$ asymmetric case is given in Algorithm \ref{alg: HOOI in d = 3} and its guarantee is provided in Theorem \ref{th:HOOI d=3}.
	
	\begin{algorithm}[!h] \caption{Higher-Order Orthogonal Iteration for Tensor Decomposition ($d = 3$)} \label{alg: HOOI in d = 3}
		\textbf{Input:} $\widetilde{\bcT}\in \bbR^{p_{1} \times p_2 \times p_{3}}$, initialization $\{\widetilde{\U}_i^{(0)}\}_{i=1}^3$ with $\widetilde{\U}_i^{(0)} \in \mathbb{O}_{p_i,r_i}$, maximum number of iterations $t_{\max}$.\\
		\textbf{Output:} $ \{\widehat{\U}_i\}_{i=1}^3, \widehat{\bcT}$. 
		\begin{algorithmic}[1]
			\State For $t=1, \ldots, t_{\max}$, do
			\begin{equation*}
			\begin{split}
			&\widetilde{\U}_1^{(t+1)} = \SVD_{r_1} \left(\cM_{1}\left( \widetilde{\bcT} \times_2 ( \widetilde{\U}_2^{(t)})^\top \times_3 (\widetilde{\U}_{3}^{(t)})^\top  \right) \right),\\
			&\widetilde{\U}_2^{(t+1)} = \SVD_{r_2} \left(\cM_{2}\left( \widetilde{\bcT} \times_1 ( \widetilde{\U}_1^{(t+1)})^\top \times_3 (\widetilde{\U}_{3}^{(t)})^\top  \right) \right),\\
			&\widetilde{\U}_3^{(t+1)} = \SVD_{r_3} \left(\cM_{3}\left( \widetilde{\bcT} \times_1 ( \widetilde{\U}_1^{(t+1)})^\top \times_2 (\widetilde{\U}_{2}^{(t+1)})^\top  \right) \right).
			\end{split}
			\end{equation*}
			\State Let $\widehat{\U}_i = \widetilde{\U}_i^{(t_{\max})}$ for $i = 1, 2, 3$ and compute 
			\begin{equation*}
			\widehat{\bcT} = \widetilde{\bcT}\times_{1} P_{\widehat{\U}_1}\times_2 P_{\widehat{\U}_2} \times_{3} P_{\widehat{\U}_3}.
			\end{equation*}
		\end{algorithmic}
	\end{algorithm}
	
	\begin{Theorem}[Tensor Perturbation Bounds for HOOI (d = 3)]\label{th:HOOI d=3}
		Consider the perturbation model \eqref{eq: perturbation model} with $\widetilde{\bcT}, \bcT, \bcZ \in \mathbb{R}^{p_1\times p_2\times p_3}$. Suppose $q\geq1$. Define the blockwise errors as in \eqref{eq: blockwise error d = 3} \eqref{eq: xi} and denote the initialization errors of $\{ \widetilde{\U}_k^{(0)} \}_{k=1}^3$ as $\tilde{e}_0 := \max_{k=1,2,3}\|\widetilde{\U}^{(0)\top}_{k\perp} \U_k\|$, $e_0 := \max_{k=1,2,3}\|\widetilde{\U}^{(0)\top}_{k\perp} \U_k\|_q$. Assume the initialization error and the signal strength satisfy 
		\begin{equation} \label{ineq: signal strength condition d = 3}
		\tilde{e}_0 \leq \sqrt{2}/2\, \text{ and } \,\lambda \geq (20 + 28 \sqrt{2})\xi.
		\end{equation} 
		Let $\widetilde{\bcT}^{(t)}:= \widetilde{\bcT} \times_{1} P_{\widetilde{\U}^{(t)}_1} \times_2 P_{\widetilde{\U}^{(t)}_2} \times_{3} P_{\widetilde{\U}^{(t)}_3}$ be the estimator of $\bcT$ after $t$ steps in Algorithm \ref{alg: HOOI}. Then with inputs $\widetilde{\bcT}$, $\{ \widetilde{\U}_i^{(0)} \}_{i=1}^3$, the mode-$k$ singular subspace updates in Algorithm \ref{alg: HOOI in d = 3} after $t$ iterations satisfy
		\begin{equation} \label{ineq: max singular bound d = 3}
		\max_{k = 1,2,3} \left\|\sin\Theta\left(\widetilde{\U}^{(t)}_k, \U_k\right)\right\|_q \leq \frac{8\tau_1}{\lambda} + \frac{e_0}{2^{t}},
		\end{equation} and the $t$-step tensor estimation satisfies
		\begin{equation} \label{ineq: t step tensor reconstruction}
		\left\|\widetilde{\bcT}^{(t)} - \bcT\right\|_{\tHS} \leq \left(  1+ \frac{6}{1 - \left( 8 \frac{\tau_1}{\lambda} + \frac{\tilde{e}_0}{2^{t-1}}   \right)^2} \right)\xi.
		\end{equation}
		
		Moreover, when $t_{\max} \geq C \log(e_0 \lambda /\tau_{1})\vee 1$ for some $C > 0$, the outputs of the estimated mode-$k$ singular subspaces of Algorithm \ref{alg: HOOI in d = 3} satisfy
		\begin{equation} \label{ineq: d=3-max-upper bound of singular space}
		\max_{k = 1,2,3} \left\|\sin\Theta\left(\widehat{\U}_k, \U_k\right)\right\|_q \leq  \frac{9\tau_1}{\lambda},
		\end{equation}
		\begin{equation}
		\begin{split} \label{ineq: singular space upp bnd d = 3}
		\left\|\sin\Theta\left(\widehat{\U}_k, \U_k\right)\right\|_q = \left\| \widehat{\U}_{k\perp}^{\top} \U_k \right\|_q \leq 4 \left(\frac{\tau_{1k}}{\lambda} +  \frac{18\tau_1 \tau_2}{\lambda^2} + \frac{81\tau_1^2 \tau_3}{\lambda^3} \right), \quad k=1,2,3,
		\end{split}
		\end{equation}
		and the output of tensor reconstruction $\widehat{\bcT}$ satisfies
		\begin{equation}\label{ineq:reconstruction-error-bound}
		\begin{split}
		\left\|\widehat{\bcT} - \bcT\right\|_{\tHS} & \leq \left\|\llbracket \bcZ; \widehat{\U}_1^\top, \widehat{\U}_2^\top, \widehat{\U}_3^\top \rrbracket\right\|_{\tHS}+ \sum_{k=1}^3 \left\|\widehat{\U}_{k\perp}^\top \mathcal{M}_k(\bcT)\right\|_F \leq 13 \xi.
		\end{split}
		\end{equation}
	\end{Theorem}
	
	\begin{Remark}[Noise Tolerance and Least Singular Value $\lambda$] \label{rem: signal strength and noise tolerance}
		Our theory relies on a lower bound assumption of the least singular value: $\lambda \geq C\xi$, which is in the same vein as the classical matrix perturbation theory \citep{davis1970rotation,wedin1972perturbation}. 
		Moreover, in the existing results on perturbation analysis for Canonical-Polyadic (CP) decomposition, e.g., Theorem 5.1 of \cite{anandkumar2014tensor} and Theorem 1 of \cite{anandkumar2014guaranteed}, one assumes $\lambda \geq C\|\bcZ\|$. 
		Since $\|\bcZ\| = \sup_{\bcY: \|\bcY\|_{\tHS} \leq 1, \rank(\bcY) = (1,\cdots,1)} \langle \bcZ, \bcY \rangle$, $\xi$ defined in \eqref{eq: xi} can be seen as a counterpart of $\|\bcZ\|$ in Tucker decomposition. 
	\end{Remark}
	
\begin{Remark}[Initialization] \label{rem: initialization}
	In Theorem \ref{th:HOOI d=3}, we assume the initialization $\{\widetilde{\U}_i^{(0)} \}_{i=1}^3$ is warm in the sense that the maximum error $\tilde{e}_0$ is upper bounded by a constant. The constant $\sqrt{2}/2$ in this upper bound is chosen for convenience and can be replaced by any fixed constant less than $1$. 
	Our perturbation bound applies to HOOI with any initialization as long as this condition holds, although the original HOOI algorithm was proposed with the initialization scheme named HOSVD, i.e., $\hat\U_k^{(0)} = \SVD_{r_k}\left(\mathcal{M}_k\left(\widetilde{\bcT}\right)\right)$. 
	Next, we briefly discuss two specific initialization schemes: HOSVD for tensor PCA/SVD \citep{richard2014statistical,zhang2018tensor} and diagonal-deletion SVD for tensor completion \citep{xia2017statistically}. For convenience of presentation, we focus on the setting $p_1 = p_2 = p_3 = p$.
\begin{itemize}
	\item (Tensor Denoising) Suppose we observe a tensor $\widetilde{\bcT} =\bcT + \bcZ \in \bbR^{p \times p \times p}$ and aim to recover $\bcT$ from $\widetilde{\bcT}$. To this end, we can apply HOOI by inputting $\widetilde{\bcT}$. When $\bcZ$ has i.i.d. $N(0,\sigma^2)$ entries, Theorem 1 in \cite{zhang2018tensor} showed if one initializes by HOSVD, as long as $\lambda \geq C p^{3/4}$, the initialization condition $\|\sin \Theta (\widetilde{\U}_k^{(0)}, \U_k)\| \leq \frac{\sqrt{2}}{2}$ holds with high probability. \cite{zhang2018tensor} also showed the signal strength requirement $\lambda \geq C p^{3/4}$ is essential, which means HOSVD is a proper initialization in the tensor denoising model.
	\item (Tensor Completion) Suppose we observe a set of entries, selected uniformly at random and indexed by $\Omega$, from a noisy tensor $\bar{\bcT} = \bcT + \bcZ$. Denote $\bar{\bcT}_\Omega$ as 
	$$(\bar{\bcT}_\Omega)_{[i_1,\ldots, i_d]} = \left\{\begin{array}{c c}
	\bar{\bcT}_{[i_1,\dots, i_d]}, & (i_1,\cdots,i_d) \in \Omega;\\
	0 & \text{otherwise}.
	\end{array}\right.$$
	Suppose $\bcZ$ has i.i.d. $N(0,\sigma^2)$ entries. Then, it is easy to check that $\bar{\bcT}_\Omega/\rho$ is an unbiased estimator of $\bcT$, where $\rho$ is the sampling ratio. \cite{xia2017statistically} proposed to apply HOOI on $\widetilde{\bcT}:=\bar{\bcT}_\Omega/\rho$ to estimate $\bcT$. They proposed to set $\widetilde{\U}_k^{(0)}$ as the leading $r_k$ singular vectors of $\mathcal{M}_k(\widetilde{\bcT}) \mathcal{M}_k(\widetilde{\bcT})^\top$ with diagonal deletion (i.e., zero the diagonal values of $\mathcal{M}_k(\widetilde{\bcT}) \mathcal{M}_k(\widetilde{\bcT})^\top$) and showed that $\|\sin \Theta (\widetilde{\U}_k^{(0)}, \U_k)\| \leq \frac{\sqrt{2}}{2}$ holds with high probability when $|\Omega| \geq C p^{3/2}$, where $|\Omega|$ is the cardinality of $\Omega$. At the same time, they proved that HOSVD requires $|\Omega| \geq C p^{2}$ to achieve the same initialization performance and may not be an ideal initialization scheme for tensor completion. 
\end{itemize}
In addition, the random initialization is also widely considered in the literature. For example, \cite{anandkumar2014tensor} proposed to pick the best one among many random trials. It can be proved that if the number of random trials is large enough (usually polynomial in the dimension), one can find a trial such that the initialization is good enough \citep{anandkumar2014tensor}.

		
	\end{Remark}

	\begin{Remark}[Mode-$k$ Singular Subspace Linear Convergence Property]
		The upper bound in \eqref{ineq: max singular bound d = 3} includes two parts: a fixed quantity that represents the intrinsic estimation error, and another quantity that decays linearly to $0$ with respect to iteration index $t$. The linear convergence of HOOI was observed in \cite{ishteva2011best}, while Theorem \ref{th:HOOI d=3} gives a rigorous proof for it. Note that HOOI can be viewed as a special alternating minimization method, which was shown to have asymptotic linear convergence rate in solving nonlinear least squares problems \citep{ruhe1980algorithms}. This fact also sheds light on the linear convergence of HOOI.
	\end{Remark}

	\begin{Remark}[Unilateral Perturbation Bounds for Mode-$k$ Singular Subspace] \label{rem: unilateral bound}
		Our tensor perturbation bounds on singular subspace share the same spirit as the unilateral perturbation bounds on singular subspaces of matrix SVD in \cite{cai2018rate}. Consider the matrix perturbation setting mentioned in Section \ref{sec: error terms of Z} with the additional assumption that $\T$ is rank-$r$ and has SVD $\U \bSigma\V^\top$. \cite{cai2018rate} showed that the upper bound of $\|\sin\Theta(\widehat{\U}, \U)\|$ can be written as $\frac{a_1}{\lambda} + \frac{a_2}{\lambda^2}$, which can be interpreted as the sum of first and second order perturbations. In Theorem \ref{th:HOOI d=3}, the upper bound of $\left\|\sin\Theta\left(\widehat{\U}_k, \U_k\right)\right\|_q$ can be also written as $\sum_{i=1}^3 \frac{a_i}{\lambda^i}$ which can be interpreted as summation of the first, second, and third order perturbations. This phenomenon also generalizes to order-$d$ case in Theorem \ref{th:HOOI general}.
		
		Due to the unilateral property, when the tensor dimension of each mode is at different order, the estimation error rate of singular subspace in each mode can vary significantly. For example in the tensor denoising setting, $\widetilde{\bcT} = \bcT + \bcZ$ where $\bcT \in \bbR^{p_1 \times p_2 \times p_3}$ is a multilinear rank-$(r_1, r_2, r_3)$ tensor and $\bcZ$ is a random tensor with i.i.d. standard normal entries. Let $r_{\max} = \max_i r_i$ and suppose $p_1 \ll p_2 \ll p_3$, $r_{\max} \ll p_1^{1/2}$. Consider $q = \infty$, then by random matrix theory \citep{vershynin2010introduction}, we can show $\tau_{1i} \leq C\sqrt{p_i}$, $\tau_2 \leq C(\sqrt{p_3 r_{\max} })$ and $\tau_3 \leq  C(\sqrt{p_3 r_{\max} })$ with high probability. Thus when $\lambda \geq C p_3\sqrt{\frac{r_{\max}}{p_1}} $, Theorem \ref{th:HOOI d=3} immediately implies, with high probability
		\begin{equation*}
		\left\|\sin\Theta\left(\widehat{\U}_k, \U_k\right)\right\| \leq C \frac{\sqrt{p_k}}{\lambda}, \quad k = 1, 2, 3
		\end{equation*} 
		for some $C > 0$. So we can see the perturbation of $\widehat{\U}_k$ depends on $p_k$. Also as $\lambda$ decreases, for different $k$, different order perturbations in \eqref{ineq: singular space upp bnd d = 3} could dominate in the perturbation bound of $\widehat{\U}_k$. For example, when $\lambda \asymp C p_3 \sqrt{\frac{r_{\max}}{p_2}}$, Theorem \ref{th:HOOI d=3} yields
		\begin{equation*}
		\begin{split}
		&\left\|\sin\Theta\left(\widehat{\U}_1, \U_1\right)\right\| \leq C \frac{\sqrt{p_3^2 r_{\max}}}{\lambda^2}, \\
		& \left\|\sin\Theta\left(\widehat{\U}_2, \U_2\right)\right\| \leq C \frac{\sqrt{p_2}}{\lambda} + C' \frac{\sqrt{p_3^2 r_{\max}}}{\lambda^2},\\
		& \left\|\sin\Theta\left(\widehat{\U}_3, \U_3\right)\right\| \leq C \frac{\sqrt{p_3}}{\lambda}
		\end{split}
		\end{equation*} 
		for constants $C, C' > 0$. More details about the application of HOOI perturbation bounds in order-$d$ tensor denoising and numerical studies for this unilateral property of singular subspace perturbation
		can be found in Sections \ref{sec: app tensor denoising} and \ref{sec: numerical study determinitic part}, respectively.
	\end{Remark}
	
	\begin{Remark}[Comparing Perturbation Bounds of truncated HOSVD and HOOI] \label{rem: compare of HOSVD and HOOI}
		It is worth mentioning that the power iteration in Algorithm \ref{alg: HOOI} plays an important role for refining tensor reconstruction. Without power iteration, the estimator $\widehat{\bcT}^{\HOSVD} = \widetilde{\bcT} \times_1 P_{\widetilde{\U}_1^{(0)}} \times_2 P_{\widetilde{\U}_2^{(0)}} \times_3 P_{\widetilde{\U}_3^{(0)}}$ with $\widetilde{\U}_i^{(0)} = \SVD_{r_i} (\cM_i(\widetilde{\bcT}))$ is called truncated HOSVD (T-HOSVD) in the literature \citep{de2000multilinear}. It is not hard to show $\| \widehat{\bcT}^{\HOSVD} - \bcT \|_{\tHS} \leq C \|\bcZ\|_{\tHS}$ for some $C > 0$. Since $\|\bcZ\|_{\tHS} = \underset{\|\bcY \|_{\tHS} \leq 1}\sup \langle \bcZ, \bcY \rangle$ and may be much larger than $\xi$, we can see the power iteration can greatly improve the accuracy for tensor reconstruction, and this echos the findings in literature in tensor denoising setting \citep{zhang2018tensor}. 
	\end{Remark}
	
	The following lemma provides an alternative way to bound $\left\|\llbracket \bcZ; \widehat{\U}_1^\top, \widehat{\U}_2^\top, \widehat{\U}_3^\top \rrbracket\right\|_{\tHS}$ appearing in the reconstruction error bound \eqref{ineq:reconstruction-error-bound}.
	\begin{Lemma} \label{lm: noise projection error}
		Suppose $\bcZ\in \mathbb{R}^{p_1 \times \cdots \times p_d}$ is a general order-$d$ tensor and $\U_k,\widehat{\U}_k \in \mathbb{O}_{p_k, r_k}$ are general matrices with orthonormal columns. For any subset $\Omega\subseteq \{1, \ldots, d\}$, we further define projections of $\bcZ$ on $\Omega$ as follows,
		\begin{equation*}
		\theta_{\Omega} = \left\|\bcZ \times_{k \in \Omega} \U_k^\top \times_{k \in \Omega^c} \U_{k\perp}^\top \right\|_{\tHS}. 
		\end{equation*}
		Then, 
		\begin{equation*}
		\left\|\llbracket \bcZ; \widehat{\U}_1^\top, \widehat{\U}_2^\top,\ldots, \widehat{\U}_d^\top \rrbracket\right\|_{\tHS} \leq \sum_{\Omega\subseteq \{1,\ldots, d\}} \theta_{\Omega} \prod_{k\in \Omega^c}\left\|\sin\Theta(\widehat{\U}_k, \U_k)\right\|.
		\end{equation*}
	\end{Lemma}
	
	We end this section by introducing a deterministic rate matching lower bound for tensor reconstruction. Since the statement of the lower bound is relative simple, we state it in general order-$d$ setting. In particular, we consider the following class of $(\bcT, \bcZ)$ pairs of $p_1 \times \cdots \times p_d$ tensors and perturbations,
	\begin{equation*}
	\mathcal{F}_r(\xi) = \left\{ \left(\bcT, \bcZ  \right): \rank(\bcT) = (r_1, \cdots, r_d) \leq \r,\underset{\substack{\|\bcY \|_{\tHS} \leq 1,\\  \rank(\bcY) \leq (r_1, \ldots, r_d)}}\sup \langle \bcZ, \bcY \rangle \leq \xi  \right\},
	\end{equation*} here $\r = (r,\ldots, r)$ and the comparison is entrywise.
	
	\begin{Theorem}[Tensor Reconstruction Lower Bound under Perturbation]\label{th: pertur lower bound}
		Consider perturbation model \eqref{eq: perturbation model}, we have the following deterministic lower bound for reconstructing $\bcT$, 
		\begin{equation*}
		\inf_{\widehat{\bcT}} \sup_{(\bcT, \bcZ) \in \mathcal{F}_r(\xi)} \|\widehat{\bcT} - \bcT\|_{\tHS} \geq \frac{\sqrt{2}}{2} \xi.
		\end{equation*}
	\end{Theorem}
	
	\begin{Remark}[Optimality of HOOI and one-step HOOI for Tensor Reconstruction] \label{rem: one-step optimality of HOOI}
		When tensor order $d$ is fixed, combining Theorem \ref{th:HOOI d=3} and \ref{th: pertur lower bound}, we have shown that HOOI with good initialization is optimal for tensor reconstruction in the class $(\bcT, \bcZ) \in \mathcal{F}_r(\xi)$. At the same time, from \eqref{ineq: t step tensor reconstruction}, we see the error rate of tensor reconstruction is optimal even after one iteration of HOOI i.e., $t_{max} = 1$ and more iterations can improve the coefficient in front of $\xi$. This suggests that in some applications where running HOOI until convergence is prohibitive, we can just run it for one iteration to get a fairly good reconstruction. See more in Section \ref{sec: numerical study: comparison} about a numerical comparison of HOOI and one-step HOOI.
		
		Apart from the optimality of our perturbation bound in tensor reconstruction, it is also interesting to study whether the perturbation bounds in \eqref{ineq: d=3-max-upper bound of singular space}, \eqref{ineq: singular space upp bnd d = 3} for singular subspaces are optimal or not and we leave it as an interesting future work.
	\end{Remark}

	\section{A Blockwise Perturbation Bound of Higher-order Orthogonal Iteration for Tensor Decomposition}\label{sec:main}
	
	In this section, we present the main results of perturbation bounds of HOOI given in Algorithm \ref{alg: HOOI}. In contrast with Theorem \ref{th:HOOI d=3}, Theorem \ref{th:HOOI general} in this section covers the general order-$d$ perturbation setting with $\widetilde{\bcT}$ having symmetric index groups $\{\Omega_i \}_{i=1}^m$. Before stating the theorem, we first define the blockwise errors of $\bcZ$ in this general setting. Let $\cS_{i}^{(-\bar{k})} := \left\{ S \subseteq [d]\setminus \{\bar{k}\}: |S| = i \right\} $ be the collection of all possible index sets with $i$ elements from $[d]\setminus \{\bar{k} \}$ and $\cS_0^{(-\bar{k})}:= \emptyset$. For $S \in \cS_{i}^{(-\bar{k})}$, we let $S^c = ([d]\setminus \{\bar{k}\}) \setminus S$. Now we define the blockwise errors of $\bcZ$ as 
	\begin{equation} \label{eq: blockwise error order d}
	\begin{split}
	&  \tau_1 = \max_{k=1,\ldots, m} \tau_{1k}, \quad \tau_{1k} = \left\|  \left(\mathcal{M}_{\bar{k}}(\bcZ \times_{i \neq \bar{k} }\U_{i'}^\top ) \right)_{\max(r_k)}\right\|_q, \quad k=1,\ldots, m;\\
	& \tau_j = \max_{k \in [m]}\Big\{ \max_{S \in \cS^{(-\bar{k})}_{j-1}} \sup_{ \substack{\V_{i'} \in \bbR^{(p_{i'} - r_{i'}) \times r_{i'}},\\ \|\V_{i'}\|_q \leq 1, i \in S } } \left\|\left(\mathcal{M}_{\bar{k}}(\bcZ \times_{i \in S}\left(\U_{i'\perp} \V_{i'} \right)^\top \times_{i \in S^c}\U_{i'}^\top \right)_{\max(r_k)}\right\|_q \Big\},\\
	&\quad \text{ for } j = 2, \ldots, m.
	\end{split}
	\end{equation}
	Finally, $\xi$ in this setting is defined as
	\begin{equation} \label{eq: blockwise error order d xi}
	\xi := \underset{\|\bcY \|_{\tHS} \leq 1,  \rank(\bcY) \leq (r_{1'}, \ldots, r_{d'})}\sup \langle \bcZ, \bcY \rangle.
	\end{equation}
	\begin{Theorem}[General Perturbation Bounds for Tensor Power Iteration]\label{th:HOOI general}
		Consider the perturbation model \eqref{eq: perturbation model} with $\widetilde{\bcT}, \bcT, \bcZ \in \mathbb{R}^{p_{1'}\times \cdots \times p_{d'}}$, symmetric index groups $(\Omega_1, \ldots, \Omega_m)$ and blockwise errors in \eqref{eq: blockwise error order d} \eqref{eq: blockwise error order d xi}. Suppose $q\geq1$. Denote the initialization errors of $\{\widetilde{\U}_k^{(0)} \}_{k=1}^m$ as $\tilde{e}_0 := \max_{k=1, \ldots, m}\|\widetilde{\U}^{(0)\top}_{k\perp} \U_k\|$, $e_0 := \max_{k=1, \ldots, m}\|\widetilde{\U}^{(0)\top}_{k\perp} \U_k\|_q$. Assume the initialization error and the signal strength satisfy 
		\begin{equation}\label{ineq: signal strength condition general d}
		\tilde{e}_0 \leq \frac{\sqrt{2}}{2} \, \text{ and } \, \lambda \geq 2^{\frac{d+4}{2}} \left(1+\frac{\sqrt{2}}{2}\right)^d \xi.
		\end{equation} 
		
		Let $\widetilde{\bcT}^{(t)}:= \widetilde{\bcT} \times_{\Omega_1} P_{\widetilde{\U}^{(t)}_1} \times \cdots \times_{\Omega_m} P_{\widetilde{\U}^{(t)}_m}$ be the estimator of $\bcT$ after $t$ steps in Algorithm \ref{alg: HOOI}. Then with inputs $\widetilde{\bcT}$, $\{\widetilde{\U}_k^{(0)} \}_{k=1}^m$, $\{\Omega_i \}_{i=1}^m$, the mode-$k$ singular subspace updates in Algorithm \ref{alg: HOOI} after $t$ iterations satisfy
		\begin{equation} \label{ineq: t step singular space error}
		\max_{k \in [m]} \left\|\sin\Theta\left(\widetilde{\U}_k^{(t)}, \U_k\right)\right\|_q \leq 2^{\frac{d+3}{2}} \frac{\tau_1}{\lambda} + \frac{e_0}{2^{t}},
		\end{equation} and the $t$-step tensor estimation satisfies
		\begin{equation} \label{ineq: t step tensor reconstruction order d}
		\left\|\widetilde{\bcT}^{(t)} - \bcT\right\|_{\tHS} \leq \left(  1+ 2d \left( 1- \left( 2^{\frac{d+3}{2}} \frac{\tau_1}{\lambda} + \frac{\tilde{e}_0}{2^{t-1}} \right)^2 \right)^{- \frac{d-1}{2} }  \right)\xi.
		\end{equation}
		
		Moreover, when $t_{\max} \geq C \log(e_0\lambda/\tau_{1})\vee 1$ for some $C > 0$, the outputs of estimated mode-$k$ singular subspace of Algorithm \ref{alg: HOOI} satisfy
		\begin{equation*}
		\max_{k \in [m]} \left\|\sin\Theta\left(\widehat{\U}_k, \U_k\right)\right\|_q \leq \left(2^{\frac{d+3}{2}} + 1\right) \frac{\tau_1}{\lambda},
		\end{equation*}
		\begin{equation}\label{ineq: upper bound of singular space}
		\begin{split}
		\left\|\sin\Theta\left(\widehat{\U}_k, \U_k\right)\right\|_q & \leq \frac{2}{\left( 1- c^*(\tau_1,\lambda,d) \right)^{ \frac{d-1}{2} } } \left(\frac{\tau_{1k}}{\lambda} +  \sum_{j=1}^{d-1}\frac{{d-1 \choose j}  \left(2^{\frac{d+3}{2}}+1\right)^j \tau_1^j  \tau_{j+1} }{\lambda^{j+1}} \right).
		\end{split}
		\end{equation}
		for $k=1,\ldots,m$ where $c^*(\tau_1,\lambda,d) := \frac{\left(2^{\frac{d+3}{2}} + 2\right)^2 \tau_1^2}{\lambda^2} \leq \frac{1}{2}$,
		and the output of tensor reconstruction $\widehat{\bcT}$ satisfies
		\begin{equation} \label{ineq: tensor reconstruction}
		\begin{split}
		\left\|\widehat{\bcT} - \bcT\right\|_{\tHS} \leq &\left\|  \bcZ\times_{\Omega_1} \widehat{\U}_1^\top \times \cdots \times_{\Omega_m} \widehat{\U}_m^\top \right\|_{\tHS} + \sum_{k=1}^m |\Omega_k| \left\|\widehat{\U}_{k\perp}^\top \mathcal{M}_{\bar{k}}(\bcT) \right\|_F \\
		\leq & \left(  1+ 2d \left( 1- c^*\left(\tau_1,\lambda,d\right) \right)^{- \frac{d-1}{2} }  \right)\xi.
		\end{split}
		\end{equation}
	\end{Theorem}
	
	\begin{Remark}[Size of $c^*(\tau_1, \lambda, d)$] \label{rem: size of xi constant}
		It is easy to check $c^*(\tau_1, \lambda, d) \leq \frac{1}{2}$ based on $\tau_1 \leq \xi$ and the requirement of the signal strength $\lambda$. So we have $\left(1- c^*(\tau_1,\lambda,d) \right)^{- \frac{d-1}{2} } \leq 2^{\frac{d-1}{2}}$ in the upper bounds of $\left\|\sin\Theta\left(\widehat{\U}_k, \U_k\right)\right\|_q$ and $\left\|\widehat{\bcT} - \bcT\right\|_{\tHS}$. However, in many practical applications, such as tensor denoising to be introduced in Section \ref{sec: app tensor denoising}, tensor order $d$ is fixed and $\lambda \gg (2^{\frac{d+3}{2}} + 2) \tau_1$. In this case $c^*(\tau_1, \lambda, d)$ could be much smaller than $\frac{1}{2}$ and the scale of $\left( 1- c^*(\tau_1,\lambda,d) \right)^{- \frac{d-1}{2} }$ can be very close to $1$.
	\end{Remark}

	\begin{Remark}[Comparison with \cite{anandkumar2014tensor}] \label{rem: compare with Anandkumar}
		Compared with the perturbation bounds of power iteration for supersymmetric CP-low-rank decomposition \cite[Theorem 5.1]{anandkumar2014tensor}, our Theorem \ref{th:HOOI general} covers more general symmetric and partial symmetric multilinear low-rank decomposition settings. Also in Theorem 5.1 of \cite{anandkumar2014tensor}, the tensor reconstruction error bound of power iteration is given in terms of tensor spectral norm, which does not improve upon the guarantee by the trivial estimator $\widetilde{\bcT}$. On the other hand, the tensor reconstruction error of $\widehat{\bcT}$ in Theorem \ref{th:HOOI general} is given in Hilbert-Schmidt norm and can be significantly better than the guarantee for $\widetilde{\bcT}$ as $\|\bcZ\|_{\tHS} \gg \xi$ in most of the applications. 
	\end{Remark}
	
	\begin{Remark}[Dependence on Tensor Order $d$]
	We note that in Theorem \ref{th:HOOI general}, the constants in our condition \eqref{ineq: signal strength condition general d} and perturbation bounds \eqref{ineq: t step singular space error} and \eqref{ineq: t step tensor reconstruction order d} scales exponentially w.r.t. the tensor order $d$. We think this exponential dependence on $d$ is not sharp. In fact, in Theorem 1 of the recent work \cite{luo2021low}, they show the dependence on $d$ in \eqref{ineq: signal strength condition general d} and \eqref{ineq: t step tensor reconstruction order d} can be reduced to poly$(d)$.   
	\end{Remark}

	\begin{Remark}[A Proof Sketch of Theorem \ref{th:HOOI general}] \label{rem: proof sketch of main theorem} 
	We provide a sketch on how to prove \eqref{ineq: t step singular space error} and \eqref{ineq: t step tensor reconstruction order d}. The rest of the results \eqref{ineq: upper bound of singular space} and \eqref{ineq: tensor reconstruction} follow easily from \eqref{ineq: t step singular space error}, \eqref{ineq: t step tensor reconstruction order d} by plugging in $t_{\max} \geq \log(e_0\lambda/\tau_{1})\vee 1$. The idea is to develop the recursive error bounds of $\widetilde{\U}_k^{(t+1)}$, i.e., the estimate of $\U_k$ at iteration $t+1$, based on the error bound of $\widetilde{\U}_k^{(t)}$, i.e., the estimate at iteration $t$. The argument can be divided into three steps. It is worth mentioning that all three steps involves complex tensor algebra and this makes the proof even more difficult. 
		
		First, we denote
		\begin{equation*}
		\begin{split}
		& \tilde{e}_{t} = \max_k \tilde{e}_{t,k},\quad \tilde{e}_{t,k} = \left\|(\widetilde{\U}_{k\perp}^{(t)})^\top \U_k \right\|,\\
			& e_{t} = \max_k e_{t,k},\quad e_{t,k} = \left\|(\widetilde{\U}_{k\perp}^{(t)})^\top \U_k \right\|_q,\quad k=1,\ldots,m; t=0,1,\ldots.
		\end{split}
		\end{equation*}
		
		{\bf Step 1}. In HOOI procedure, the update for the mode-$k$ singular subspace satisfies
		\begin{equation*}
		\begin{split}
		\widetilde{\U}_k^{(t+1)} 
		= & \SVD_{r_k}\bigg( \mathcal{M}_{\bar{k}}\left(\bcT \times_{i \in \underline{\Omega}_k} \widetilde{\U}^{(t+1)\top}_{i'} \times_{\widecheck{\Omega}_k} \widetilde{\U}^{(t)\top}_k  \times_{i \in \overline{\Omega}_k} \widetilde{\U}^{(t)\top}_{i'}   \right) \\
		& \quad \quad \quad + \mathcal{M}_{\bar{k}}\left( \bcZ \times_{i \in \underline{\Omega}_k} \widetilde{\U}^{(t+1)\top}_{i'} \times_{\widecheck{\Omega}_k} \widetilde{\U}^{(t)\top}_k  \times_{i \in \overline{\Omega}_k} \widetilde{\U}^{(t)\top}_{i'}  \right)\bigg),
		\end{split}
		\end{equation*} here $\underline{\Omega}_i :=\bigcup_{j=1}^{i-1} \Omega_j, \overline{\Omega}_i := \bigcup_{j=i+1}^d \Omega_j$.
		To give an upper bound for $e_{t+1,k}$, we aim to give an upper bound for 
		\begin{equation} \label{eq: step1 upper bound}
		\left\|\left(\mathcal{M}_{\bar{k}}\left(\bcZ \times_{i \in \underline{\Omega}_k} \widetilde{\U}^{(t+1)\top}_{i'} \times_{\widecheck{\Omega}_k} \widetilde{\U}^{(t)\top}_k  \times_{i \in \overline{\Omega}_k} \widetilde{\U}^{(t)\top}_{i'}  \right) \right)_{\max(r_k)}\right\|_q
		\end{equation}
		by using $\tau_1, \ldots, \tau_m$, $e_t, e_{t+1}$. The main idea to bound \eqref{eq: step1 upper bound} is to introduce $\I = P_{\U_{i'}} + P_{\U_{i'  \perp}}$ in each mode multiplication, expand the mode products, then write the whole term into summation of many small terms.
		
		{\bf Step 2.} After getting an upper bound for \eqref{eq: step1 upper bound}, we use induction to prove the following claim,
		\begin{equation*}
		e_{t} \leq 2^{(d+3)/2}\tau_{1}/\lambda + e_0/2^t, \quad \tilde{e}_{t} \leq 2^{(d+3)/2}\tau_{1}/\lambda + \tilde{e}_0/2^t; \quad t=0,1,\ldots. 
		\end{equation*} One technical difficulty is to deal with the sequential updating of singular subspaces in HOOI and we use the induction idea again to tackle it. 
		Tools we use in this step include the singular subspace bound in \cite[Theorem 5]{luo2020schatten}. 
		
		{\bf Step 3.} The final and most challenging step involves upper bounding the tensor reconstruction error $\| \widetilde{\bcT}\times_{\Omega_1} P_{\widehat{\U}_1}\times \cdots \times_{\Omega_m} P_{\widehat{\U}_m} - \bcT\|_{\tHS}$ by the unified quantity $\xi$. By decomposing $\bcT$ onto the estimated singular subspaces, we can show that
		\begin{equation*}
		\begin{split}
		& \left\|\widetilde{\bcT} \times_{\Omega_1} P_{\widehat{\U}_1} \times \cdots \times_{\Omega_m} P_{\widehat{\U}_m} - \bcT \right\|_{\tHS} \\
		\leq & \left\| \bcZ \times_{\Omega_1} P_{\widehat{\U}_1} \times \cdots \times_{\Omega_m} P_{\widehat{\U}_m}\right\|_{\tHS} + \sum_{k=1}^d \left\|\widehat{\U}_{k'\perp}^\top \mathcal{M}_{k}(\bcT)\right\|_F.
		\end{split}
		\end{equation*}
		By definition, $\left\| \bcZ \times_{\Omega_1} P_{\widehat{\U}_1} \times \cdots \times_{\Omega_m} P_{\widehat{\U}_m}\right\|_{\tHS} \leq \xi$. We further show
		\begin{equation*}
		\begin{split}
		\left\|\widehat{\U}_{k'\perp}^\top \mathcal{M}_{k}(\bcT)\right\|_F
		\overset{(a)}\leq C(\tau_1,d,\lambda) \, \xi.
		\end{split}
		\end{equation*}
		Here $C(\tau_1,d,\lambda)$ is a quantity that depends on $\tau_1,d,\lambda$. 
		The main challenge to prove (a) is that $\widehat{\U}_{k'}$ is not the left singular subspace of $\cM_k(\widetilde{\bcT})$. So to leverage the SVD property of $\widehat{\U}_{k'\perp}$, we have to project $\widetilde{\bcT}$ onto $\widetilde{\U}_i^{(t_{\max})}$ and $\widetilde{\U}_i^{(t_{\max}-1)}$, then use the subspace perturbation bounds established before. 
		
	\end{Remark}

	Note that Theorem \ref{th:HOOI general} covers the general situation where $\bcT$ may have partial symmetric modes. We provide a corollary for the common asymmetric case, i.e., $\Omega_i = \{i\}, i = 1, \ldots, d$ in the Appendix.
	
	\section{Perturbation Bounds of Power Iteration for Tensors with Partial Low Multilinear Rank Structure} \label{sec: HOOI in partial low rank}
	
	In some applications, e.g., multilayer network analysis \citep{lei2019consistent}, the tensor $\widetilde{\bcT}$ only has low-rank structure on a subset of modes. Both the tensor power iteration algorithm and our perturbation theory can be generalized to such cases. For a better illustration, we present the modified algorithm and theory for order-3 tensor perturbation with mode $1$ to be dense. Specifically, we consider 
	\begin{equation}\label{eq: perturbation model in partial low rank}
	\widetilde{\bcT} = \bcT + \bcZ \in \bbR^{p_1 \times p_2 \times p_3},
	\end{equation}
	where $\bcT$ is the signal tensor and $\bcZ$ is the noise. We assume $\bcT$ is low-rank on modes $2$ and $3$, i.e., $\bcT = \bcS \times_2 \U_2 \times_3 \U_3$, where $\bcS \in \bbR^{p_1 \times r_2 \times r_3}$ is the core tensor and $\U_i \in \mathbb{O}_{p_i,r_i}$ for $i = 2, 3$. In this setting, we consider the modified tensor power iteration algorithm for low-rank tensor decomposition in Algorithm \ref{alg: HOOI in d = 3 partial low rank}.
	
	\begin{algorithm}[!h] \caption{Power Iteration for Tensor Decomposition in Partial Multilinear Low-Rank Setting \eqref{eq: perturbation model in partial low rank}} \label{alg: HOOI in d = 3 partial low rank}
		\textbf{Input:} $\widetilde{\bcT}\in \bbR^{p_{1} \times p_2 \times p_{3}}$, initialization $\{\widetilde{\U}_i^{(0)}\in \mathbb{O}_{p_i,r_i}\}_{i=2}^3$, maximum number of iterations $t_{\max}$.\\
		\textbf{Output:} $ \{\widehat{\U}_i\}_{i=2}^3, \widehat{\bcT}$. 
		\begin{algorithmic}[1]
			\State For $t=1, \ldots, t_{\max}$, do
			\begin{equation*}
			\begin{split}
			&\widetilde{\U}_2^{(t+1)} = \SVD_{r_2} \left(\cM_{2}\left( \widetilde{\bcT} \times_3 (\widetilde{\U}_{3}^{(t)})^\top  \right) \right)\\
			&\widetilde{\U}_3^{(t+1)} = \SVD_{r_3} \left(\cM_{3}\left( \widetilde{\bcT} \times_2 (\widetilde{\U}_{2}^{(t+1)})^\top  \right) \right).
			\end{split}
			\end{equation*}
			\State Let $\widehat{\U}_i = \widetilde{\U}_i^{(t_{\max})}$ for $i = 2, 3$ and compute 
			\begin{equation*}
			\widehat{\bcT} = \widetilde{\bcT}\times_2 P_{\widehat{\U}_2} \times_{3} P_{\widehat{\U}_3}.
			\end{equation*}
		\end{algorithmic}
	\end{algorithm}
	In the setting \eqref{eq: perturbation model in partial low rank}, we can define the blockwise errors of $\bcZ$ as follows:
	\begin{equation} \label{eq: blockwise error d = 3 partial low rank}
	\begin{split}
	& \tau_1 = \max_{k= 2,3} \tau_{1k}, \quad \tau_{12} = \left\|  \left(\mathcal{M}_2(\bcZ \times_{3}\U_{3}^\top ) \right)_{\max(r_2)}\right\|_q, \tau_{13} = \left\|  \left(\mathcal{M}_3(\bcZ \times_{2}\U_{2}^\top ) \right)_{\max(r_3)}\right\|_q \\
	& \tau_2 = \max \Big\{\max_{\substack{\V\in \mathbb{R}^{(p_{3}-r_{3})\times r_{3}}\\ \|\V\|_q\leq 1}}\left\|\left(\mathcal{M}_2(\bcZ \times_{3}(\U_{3\perp}\V)^\top ) \right)_{\max(r_2)}\right\|_q, \\
	&\quad\quad\quad\quad  \max_{\substack{\V\in \mathbb{R}^{(p_{2}-r_{2})\times r_{2}}\\ \|\V\|_q\leq 1}}\left\|\left(\mathcal{M}_3(\bcZ \times_{2}(\U_{2\perp}\V)^\top ) \right)_{\max(r_3)}\right\|_q \Big\};\\
	& \xi = \sup_{\rank(\bcY) \leq (p_1, r_2, r_3), \|\bcY\|_{\tHS} \leq 1 } \left\langle
	\bcZ, \bcY \right\rangle.
	\end{split}
	\end{equation}
	
	We have the following perturbation bounds for the outputs of Algorithm \ref{alg: HOOI in d = 3 partial low rank}.
	\begin{Theorem}[Tensor Perturbation Bounds with Partial Multilinear Low-Rank]\label{th:HOOI d=3 partial low rank}
		Consider the perturbation model \eqref{eq: perturbation model in partial low rank} with $\widetilde{\bcT}, \bcT, \bcZ \in \mathbb{R}^{p_1\times p_2\times p_3}$. Suppose $q\geq1$. Define the blockwise errors as in \eqref{eq: blockwise error d = 3 partial low rank} and denote the initialization errors of $\{ \widetilde{\U}_k^{(0)} \}_{k=2}^3$ as $\tilde{e}_0 := \max_{k=2,3}\|\widetilde{\U}^{(0)\top}_{k\perp} \U_k\|$, $e_0 := \max_{k=2,3}\|\widetilde{\U}^{(0)\top}_{k\perp} \U_k\|_q$. Assume the initialization error and the signal strength satisfy 
		\begin{equation} \label{ineq: signal strength condition d = 3 partial low rank}
		\tilde{e}_0 \leq \sqrt{2}/2\, \text{ and } \,\lambda \geq 16\xi.
		\end{equation} 
		Then with inputs $\widetilde{\bcT}$, $\{ \widetilde{\U}_i^{(0)} \}_{i=2}^3$, the mode-$k$ singular subspace updates in Algorithm \ref{alg: HOOI in d = 3 partial low rank} after $t$ iterations satisfy
		\begin{equation} \label{ineq: max singular bound d = 3 partial low rank}
		\max_{k = 2,3} \left\|\sin\Theta\left(\widetilde{\U}^{(t)}_k, \U_k\right)\right\|_q \leq \frac{4\sqrt{2}\tau_1}{\lambda} + \frac{e_0}{2^{t}}.
		\end{equation} 
		Moreover, when $t_{\max} \geq C\log(e_0 \lambda /\tau_{1})\vee 1$ for some $C > 0$, the outputs of estimated mode-$k$ singular subspace of Algorithm \ref{alg: HOOI in d = 3 partial low rank} satisfy
		\begin{equation*}
		\max_{k = 2,3} \left\|\sin\Theta\left(\widehat{\U}_k, \U_k\right)\right\|_q \leq  \frac{(4\sqrt{2} + 1)\tau_1}{\lambda},
		\end{equation*}
		\begin{equation}
		\begin{split} \label{ineq: singular space upp bnd d = 3 partial low rank}
		\left\|\sin\Theta\left(\widehat{\U}_k, \U_k\right)\right\|_q = \left\| \widehat{\U}_{k\perp}^{\top} \U_k \right\|_q \leq 2\sqrt{2} \left(\frac{\tau_{1k}}{\lambda} +  \frac{(4\sqrt{2} + 1)\tau_1 \tau_2}{\lambda^2} \right), \quad k=2,3.
		\end{split}
		\end{equation}
		The output of tensor reconstruction $\widehat{\bcT}$ satisfies
		\begin{equation*}
		\begin{split}
		\left\|\widehat{\bcT} - \bcT\right\|_{\tHS} & \leq \left\| \bcZ \times_2 P_{\widehat{\U}_2} \times_3 P_{\widehat{\U}_3} \right\|_{\tHS}+ \sum_{k=2}^3 \left\|\widehat{\U}_{k\perp}^\top \mathcal{M}_k(\bcT)\right\|_F \leq (4\sqrt{2} + 1) \xi.
		\end{split}
		\end{equation*}
	\end{Theorem}
	The proof of Theorem \ref{th:HOOI d=3 partial low rank} follows the proof of Theorem \ref{th:HOOI general} easily. For simplicity, we omit it here.
	
	\section{Implications in Statistics and Machine Learning} \label{sec: applications}
	In this section, we consider a couple of applications of the HOOI perturbation bounds we developed in statistics and machine learning. Specifically, here we consider the perturbation model \eqref{eq: perturbation model} and assume $\bcZ_{i_1, \ldots, i_d}$'s are independent, mean-zero $\sigma$-subgaussian, where $\sigma > 0$ is the subgaussianity parameter. More precisely,
	\begin{equation*}
	\bbE \exp(\lambda \bcZ_{i_1, \ldots, i_d}) \leq \exp(C \lambda^2 \sigma^2 ),\, \text{for all} \, (i_1, \ldots, i_d) \in [p_1] \times \cdots \times [p_d] \text{ and all } \lambda \in \bbR,
	\end{equation*} where $C > 0$ is some absolute constant. For convenience, we let $p_{\max} = \max_i p_i, p_{\min} = \min_i p_i, r_{\max} = \max_i r_i, r_{\min} = \min_i r_i$.
	
	In this setting, the quantity $\xi$ is in fact closely related to the Gaussian width \citep{gordon1988milman} studied in literature that measures the size or complexity of a given set. Recall the Gaussian width of a set $\cS \subset \bbR^{p_1 \times \cdots \times p_d}$ is defined to be
	\begin{equation*}
	w(\cS) := \bbE \left(  \sup_{\bcY \in \cS} \left\langle \bcB, \bcY \right\rangle  \right),
	\end{equation*} where $\bcB \in \bbR^{p_1 \times \cdots  \times p_d}$ is a tensor whose entries are independent $N(0,1)$ random variables. In view of $w(\cS)$, we can regard $\xi$ as the Gaussian width with no expectation and $\cS = \{\bcY : \|\bcY\|_{\tHS} \leq 1, \rank(\bcY) \leq (r_1, \ldots, r_d) \}$. It can be shown in the case that $\bcZ$ has i.i.d. $N(0,1)$ entries, $\xi$ and $w(\cS)$ are the same up to constant with high probability \citep{raskutti2019convex}.
	
	In the following subsections, we consider two particular structures of $\bcT$, one is pure low multilinear rank structure, namely tensor denoising/tensor PCA/tensor SVD studied in literature \citep{richard2014statistical,zhang2018tensor,perry2020statistical,hopkins2015tensor} and another one is tensor co-clustering/block structure \citep{chi2020provable,wang2019multiway}.

	\subsection{HOOI for Tensor Denoising} \label{sec: app tensor denoising}
	In tensor denoising, we assume $\bcT$ has the following structure,
	\begin{equation}\label{eq: tensor denoising T}
	\bcT = \bcS \times_1 \U_1 \times \cdots \times_d \U_d,
	\end{equation} 
	where $\bcS \in \bbR^{r_1 \times \cdots \times r_d}$ is the core tensor and $\{\U_i \in \mathbb{O}_{p_i, r_i} \}_{i=1}^d $ are loading matrices. With the established tensor perturbation bounds, we can establish the following theoretical guarantee for the performance of HOOI on tensor denoising with a very short proof.
	\begin{Theorem}[Tensor Denoising: General Order $d$]\label{th: tensor denoising}
		Consider the tensor denoising problem $``\eqref{eq: perturbation model}+\eqref{eq: tensor denoising T}"$ and Algorithm \ref{alg: HOOI} with inputs $\widetilde{\bcT}$, $\Omega_i = \{i \}$, initialization $\{ \widetilde{\U}_i^{(0)} \}_{i=1}^d$ and $t_{\max} = C \left( \log \frac{\lambda/\sigma}{\sqrt{p_{\max}}} \vee 1  \right)$ for some $C > 0$, where $\lambda = \min_k \sigma_{r_k}\left(\mathcal{M}_k(\bcS)\right)$ is the minimal singular value of each matricization of $\bcS$. Assume $r_{\max} \leq p_{\min}^{1/(d-1)}$ and $\max_i\|\widetilde{\U}_{i\perp}^{(0)\top} \U_i\| \leq \sqrt{2}/2$. Then if $\lambda/\sigma \geq 2^{(d+4)/2} (1 + \sqrt{2}/2)^{d} \sqrt{p_{\max} r_{\max}}$, with probability at least $1 - \exp(-c p_{\min})$, the output $\widehat{\U}_k, \widehat{\bcT}$ satisfy
		\begin{equation*}
		\left\|\sin\Theta\left(\widehat{\U}_k, \U_k\right)\right\| \leq C \left( 1- \frac{c}{r_{\max}} \right)^{-(d-1)/2}  \left( \frac{\sqrt{p_k}}{\lambda /\sigma } + \frac{ \sqrt{p_{\max}^2 r_{\max}}   }{\left(\lambda /\sigma  \right)^2}  \right)
		\end{equation*} 
		and 
		\begin{equation} \label{ineq: HOOI in tensor denoising}
		\|\widehat{\bcT} - \bcT\|_{\tHS} \leq C \left(1 + 2d\left(1-\frac{c}{r_{\max}} \right)^{-\frac{d-1}{2}} \right) \sigma \sqrt{ \sum_{i=1}^d p_i r_i },
		\end{equation} for some constants $c, C > 0$.
	\end{Theorem}
	
	When $d$ is a constant, the upper bound for tensor reconstruction error matches the lower bound in \cite[Theorem 3]{zhang2018tensor}, which shows HOOI achieves the optimal tensor reconstruction error in the tensor denoising problem.

	\subsection{Tensor Co-clustering/Block Model} \label{sec: app-tensor block model}
	
	Co-clustering is among the most important unsupervised learning methods that reveals the checkerbox-like association pattern in data. A number of algorithms have been proposed \citep{wang2015multi,wu2016general,kolda2008scalable,papalexakis2012k,sun2009multivis,jegelka2009approximation} for tensor co-clustering in the literature, however most of the work does not provide statistical guarantees for recovering the underlying co-clustering structure. Very recently, \cite{chi2020provable} and \cite{wang2019multiway} studied the performance of co-clustering estimation and cocluster recovery based on convex relaxation and combinatorial search algorithms. By using the tools of perturbation bounds of HOOI given in Section \ref{sec:main}, we are able to provide the first guarantee for co-clustering estimation and cocluster recovery based on computational efficient HOOI algorithm. Compared to the convex relaxation approach \citep{chi2020provable}, HOOI has a better guarantee for tensor reconstruction and it also gives guarantee for cocluster membership recovery. Specifically, in the tensor co-clustering/block model, we assume $\bcT$ has the following structure,
	\begin{equation} \label{eq: tensor block model}
	\bcT = \bcB \times_1 \bPi_1 \times \cdots \times_d \bPi_d,
	\end{equation} 
	where $\bPi_i \in \mathbb{M}_{p_i, r_i}$ and $\mathbb{M}_{p_i, r_i}$ is the collection of all $p_i \times r_i$ membership matrices with each row has exactly one $1$ and $(r_i - 1)$ $0$'s. For any $\bPi_i$, the cocluster membership of node $k$ is denoted by $g_i^{(k)} \in [r_i]$, which satisfies $(\bPi_i)_{[k,g_i^{(k)}]} = 1$. Let $G^{(j)}_i \equiv G^{(j)}_i(\bPi_i) = \{k \in [p_i]: g_i^{(k)} = j \}$ be the set of $i$th mode node indices that belongs to cocluster $j$ and $p^{(j)}_i = |G^{(j)}_i|$ for all $j \in [r_i]$. For simplicity, we assume the cocluster sizes for each cluster are on the same order for every mode, i.e.,
	\begin{equation}\tag{A1} \label{assu: cocluster size assumption}
	p_i^{(1)} \asymp p_i^{(2)} \asymp \cdots \asymp p_{i}^{(r_i)} \asymp \frac{p_i}{r_i}, \text{ for } i = 1, \ldots, d.
	\end{equation}

	We consider two cocluster membership recovery error metrics:
	\begin{enumerate}
		\item Let $\|\M\|_0$ be the number of nonzero entries in matrix $\M$. Suppose $\E_{r_i}$ is the set of all $r_i \times r_i$ permutation matrices. Define
		\begin{equation} \label{eq: misclassification error rate}
		err(\widehat{\bPi}_i, \bPi_i) = \frac{1}{p_i} \min_{\J \in \E_{r_i}} \| \widehat{\bPi}_i \J - \bPi_i  \|_0
		\end{equation} 
		as the misclassification rate of $\widehat{\bPi}_i$.
		\item Define
		\begin{equation*}
		\widetilde{err}(\widehat{\bPi}_i, \bPi_i) = \min_{\J \in \E_{r_i}} \max_{1 \leq j \leq r_i} \frac{1}{p^{(j)}_i} \| (\widehat{\bPi}_i \J)_{[G^{(j)}_i,:]} - (\bPi_i)_{[G^{(j)}_i,:]}  \|_0.
		\end{equation*} 
		Intuitively speaking, $\widetilde{err}(\widehat{\bPi}_i, \bPi_i)$ measures the worst relative misclassification rates over all communities.
	\end{enumerate}
	It is easy to check that $0 \leq err(\widehat{\bPi}_i, \bPi_i) \leq \widetilde{err}(\widehat{\bPi}_i, \bPi_i) \leq 2$.

	The following lemma gives a Tucker decomposition of $\bcT$ in the tensor block model \eqref{eq: tensor block model}. This decomposition bridges the tensor co-clustering model and the tensor signal-plus-noise model, which explains why HOOI would work for tensor co-clustering.
	\begin{Lemma} \label{lem: why HOOI work in tensor block}
		Suppose $\bcT$ has the tensor co-clustering/block structure \eqref{eq: tensor block model}, where $\bcB$ is a multilinear rank-$(r_1, \ldots, r_d)$ tensor. Assume $\bcS \times_1 \V_1 \times \cdots \times_d \V_d$ with $\V_i \in \mathbb{O}_{r_i, r_i}$ is the Tucker decomposition of $\bcB \times_1 (\bPi_1^\top \bPi_1)^{\frac{1}{2}} \times \cdots \times (\bPi_d^\top \bPi_d)^{\frac{1}{2}}$. Then 
		\begin{equation*}
		\bcT = \bcS \times_1 \U_1 \times \cdots \times_d \U_d
		\end{equation*} 
		with $\U_i = \bPi_i(\bPi_i^\top \bPi_i)^{-\frac{1}{2}}\V_i \in \mathbb{O}_{p_i, r_i}$ for $i = 1, \ldots, d$.
	\end{Lemma}
	
	The following Theorem \ref{th: tensor block model} gives the theoretical guarantee on the performance of Algorithm \ref{alg: Commu Dete of tensor block model} for tensor reconstruction and cocluster membership recovery. 
	
	\begin{algorithm}[!h] \caption{HOOI for Tensor Co-clustering/Block Model} \label{alg: Commu Dete of tensor block model}
		\textbf{Input:} Tensor $\widetilde{\bcT}\in \bbR^{p_1 \times \cdots \times p_d}$, indices group $\Omega_i = \{i \}$, initialization $\widetilde{\U}^{(0)}_i \in \bbR^{p_i \times r_i}$ for $i=1, \ldots, d$, maximum number of iterations $t_{\max}$, \\
		\textbf{Output:} $\widehat{\bPi}_i \in \mathbb{M}_{p_i, r_i}$,$i =1, \ldots, d$ and $\widehat{\bcT}$. 
		\begin{algorithmic}[1]
			\State Apply Algorithm \ref{alg: HOOI} with input $\widetilde{\bcT}$, $\{\Omega_i\}_{i=1}^d$, $\{\widetilde{\U}^{(0)}_i\}_{i=1}^d$, maximum number of iterations $t_{\max}$ and get outputs $\{\widehat{\U}_i\}_{i=1}^d$ and $\widehat{\bcT}$.
			\State For each mode $i$, apply $\epsilon$-approximation K-means \citep{kumar2004simple} on $\widehat{\U}_i$, i.e., compute $\widehat{\bPi}_i, \in \mathbb{M}_{p_i, r_i}, \widehat{\X}_i \in \bbR^{r_i \times r_i}$ such that
			\begin{equation} \label{ineq: spectral clustering}
			\| \widehat{\bPi}_i \widehat{\X}_i - \widehat{\U}_i \|_F^2 \leq (1 + \epsilon) \min_{\bPi \in \mathbb{M}_{p_i, r_i}, \X \in \bbR^{r_i \times r_i}} \| \bPi \X - \widehat{\U}_i \|_F^2.
			\end{equation}
		\end{algorithmic}
	\end{algorithm}

	\begin{Theorem}[HOOI for Tensor Co-clustering/Block Model]\label{th: tensor block model}
		Consider the tensor co-clustering/block model ``$\eqref{eq: perturbation model}+\eqref{eq: tensor block model}$" and the Algorithm \ref{alg: Commu Dete of tensor block model} with inputs $\widetilde{\bcT}$, initializations $\{ \widetilde{\U}_i^{(0)} \}_{i=1}^d$ and $t_{\max} = C \left( \log \left(\frac{\lambda/\sigma}{\sqrt{p_{\max}}}\sqrt{\frac{\prod_{i=1}^d p_i }{\prod_{i=1}^d r_i}} \right) \vee 1  \right)$ for some $C > 0$, where $\lambda = \min_k \sigma_{r_k}\left(\mathcal{M}_k(\bcB)\right)$ is the minimal singular value at each matricization of the core tensor parameter $\bcB$. Assume $r_{\max} \leq p_{\min}^{1/(d-1)}$, $\max_{i}\| \widetilde{\U}^{(0)\top}_{i\perp} \U_i \| \leq \frac{\sqrt{2}}{2}$, and \eqref{assu: cocluster size assumption} holds. Then if
		\begin{equation*}
		\lambda/\sigma \geq C \sqrt{ \frac{p_{\max}^2 r_{\max} \prod_{i=1}^d r_i }{p_{\min} \prod_{i=1}^d p_i }},
		\end{equation*} 
		for sufficiently large constant $C > 0$, with probability at least $1 - \exp(-c p_{\min})$ for some $c > 0$, $\widehat{\U}_k$, $\widehat{\bcT}$ satisfy
		\begin{equation*}
		\| \sin \Theta(\widehat{\U}_k, \U_k)  \| \leq C(d)\frac{ \sqrt{p_k} }{\lambda/\sigma} \sqrt{\frac{\prod_{i=1}^d r_i}{\prod_{i=1}^d p_i}}, \quad k=1, \ldots, d
		\end{equation*}
		\begin{equation*}
		\|\widehat{\bcT} - \bcT\|_{\tHS} \leq C(d) \sigma \sqrt{ \sum_{i=1}^d p_i r_i };
		\end{equation*}
		and we also have the following upper bound on cocluster recovery error,
		\begin{equation*}
		\begin{split}
		& err(\widehat{\bPi}_i, \bPi_i) \leq C_1(d,\epsilon) \frac{ p'_{i,\max}r_{i}}{(\lambda/\sigma)^2} \frac{\prod_{i=1}^d r_i}{\prod_{i=1}^d p_i},\\
		& \widetilde{err}(\widehat{\bPi}_i, \bPi_i) \leq C_2(d,\epsilon) \frac{ p_i r_i }{(\lambda/\sigma)^2} \frac{\prod_{i=1}^d r_i}{\prod_{i=1}^d p_i}.
		\end{split}
		\end{equation*} 
		Here $C(d),C_1(d,\epsilon ), C_2(d,\epsilon )>0$ are some constants depending only on $d$ and $\epsilon$, $p'_{i,\max}$ is the second largest cocluster size at mode $i$.
	\end{Theorem}
	
	Note that our cocluster recovery guarantee is new for polynomial-time algorithms. When $p = p_1 = \cdots = p_d$, the best polynomial time algorithm guarantee for tensor reconstruction is $p^{d-1}$ in \cite{chi2020provable} and our result can be significantly better.

	\section{Numerical Studies} \label{sec: numerical study}
	
	In this section, we first provide numerical studies to support the main theoretical results in Section \ref{sec:main} and then compare HOOI with other existing algorithms for tensor decomposition in applications of tensor denoising and tensor co-clustering. Throughout the simulation, we consider order-$3$ tensor perturbation setting $\widetilde{\bcT} = \bcT + \bcZ$ with $\bcZ$ being the noise tensor with i.i.d. $N(0, \sigma^2)$ entries. Without particular specification, we set $p = p_1 = p_2 = p_3, r = r_1 = r_2 = r_3$. The error metrics we consider for tensor reconstruction and mode-$k$ singular subspace estimation are root mean square error (RMSE) $\| \widehat{\bcT} - \bcT \|_{\tHS}$ and $\|\sin \Theta(\widehat{\U}_k, \U_k) \|$, respectively. All simulations are repeated $100$ times and the average statistics are reported.
	
	\subsection{Perturbation Bounds of HOOI with good initialization} \label{sec: numerical study determinitic part}
	In this simulation, we study the perturbation bounds of HOOI with randomly generated good initialization. Let $\bcT = \bcS \times_1 \U_1 \times_2 \U_2 \times_3 \U_3$, where $\U_i \in \bbR^{p_i \times r}$ is generated uniformly at random from $\mathbb{O}_{p_i, r}$ and $\bcS \in \bbR^{r \times r \times r}$ is a diagonal tensor with diagonal values $\{ i\lambda\}_{i=1}^r$. The initializations of $\U_i$ of Algorithm \ref{alg: HOOI} are $\widetilde{\U}_i^{(0)} =  \frac{1}{\sqrt{2}} \U_i + \frac{1}{\sqrt{2}} \U_i'$, where $\U_i' = \U_{i\perp} \mathbf{O}$ for some random orthogonal matrix $\mathbf{O} \in \mathbb{O}_{p_i - r,r}$. It is easy to check that $\| \sin \Theta (\U_i, \widetilde{\U}_i^{(0)})\| = \frac{\sqrt{2}}{2}$ for $i = 1,2,3$.
	
\begin{figure}
		\centering
		\subfigure[Tensor Reconstruction]{\includegraphics[height = 0.25\textwidth]{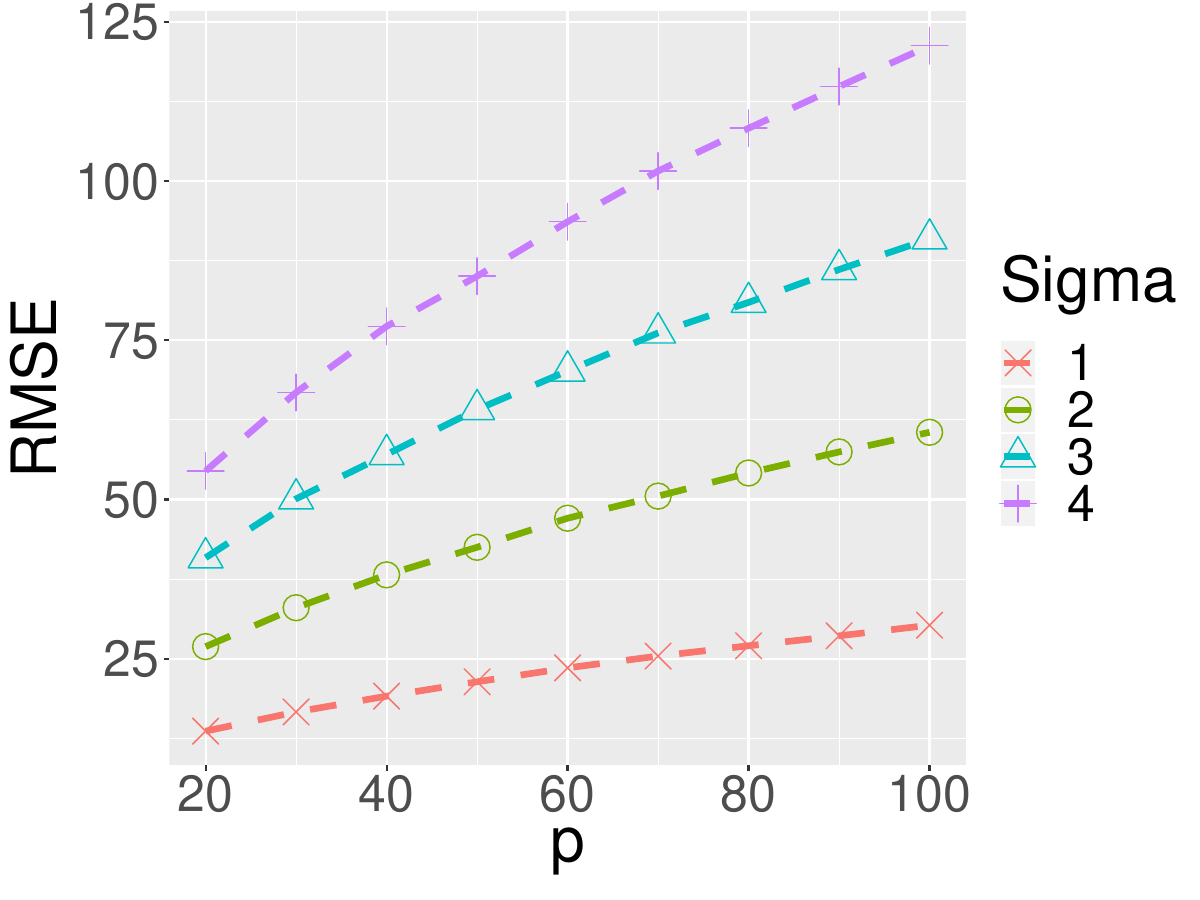}}
		\subfigure[Convergence of Singular Subspace]{\includegraphics[height = 0.25\textwidth]{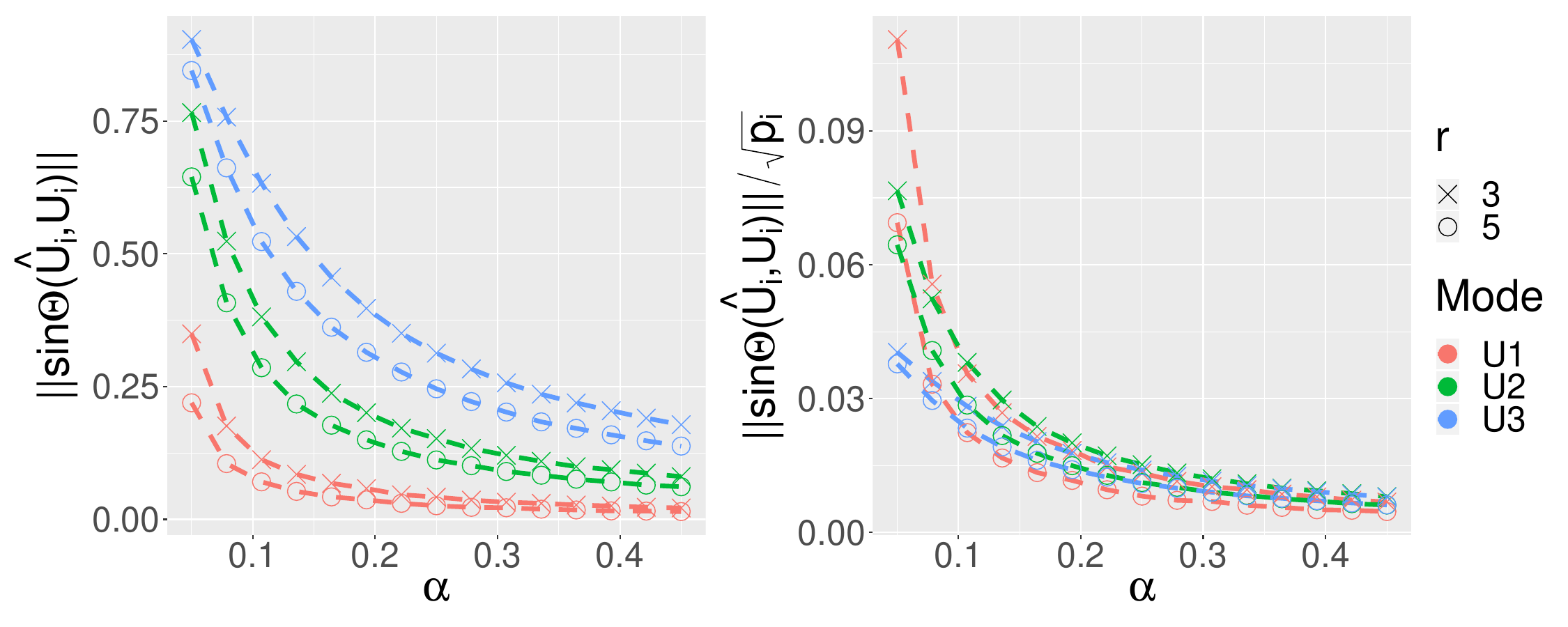}}
		\caption{HOOI with good initialization. (a) Tensor reconstruction error $\|\widehat{\bcT} - \bcT\|_{\tHS}$ for $p \in \{20,30,\ldots,100\}, r=5$, $\sigma \in \{1,2,3,4\}$ and $\lambda = 5 \sqrt{pr}\sigma$; (b) Mode-$k$ singular subspace estimation with and without rescaling under $p_1 = 10, p_2 = 100, p_3 = 500$, $r \in \{3,5\}, \sigma = 1$ and $\lambda = \alpha \cdot p_3 \frac{\sqrt{r_1}}{\sqrt{p_1}}$ with varying $\alpha$} \label{fig: simulation1}
	\end{figure}

	First for tensor reconstruction, let $p \in \{20,30,\ldots,100\}, r=5$, $\sigma \in \{1,2,3,4\}$ and $\lambda = 5 \sqrt{pr}\sigma$. We can check that with high probability, $\|\bcZ\|_{\tHS} \leq C p^{\frac{3}{2}} \sigma$ and $\xi \leq C\sqrt{pr}\sigma$ for some $C > 0$ following the same proof as Theorem \ref{th: tensor denoising}. 
	In Figure \ref{fig: simulation1}(a), the RMSE of tensor reconstruction of HOOI is presented. We find as the perturbation results in Section \ref{sec:main} suggest, $\|\widehat{\bcT}-\bcT\|_{\tHS}$ can be much smaller than $\|\bcZ\|_{\tHS}$. This demonstrates the superior performance of the HOOI estimator compared to the trivial estimator $\widetilde{\bcT}$. At the same time, the RMSE for tensor reconstruction increases as $p$ and $\sigma$ become bigger and this matches our theoretical findings in Theorem \ref{th:HOOI general} that the error bound of HOOI for $\| \widehat{\bcT} - \bcT \|_{\tHS}$ is $O(\xi)$, which increases as $p, \sigma$ increase.
	
	Next we demonstrate the unilateral perturbation bounds for mode-$k$ singular subspace estimation. Specifically, we consider $p_1 = 10, p_2 = 100, p_3 = 500$, $r \in \{ 3,5\}, \sigma = 1$ and $\lambda = \alpha \cdot p_3 \frac{\sqrt{r}}{\sqrt{p_1}}$ with varying $\alpha$. The errors of the mode-1, mode-2, mode-3 estimated singular subspaces with and without rescaling are provided in Figure \ref{fig: simulation1}(b). We can see from Figure \ref{fig: simulation1}(b) left panel the errors of estimated singular subspaces converge to different values depending on the corresponding mode size $p_i$, and a further rescaling of estimation error by $\sqrt{p_i}$ makes them roughly on the same level (see Figure \ref{fig: simulation1}(b) right panel). This matches the unilateral property of the singular subspace perturbation results in Remark \ref{rem: unilateral bound} that when $\lambda = O(p_3 \frac{\sqrt{r}}{\sqrt{p_1}})$, $\|\sin \Theta(\widehat{\U}_k, \U_k) \|\leq C\frac{\sqrt{p_i}}{\lambda/\sigma}$ for some $C > 0$, and this upper bound increases linearly with respect to $\sqrt{p_i}$.

	\subsection{Comparison of HOOI with other Algorithms} \label{sec: numerical study: comparison}
	In this section, we do a comparison of HOOI with truncated HOSVD (T-HOSVD) \citep{de2000multilinear} and sequentially truncated HOSVD (ST-HOSVD) \citep{vannieuwenhoven2012new} in the tensor denoising and tensor co-clustering applications. We also include one-step HOOI (O-HOOI), since it might be useful as a surrogate of HOOI in large scale tensor decomposition settings as we mentioned in Remark \ref{rem: one-step optimality of HOOI}. The initialization we consider for HOOI and O-HOOI are ST-HOSVD with natural truncation order, i.e., $\widetilde{\U}^{(0)}_i = \SVD_{r} (\cM_i(\widetilde{\bcT} \times_{j < i} \widetilde{\U}_j^{(0)} ))$. In Table \ref{tab: time complexity compare}, we give the time complexity of HOOI, O-HOOI, ST-HOSVD and T-HOSVD. We can see that as long as $dr \leq p$, a common case in practice, the time complexity of O-HOOI and ST-HOSVD are on the same order, and they could be faster than full HOOI and HOSVD in general.

	\begin{table}
		\centering
		\begin{tabular}{c | c | c | c |c}
			\hline
			&  HOOI & O-HOOI & ST-HOSVD & T-HOSVD \\
			\hline
			Complexity & $O(p^{d+1} + t_{\max} d r p^d )$ & $O(p^{d+1} + d r p^d)$ & $O(p^{d+1})$ & $O(d p^{d+1})$\\
			\hline
		\end{tabular}
		\caption{Time Complexity of HOOI, O-HOOI, ST-HOSVD and T-HOSVD under setting $p_1 = \cdots = p_d = p, r_1 = \cdots = r_d = r$, $r \ll p$. HOOI and O-HOOI is initialized by ST-HOSVD.} \label{tab: time complexity compare}
	\end{table}

\begin{figure}
		\centering
		\subfigure[Tensor Reconstruction]{\includegraphics[height = 0.26\textwidth]{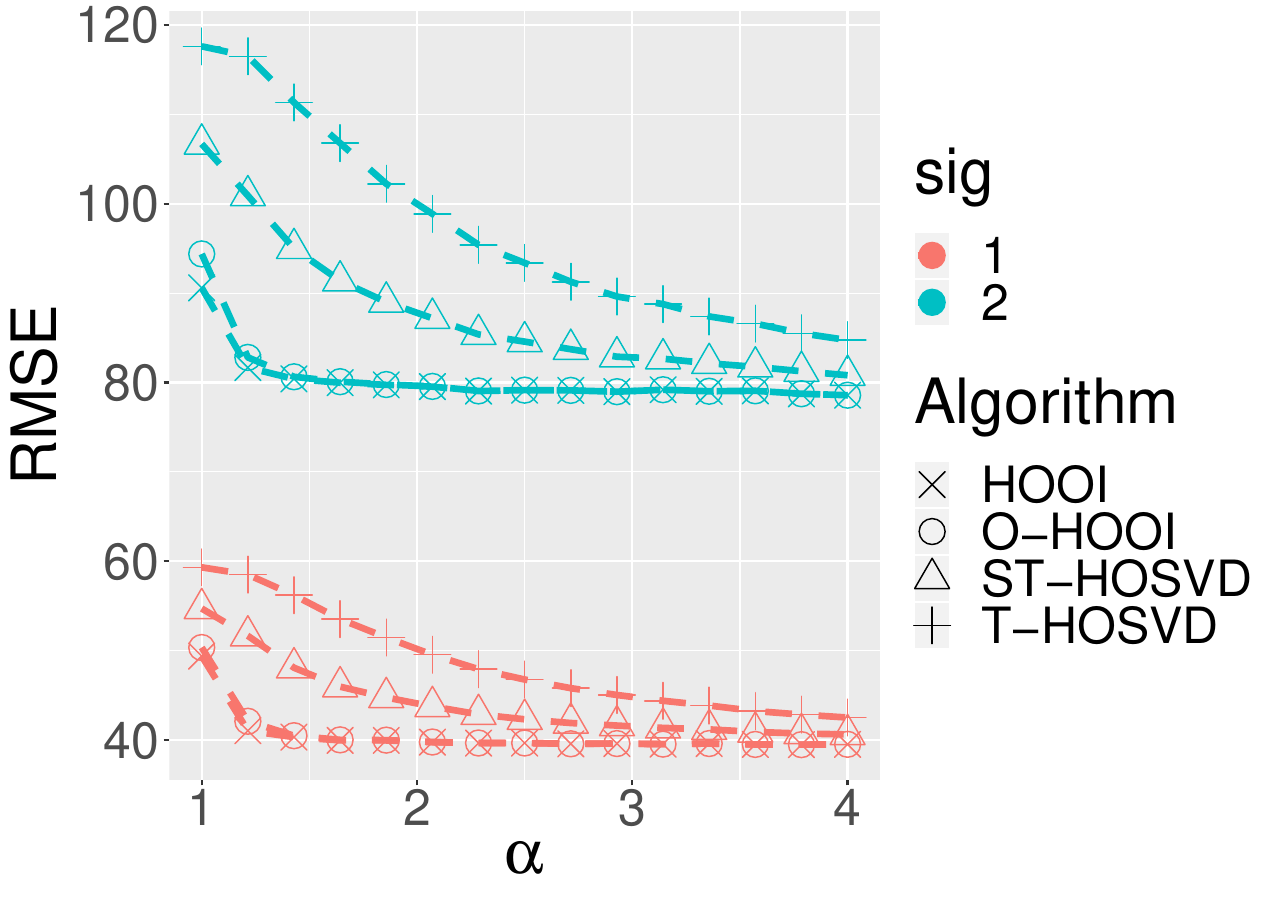}}
		\subfigure[Singular Subspace Estimation]{\includegraphics[height = 0.26\textwidth]{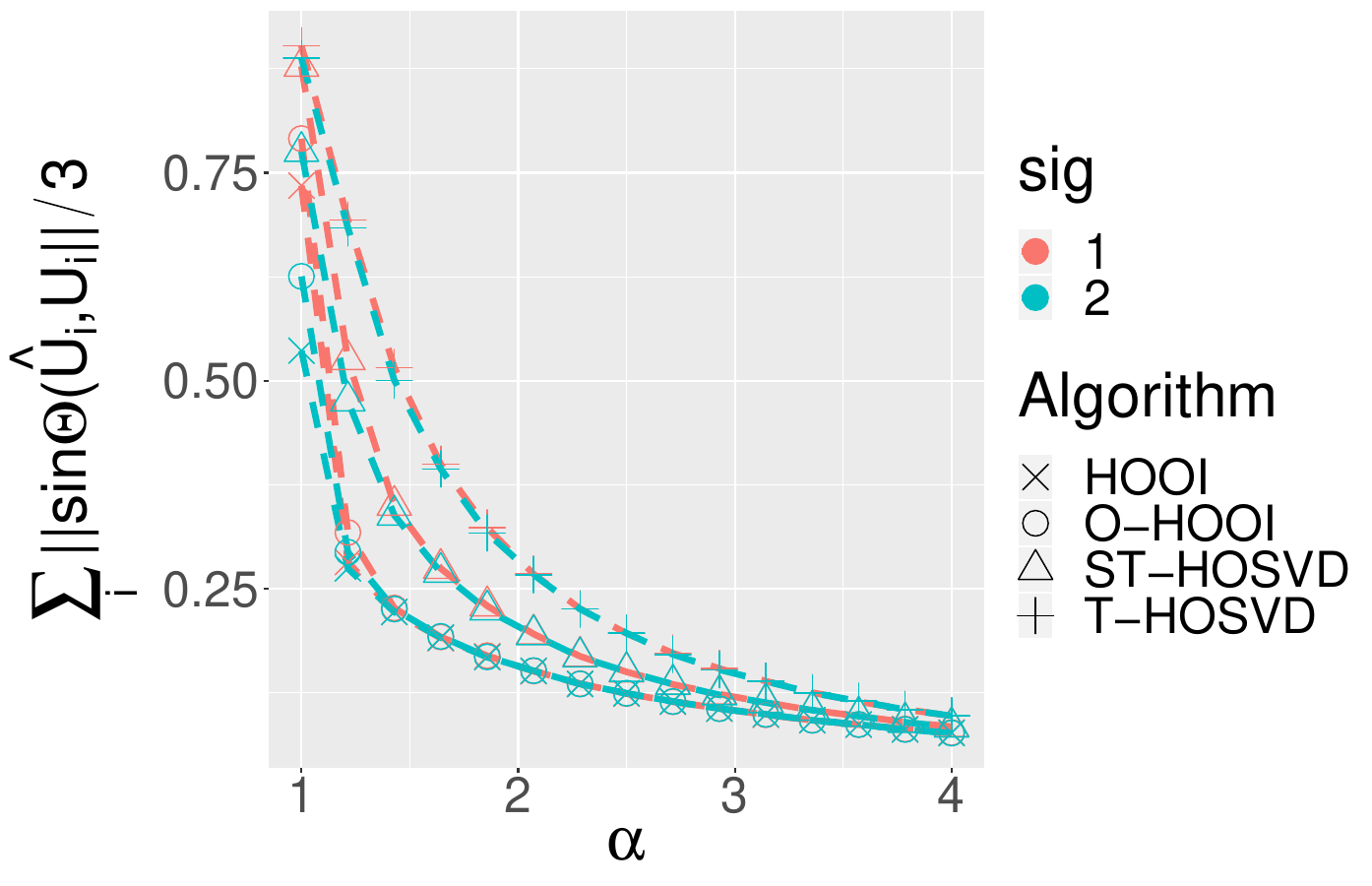}}
		\caption{Comparison of HOOI, one-step HOOI (O-HOOI), truncated HOSVD (T-HOSVD), sequentially truncated HOSVD (ST-HOSVD) in tensor denoising under $p = 100$, $r= 5$, $\sigma \in \{1,2 \}$, $\lambda = \alpha \cdot p^{\frac{3}{4}} \sigma$ with $\alpha \in [1,4]$. (a) Tensor reconstruction; (b) Averaged singular subspace estimation} \label{fig: simulation2}
	\end{figure}

	In tensor denoising, the generating process of $\bcT$ is the same as before. Let $p = 100$, $r= 5$, $\sigma \in \{1,2 \}$, $\lambda = \alpha \cdot p^{\frac{3}{4}} \sigma$ with varying $\alpha$. The comparison of these algorithms for tensor reconstruction and singular subspace estimation are given in Figure \ref{fig: simulation2}. First, we find that HOOI is best in both tensor reconstruction and singular subspace estimation among four algorithms. Meanwhile, O-HOOI is slightly worse than HOOI for small $\alpha$ and has very close performance with HOOI when $\alpha$ is relative large, which suggests that in some computationally heavy applications, we can just run HOOI for one iteration to achieve reasonable estimation. Part of this phenomenon can be explained by the one-iteration optimality of HOOI for tensor reconstruction as we discussed in Remark \ref{rem: one-step optimality of HOOI}. 
	On the other hand, HOOI and O-HOOI are often much better than T-HOSVD and ST-HOSVD for both tensor reconstruction and mode-$k$ singular subspace estimation within a wide range of settings. 
	
\begin{figure}
		\centering
		\subfigure[HOOI for Cocluster Recovery]{\includegraphics[height = 0.29\textwidth]{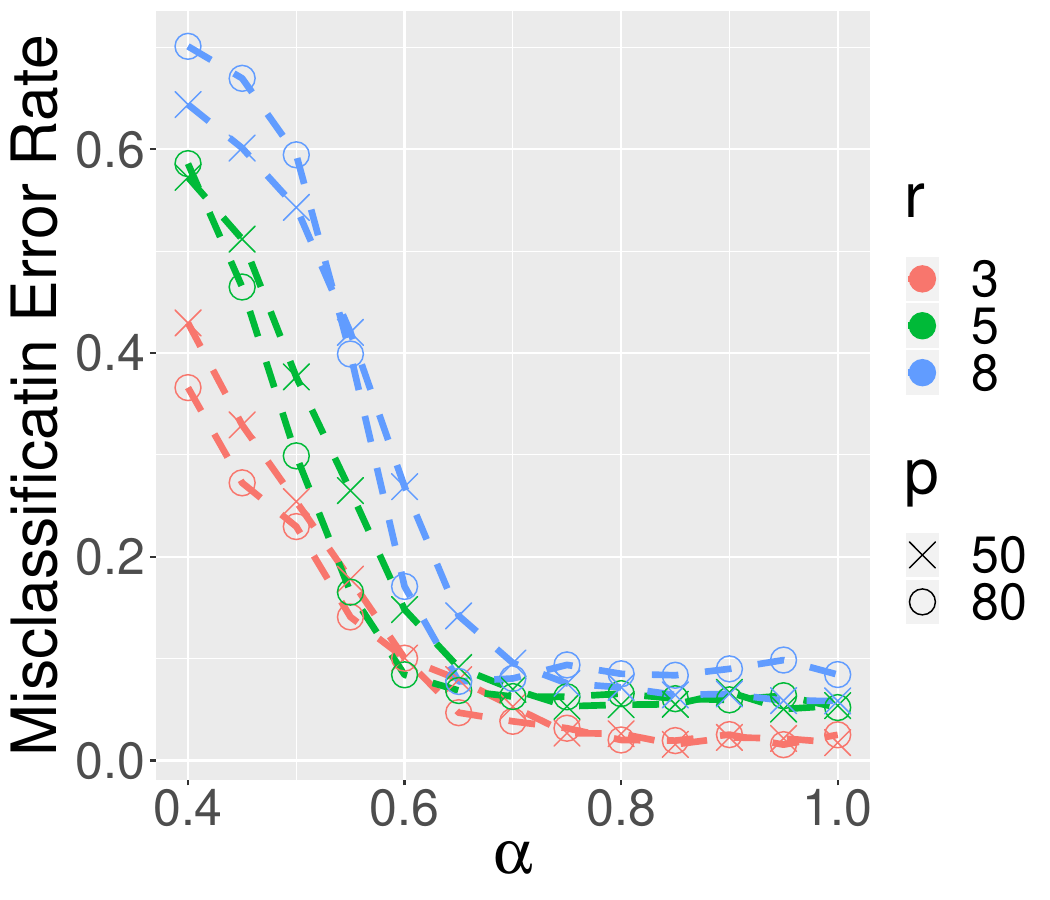}}
		\hskip.5cm
		\subfigure[Comparison in Cocluster Recovery]{\includegraphics[height = 0.29\textwidth]{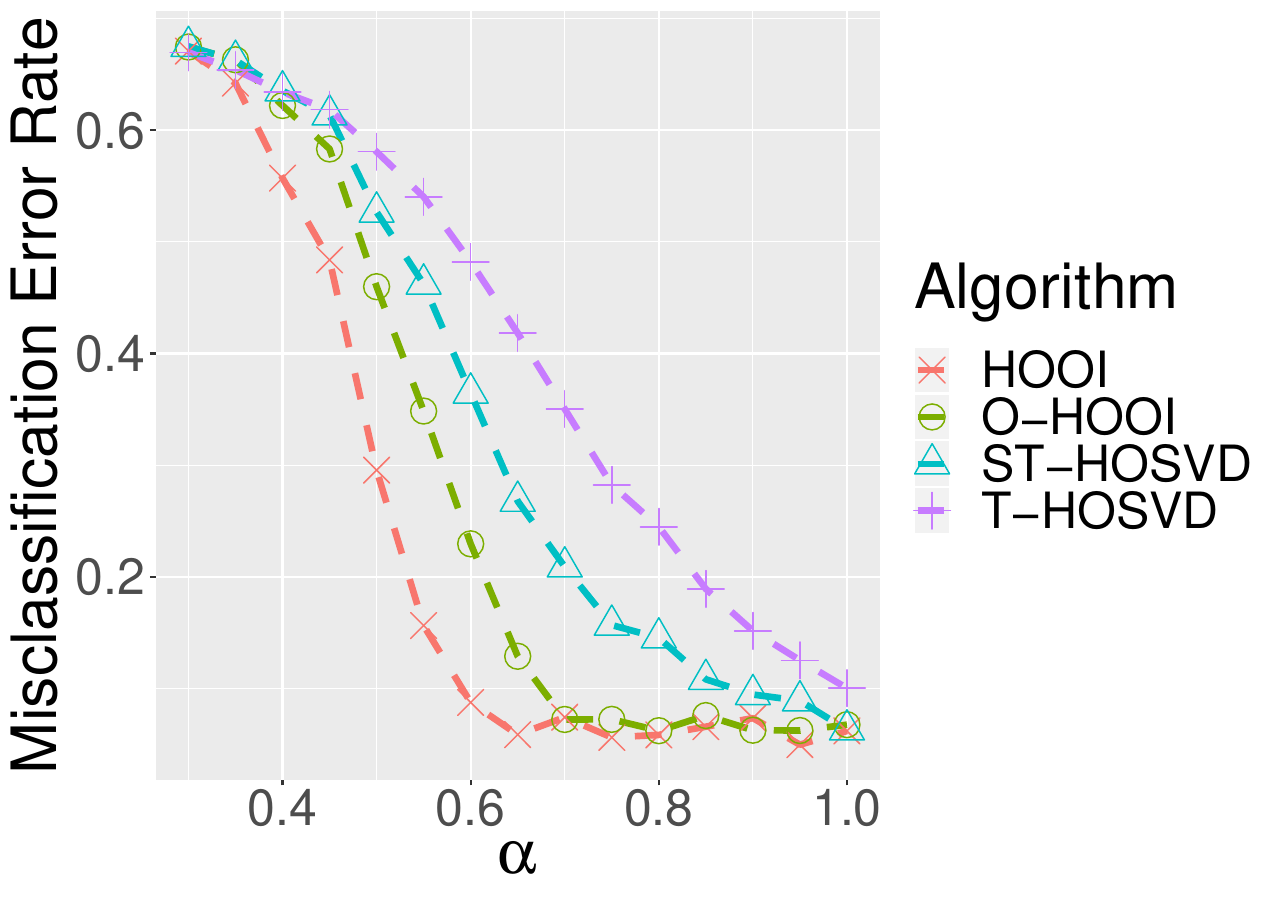}}
		\caption{Tensor cocluster recovery under $\sigma = 1$, $\lambda = \alpha \cdot \frac{r^{3/2}}{p^{3/4}} \sigma$ with varying $\alpha$. (a) HOOI on tensor cocluster recovery under $p \in \{50,80 \}$, $r \in \{3,5,8 \}$. (b) Comparison of HOOI, one-step HOOI (O-HOOI), truncated HOSVD (T-HOSVD), sequentially truncated HOSVD (ST-HOSVD) in cocluster recovery under $p = 80, r = 5$.} \label{fig: simulation3}
	\end{figure}

	Finally, we study the performance of the HOOI-based Algorithm \ref{alg: Commu Dete of tensor block model} in tensor co-clustering recovery and do a comparison of it with T-HOSVD, ST-HOSVD, and O-HOOI based clustering algorithms. In this simulation, we generate $\bcT = \bcB \times_1 \bPi_1 \times_2 \bPi_2 \times_3 \bPi_3$ such that $\{\bPi_i\}_{i=1}^3$ have balanced cluster size and $\bcB = \frac{\bcB_0}{\min_i \sigma_r\left(\cM_i(\bcB_0) \right)} \lambda $ with $\bcB_0 \overset{i.i.d.}\sim N(0,1)$. The error metric we consider is the average cocluster membership misclassification error rate in \eqref{eq: misclassification error rate}. The performance of Algorithm \ref{alg: Commu Dete of tensor block model} under $p \in \{50,80 \}, r \in \{3,5,8 \}, \sigma = 1, \lambda = \alpha \cdot \frac{r^{3/2}}{p^{3/4}} \sigma$ is presented in Figure \ref{fig: simulation3}(a). We can see the misclassification error decreases as the signal strength increases and cocluster number decreases. The comparison of Algorithm \ref{alg: Commu Dete of tensor block model} and T-HOSVD, ST-HOSVD, O-HOOI based spectral clustering is given in Figure \ref{fig: simulation3}(b) under the same setting with $r = 5$. Again, HOOI-based algorithm has the best performance in cocluster recovery. O-HOOI and ST-HOSVD perform similarly here and both of them are much better than T-HOSVD.
	
	\section{Conclusion and Discussion}\label{sec: conclusion and discussion}
	
	In this paper, we provide the first sharp blockwise perturbation bounds of HOOI for tensors with guarantees for both tensor reconstruction and mode-$k$ singular subspace estimation. Furthermore, we show both HOOI and one-step HOOI with good initialization is optimal in terms of tensor reconstruction by providing rate matching lower bound. Finally, we support our theoretical results with extensive numerical studies and apply them to tensor denoising and tensor co-clustering applications. Apart from the applications mentioned above, the main perturbation results can be applied to many other applications where tensor ``spectral method" HOOI is applicable, such as tensor completion \citep{yuan2016tensor,yuan2017incoherent,xia2017statistically,xia2017polynomial}, hypergraphic stochastic block model \citep{ghoshdastidar2014consistency,ghoshdastidar2017consistency,ke2019community,chien2019minimax,ahn2018hypergraph,kim2018stochastic}, multilayer network \citep{lei2019consistent,jing2020community}, MPCA \citep{lu2008mpca}, latent variable model \citep{anandkumar2014tensor}, etc. In tensor completion \citep{xia2017statistically} and many other applications, more specialized initializers can achieve better performance than HOSVD -- the classic initializers for HOOI in the literature \citep{de2000best}. Our tensor perturbation bounds still apply to these cases as our theoretical analysis admits all initializers satisfying certain mild conditions.
	
	At the same time, due to the NP hardness of computing many tensor quantities \citep{hillar2013most}, the Alternating Least Square (ALS) and Power iteration have been the ``workhorse" algorithms in computing low-rank tensor approximation and solving many other tensor problems \citep{kolda2009tensor}. Our induction proof idea in Theorem \ref{th:HOOI general} could also shed light on how to analyze other iterative ALS/Power iteration procedures for tensor problems \citep{zhou2013tensor,lu2008mpca,wang2019multiway,xu2005concurrent,yan2005discriminant,yan2006multilinear,lee2020tensor}.
	
	The convergence rate of HOOI is another important topic related to the results in this paper. \cite{savas2010quasi,elden2009newton} observed that HOOI converge fast when the target tensor has fast-decaying multilinear singular values (in other words, the tensor is approximately Tucker low-rank). On the other hand, HOOI may converge slowly when the target tensor is high-rank or approximately sparse \citep{savas2010quasi,elden2009newton}. In this paper, we studied how HOOI converges when the target tensor is approximately Tucker low-rank. It is interesting to further explore how HOOI converges in the less ideal settings, such as for the high-rank or approximately sparse tensors. 

	In addition to the perturbation results for HOOI, it is also interesting to develop perturbation results for the second-order algorithms, such as (quasi-)Newton-Grassmann method \citep{elden2009newton,savas2010quasi}, geometric Newton method \citep{ishteva2009differential}, Riemannian trust region scheme \citep{ishteva2011best} since these second-order methods may take much fewer iterations than HOOI to converge. Also, this paper mainly focuses on Tucker format of tensor decomposition. Although Tucker format has many advantages, in ultra higher-order tensor problems, the storage cost of the core tensor in Tucker format scales exponentially with respect to the tensor order and it is more desirable to consider other low-rank tensor decomposition formats, such as the hierarchical Tucker decomposition \citep{ballani2013projection,grasedyck2010hierarchical,hackbusch2009new} and Tensor-Train decomposition \citep{oseledets2011tensor, oseledets2009breaking,zhou2020optimal}. It is interesting to develop the perturbation bounds for algorithms on Hierarchical Tucker or Tensor Train tensor decompositions.
	
\acks{The authors would like to thank the Action Editor Animashree Anandkumar and the anonymous referees for constructive suggestions. Y.L. would like to thank Rungang Han for helpful discussions during the project. The research of Y. L. and A. R. Zhang was in part by NSF Grants CAREER-1944904, NSF DMS-1811868, NSF DMS-2023239, NIH Grant R01 GM131399, and Wisconsin Alumni Research Foundation (WARF); the research of G. R. was supported in part by NSF Grant DMS-2015285; the research of M. Y. was supported in part by NSF Grant DMS-1811767.}

\appendix



\section{Tensor Perturbation Bounds for HOOI in Asymmetric Case} \label{sec: HOOI in asymmetric case}
In this section, we present a Corollary of Theorem \ref{th:HOOI general} in the case where $\Omega_i = \{ i \}, i = 1, \ldots, d$, which appears frequently in practice.
\begin{Corollary}[Tensor Perturbation Bounds for HOOI in Asymmetric Case]\label{th:HOOI d diff}
	Consider the perturbation model \eqref{eq: perturbation model} with $\widetilde{\bcT}, \bcT, \bcZ \in \mathbb{R}^{p_1\times \cdots \times p_d}$ and $\Omega_i = \{ i\}$, $i = 1, \ldots, d$.
	Define  $\cS_{i}^{(-k)} := \left\{ S \subseteq [d]\setminus \{k\}: |S| = i \right\} $ as the set of all possible index sets with $i$ elements from $[d]\setminus \{k \}$ and $\cS_0^{(-k)}:= \emptyset $. For $S \in \cS_{i}^{(-k)}$, let $S^c = ([d]\setminus \{k\}) \setminus S$. Now we define the blockwise errors as
	\begin{equation*}
	\begin{split}
	&  \tau_1 = \max_{k=1,\ldots, d} \tau_{1k}, \quad \tau_{1k} = \left\|  \left(\mathcal{M}_k(\bcZ \times_{i \neq k}\U_{i}^\top ) \right)_{\max(r_k)}\right\|_q, \quad k=1,\ldots, d;\\
	& \tau_j = \max_{k = 1, \ldots,d}\Big\{ \max_{S \in \cS^{(-k)}_{j-1}} \sup_{\substack{  \V_i \in \bbR^{(p_i - r_i) \times r_i},\\ \|\V_i\|_q \leq 1, i \in S }} \left\|\left(\mathcal{M}_k(\bcZ \times_{i \in S}\left(\U_{i\perp} \V_i \right)^\top \times_{i \in S^c}\U_{i}^\top \right)_{\max(r_k)}\right\|_q \Big\},\\
	&\quad \text{ for } j = 2, \ldots, d;
	\end{split}
	\end{equation*}
	Denote the initialization errors of $\{ \widetilde{\U}_k^{(0)} \}_{k=1}^d$ as $\tilde{e}_0 := \max_{k=1, \ldots, d}\|\widetilde{\U}^{(0)\top}_{k\perp} \U_k\|$, $e_0 := \max_{k=1, \ldots, d}\|\widetilde{\U}^{(0)\top}_{k\perp} \U_k\|_q$. Assume the initialization error and the signal strength satisfy 
	\begin{equation*}
	\tilde{e}_0 \leq \frac{\sqrt{2}}{2} \, \text{ and } \, \lambda \geq 2^{\frac{d+4}{2}} \left(1+\frac{\sqrt{2}}{2}\right)^d \xi,
	\end{equation*}
	where $\xi := \underset{\|\bcY \|_{\tHS} \leq 1,  \rank(\bcY) \leq (r_{1}, \ldots, r_{d})}\sup \langle \bcZ, \bcY \rangle$. 
	
	Then with inputs $\widetilde{\bcT}$, $\{\widetilde{\U}_k^{(0)} \}_{k=1}^d$, $\{\Omega_i \}_{i=1}^d$, the estimated mode-$k$ singular subspaces updates in Algorithm \ref{alg: HOOI} after $t$ iterations satisfy
	\begin{equation*}
	\max_{k \in [d]} \left\|\sin\Theta\left(\widetilde{\U}_k^{(t)}, \U_k\right)\right\|_q \leq 2^{\frac{d+3}{2}} \frac{\tau_1}{\lambda} + \frac{e_0}{2^{t}}.
	\end{equation*}
	Moreover, when $t_{\max} \geq \log(e_0\lambda/\tau_{1})\vee 1$, the outputs of estimated mode-$k$ singular subspaces of Algorithm \ref{alg: HOOI} satisfy
	\begin{equation*}
	\max_{k \in [d]} \left\|\sin\Theta\left(\widehat{\U}_k, \U_k\right)\right\|_q \leq \left(2^{\frac{d+3}{2}} + 1\right) \frac{\tau_1}{\lambda},
	\end{equation*}
	\begin{equation*}
	\begin{split}
	\left\|\sin\Theta\left(\widehat{\U}_k, \U_k\right)\right\|_q\leq \frac{2}{\left(1- c^*(\tau_1,\lambda,d) \right)^{\frac{d-1}{2} }} \left(\frac{\tau_{1k}}{\lambda} +  \sum_{j=1}^{d-1} \frac{{d-1 \choose j}  \left(2^{\frac{d+3}{2}}+1\right)^j \tau_1^j  \tau_{j+1} }{\lambda^{j+1}} \right), 
	\end{split}
	\end{equation*}
	for $k=1,\ldots,d$, where $c^*(\tau_1,\lambda,d) = \left(2^{\frac{d+3}{2}} + 2\right)^2 \frac{\tau_1^2}{\lambda^2} \leq \frac{1}{2}$, and the output of tensor reconstruction $\widehat{\bcT}$ satisfies
	\begin{equation*}
	\begin{split}
	\left\|\widehat{\bcT} - \bcT\right\|_{\tHS} & \leq \left\|\llbracket \bcZ; \widehat{\U}_1^\top, \ldots, \widehat{\U}_d^\top \rrbracket\right\|_{\tHS} + \sum_{k=1}^d \left\|\widehat{\U}_{k\perp}^\top \mathcal{M}_k(\bcT)\right\|_F \\
	& \leq \left(1 + 2d \left( 1- c^*(\tau_1,\lambda,d) \right)^{- \frac{d-1}{2} } \right) \xi.
	\end{split}
	\end{equation*}
\end{Corollary}

\section{Additional Proofs} \label{sec: additional proofs}

\subsection{Proof of Lemma \ref{lem: character of xi}}
For generality, here we present the proof for order $d$ case. First we show the first equivalent characterization.
\begin{equation*}
\begin{split}
\underset{\|\bcY \|_{\tHS} \leq 1,  \rank(\bcY) \leq (r_1, \ldots, r_d)}\sup \langle \bcZ, \bcY \rangle &\geq \underset{ \| \bcS \|_{\tHS} \leq 1, \U_i \in \mathbb{O}_{p_i, r_i} } \sup \langle \bcZ, \bcS \times_1 \U_1 \times \cdots \times_d \U_d \rangle \\
& = \underset{ \| \bcS \|_{\tHS} \leq 1, \U_i \in \mathbb{O}_{p_i, r_i} } \sup \langle \bcZ \times_1 \U_1^\top \times \cdots \times \U_d^\top, \bcS \rangle\\
& = \underset{\U_i \in \mathbb{O}_{p_i, r_i} } \sup \| \bcZ \times_1 \U_1^\top \times \cdots \times_d \U_d^\top \|_{\tHS}.
\end{split}
\end{equation*}
On the other hand by Theorem 2 of \cite{de2000multilinear}, we have 
\begin{equation*}
\bcY = \left(\bcY \times_1 \widetilde{\U}_1^\top \times \cdots \times \widetilde{\U}_d^\top \right) \times_1 \widetilde{\U}_1 \times \cdots \times_d \widetilde{\U}_d,  
\end{equation*} where $\widetilde{\U}_i \in \bbR^{p_i \times r_i}$ is the left singular space of $\mathcal{M}_i(\bcY)$. Since $\|\bcY \times_1 \widetilde{\U}_1^\top \times \cdots \times \widetilde{\U}_d^\top\|_{\tHS} \leq 1$, 
\begin{equation*}
\underset{\|\bcY \|_{\tHS} \leq 1,  \rank(\bcY) \leq (r_1, \ldots, r_d)}\sup \langle \bcZ, \bcY \rangle \leq \underset{ \| \bcS \|_{\tHS} \leq 1, \U_i \in \mathbb{O}_{p_i, r_i} } \sup \langle \bcZ, \bcS \times_1 \U_1 \times \cdots \times_d \U_d \rangle.
\end{equation*}
So, we have proved the lemma. 
\quad $\blacksquare$

\subsection{Proof of Lemma \ref{lm: noise projection error}}
This proof idea of this lemma is to project $\bcZ$ onto orthogonal subspaces $\U_k$ and $\U_{k \perp}$ at each modes.

\begin{equation*}
\begin{split}
&\bcZ \times_1 \widehat{\U}_1^\top \times \cdots \times \widehat{\U}_d^\top \\
=  & \left( \bcZ \times_1 (\Proj_{\U_1} + \Proj_{\U_{1\perp}}) \times \cdots \times_d (\Proj_{\U_d} + \Proj_{\U_{d\perp}}) \right)  \times_1 \widehat{\U}_1^\top \times \cdots \times_d \widehat{\U}_d^\top\\
=  &\left( \sum_{\Omega \subseteq [d] } \bcZ \times_{k \in \Omega} \Proj_{\U_k} \times_{k \in \Omega^c} \Proj_{\U_{k\perp}}  \right)\times_1 \widehat{\U}_1^\top \times \cdots \times_d \widehat{\U}_d^\top\\
=  & \sum_{\Omega \subseteq [d] } \bcZ \times_{k \in \Omega} \widetilde{\U}_k^\top \Proj_{\U_k} \times_{k \in \Omega^c} \widetilde{\U}_k^\top \Proj_{\U_{k\perp}}. 
\end{split}
\end{equation*}
So by triangle inequality, we have 
\begin{equation*}
\begin{split}
\left\| \llbracket \bcZ; \widehat{\U}_1^\top, \ldots, \widehat{\U}_d^\top \rrbracket  \right\|_{\tHS} & \leq \sum_{\Omega \subseteq [d] } \left\| \bcZ \times_{k \in \Omega} \widetilde{\U}_k^\top \Proj_{\U_k} \times_{k \in \Omega^c} \widetilde{\U}_k^\top \Proj_{\U_{k\perp}} \right\|_{\tHS}\\
& \leq \sum_{\Omega \subseteq [d] } \left\| \bcZ \times_{k \in \Omega} \U_k^\top \times_{k \in \Omega^c} \U_{k\perp}^\top \right\|_{\tHS} \prod_{k \in \Omega} \| \widetilde{\U}_k^\top \U_k \| \prod_{k \in \Omega^c} \| \widetilde{\U}_k^\top \U_{k \perp} \|\\
& \leq \sum_{\Omega \subseteq [d] } \left\| \bcZ \times_{k \in \Omega} \U_k^\top \times_{k \in \Omega^c} \U_{k\perp}^\top \right\|_{\tHS} \prod_{k \in \Omega^c} \| \sin \Theta(\widehat{\U}_k, \U_k) \|\\
& =  \sum_{\Omega \subseteq \{1,\ldots, d\}} \theta_{\Omega} \prod_{k\in \Omega^c}\left\|\sin\Theta(\widehat{\U}_k, \U_k)\right\|.
\end{split}
\end{equation*}
Here the second inequality is due to the fact that $\| \bcZ \times_i \A \B \|_{\tHS} =  \| \A \B \mathcal{M}_i (\bcZ) \|_{\tHS} \leq \|\A\| \| \B \mathcal{M}_i (\bcZ) \|_{\tHS} = \|\A\| \| \bcZ \times_i \B \|_{\tHS}$ and we apply iteratively for each mode with $\A = \widetilde{\U}_k^\top \U_k$ (or $\widetilde{\U}_k^\top \U_{k\perp}$) and $\B = \U_k^\top$ (or $\U_{k\perp}^\top$). The third inequality is due to that $\| \widetilde{\U}_k^\top \U_k \| \leq 1$ and $\left\| \widetilde{\U}_k^\top \U_{k \perp} \right\| = \left\|\sin\Theta(\widehat{\U}_k, \U_k)\right\|$. \quad $\blacksquare$

\subsection{Proof of Theorem \ref{th: pertur lower bound}}
The proof is done by construction. Let's denote $\bcI_r \in (\bbR^r)^{\otimes d}$ as the order-$d$ identity tensor with entries $(i,i,\ldots, i)$ to be $1$ and others are $0$. We construct 
\begin{equation*}
\bcZ_1 = \frac{\xi}{\sqrt{r}} \bcI_r \times_1 \left(  \begin{array}{c}
\0_{r \times r}  \\
\I_r \\
\0_{(p_1-2r) \times r}
\end{array} \right) \times \cdots \times_d \left(  \begin{array}{c}
\0_{r \times r}  \\
\I_r \\
\0_{(p_d-2r) \times r}
\end{array} \right),
\end{equation*}where $\0_{m \times n}$ denotes a $m \times n$ matrix with all entries to be $0$. 

It is easy to check that $\underset{\|\bcY \|_{\tHS} \leq 1,  \rank(\bcY) \leq (r, \ldots, r)}\sup \langle \bcZ_1, \bcY \rangle \leq \|\bcZ_1\|_{\tHS} = \xi $. Similarly we construct 
\begin{equation*}
\bcT_1 = \frac{\xi}{\sqrt{r}} \bcI_r \times_1 \left(  \begin{array}{c}
\I_r \\
\0_{r \times r} \\
\0_{(p_1-2r) \times r}
\end{array} \right) \times \cdots \times_d \left(  \begin{array}{c}
\I_r  \\
\0_{r \times r} \\
\0_{(p_d-2r) \times r}
\end{array} \right).
\end{equation*}

Also we let $\bcZ_2 = \bcT_1$ and $\bcT_2 = \bcZ_1$, and it is easy to check $(\bcT_1, \bcZ_1), (\bcT_2, \bcZ_2) \in \mathcal{F}_r(\xi)$. At the same time, we have $\bcZ_1 + \bcT_1 = \bcT_2 + \bcZ_2$. Thus
\begin{equation*}
\begin{split}
\inf_{\widehat{\bcT}} \sup_{(\bcT, \bcZ) \in \mathcal{F}_r(\xi)} \| \widehat{\bcT} - \bcT \|_{\tHS} &\geq \inf_{\widehat{\bcT}} \max \left\{ \| \widehat{\bcT} - \bcT_1 \|_{\tHS}, \| \widehat{\bcT} - \bcT_2 \|_{\tHS}  \right\} \\
& \geq \inf_{\widehat{\bcT}}\frac{1}{2} \left(\| \widehat{\bcT} - \bcT_1 \|_{\tHS} + \| \widehat{\bcT} - \bcT_2 \|_{\tHS}\right)\\
& \geq \frac{1}{2} \|\bcT_1 - \bcT_2\|_{\tHS} = \frac{\sqrt{2}}{2} \xi.  
\end{split}
\end{equation*}
\quad $\blacksquare$

\subsection{Proof of Theorem \ref{th:HOOI general}}
The proof is long and nontrivial. The main idea of the proof is to develop the recursive error bound of $\widetilde{\U}_k^{(t+1)}$ i.e., the estimate of $\U_k$ at iteration $t+1$, based on the error bound of $\widetilde{\U}_k^{(t)}$, i.e., the estimate at iteration $t$. The outline of the proof is the following: after a briefly introduction of notations, the main proof could be divided into three steps. 

\begin{itemize}
	\item Step 1: Recall in HOOI procedure, the update for the mode-$k$ singular subspace satisfies
	\begin{equation*}
	\begin{split}
	\widetilde{\U}_k^{(t+1)} 
	= & \SVD_{r_k}\bigg( \mathcal{M}_{\bar{k}}\left(\bcT \times_{i \in \underline{\Omega}_k} \widetilde{\U}^{(t+1)\top}_{i'} \times_{\widecheck{\Omega}_k} \widetilde{\U}^{(t)\top}_k  \times_{i \in \overline{\Omega}_k} \widetilde{\U}^{(t)\top}_{i'}   \right) \\
	& \quad \quad \quad + \mathcal{M}_{\bar{k}}\left( \bcZ \times_{i \in \underline{\Omega}_k} \widetilde{\U}^{(t+1)\top}_{i'} \times_{\widecheck{\Omega}_k} \widetilde{\U}^{(t)\top}_k  \times_{i \in \overline{\Omega}_k} \widetilde{\U}^{(t)\top}_{i'}  \right)\bigg).
	\end{split}
	\end{equation*}
	To give an upper bound for $e_{t+1,k}$, we aim to give an upper bound for 
	\begin{equation*} 
	\left\|\left(\mathcal{M}_{\bar{k}}\left(\bcZ \times_{i \in \underline{\Omega}_k} \widetilde{\U}^{(t+1)\top}_{i'} \times_{\widecheck{\Omega}_k} \widetilde{\U}^{(t)\top}_k  \times_{i \in \overline{\Omega}_k} \widetilde{\U}^{(t)\top}_{i'}  \right) \right)_{\max(r_k)}\right\|_q
	\end{equation*}
	by using $\tau_1, \ldots, \tau_m$, $e_t, e_{t+1}$ in this step..
	\item Step 2: After getting an upper bound for \eqref{eq: step1 upper bound}, in this step we use induction to prove the following claim,
	\begin{equation*}
	\begin{split}
		e_{t} \leq 2^{(d+3)/2}\tau_{1}/\lambda + e_0/2^t; t=0,1,\ldots. \\
		\tilde{e}_{t} \leq 2^{(d+3)/2}\tau_{1}/\lambda + \tilde{e}_0/2^t; t=0,1,\ldots. 
	\end{split}
	\end{equation*}
	\item Step 3: Derive the error bound for $\|\widehat{\bcT} - \bcT\|_{\tHS}$ by the unified quantity $\xi$.
\end{itemize}

For convenience, in this proof we denote
$$\T_k = \mathcal{M}_{\bar{k}}(\bcT), \quad \widetilde{\T}_k = \mathcal{M}_{\bar{k}}(\widetilde{\bcT}), \quad \Z_k = \mathcal{M}_{\bar{k}}(\bcZ),\quad k=1,\ldots,m.$$ 
Suppose 
\begin{equation}\label{ineq:e_t-f_t}
\begin{split}
	e_{t} = \max_k e_{t,k},\quad e_{t,k} = \left\|(\widetilde{\U}_{k\perp}^{(t)})^\top \U_k \right\|_q,\quad k=1,\ldots,m; t=0,1,\ldots.\\
	\tilde{e}_{t} = \max_k \tilde{e}_{t,k},\quad \tilde{e}_{t,k} = \left\|(\widetilde{\U}_{k\perp}^{(t)})^\top \U_k \right\|,\quad k=1,\ldots,m; t=0,1,\ldots.
\end{split}
\end{equation}

{\noindent \bf Step 1}. Recall the procedure of HOOI that
\begin{equation}\label{ineq:HOOI-0}
\begin{split}
\widetilde{\U}_k^{(t+1)} = & \SVD_{r_k}\left(\mathcal{M}_{\bar{k}}\left(\widetilde{\bcT}\times_{i \in \underline{\Omega}_k} \widetilde{\U}^{(t+1)\top}_{i'} \times_{i \in \widecheck{\Omega}_k} \widetilde{\U}^{(t)\top}_k  \times_{i \in \overline{\Omega}_k} \widetilde{\U}^{(t)\top}_{i'}   \right) \right) \\
= & \SVD_{r_k}\bigg( \mathcal{M}_{\bar{k}}\left(\bcT \times_{i \in \underline{\Omega}_k} \widetilde{\U}^{(t+1)\top}_{i'} \times_{\widecheck{\Omega}_k} \widetilde{\U}^{(t)\top}_k  \times_{i \in \overline{\Omega}_k} \widetilde{\U}^{(t)\top}_{i'}   \right) \\
& + \mathcal{M}_{\bar{k}}\left( \bcZ \times_{i \in \underline{\Omega}_k} \widetilde{\U}^{(t+1)\top}_{i'} \times_{\widecheck{\Omega}_k} \widetilde{\U}^{(t)\top}_k  \times_{i \in \overline{\Omega}_k} \widetilde{\U}^{(t)\top}_{i'}  \right)\bigg).
\end{split}
\end{equation}
Notice that $\rank \left( \mathcal{M}_{\bar{k}}\left(\bcT \times_{i \in \underline{\Omega}_k} \widetilde{\U}^{(t+1)\top}_{i'} \times_{\widecheck{\Omega}_k} \widetilde{\U}^{(t)\top}_k  \times_{i \in \overline{\Omega}_k} \widetilde{\U}^{(t)\top}_{i'}   \right) \right) \leq r_k$, to apply Theorem 5 in \cite{luo2020schatten}, the key is to give an upper bound for $$\left\|\left(\mathcal{M}_{\bar{k}}\left(\bcZ \times_{i \in \underline{\Omega}_k} \widetilde{\U}^{(t+1)\top}_{i'} \times_{\widecheck{\Omega}_k} \widetilde{\U}^{(t)\top}_k  \times_{i \in \overline{\Omega}_k} \widetilde{\U}^{(t)\top}_{i'}  \right) \right)_{\max(r_k)}\right\|_q.$$ 

To simplify the notation, for $S_j \in \cS_j^{(-\bar{k})}$, we let $S_{j1} = S_j \bigcap \underline{\Omega}_k, S_{j2} = S_j \bigcap ( \widecheck{\Omega}_k \bigcup \overline{\Omega}_k )$ and $S_{j1}^c = S_j^c \bigcap \underline{\Omega}_k, S_{j2}^c = S_j^c \bigcap ( \widecheck{\Omega}_k \bigcup \overline{\Omega}_k )$. First we can introduce $\I = \left(\Proj_{\U_{i'}} + \Proj_{\U_{i'\perp}}\right) $ in the expression,
\begin{equation}\label{eq: equiva-expre}
	\begin{split}
		& \left\|\left(\mathcal{M}_{\bar{k}}\left(\bcZ \times_{i \in \underline{\Omega}_k} \widetilde{\U}^{(t+1)\top}_{i'} \times_{\widecheck{\Omega}_k} \widetilde{\U}^{(t)\top}_k  \times_{i \in \overline{\Omega}_k} \widetilde{\U}^{(t)\top}_{i'}  \right) \right)_{\max(r_k)}\right\|_q \\
= & \Big\|\Big(\mathcal{M}_{\bar{k}}\Big(\bcZ \times_{i \in \underline{\Omega}_k} \widetilde{\U}^{(t+1)\top}_{i'} \left(\Proj_{\U_{i'}} + \Proj_{\U_{i'\perp}}\right) \times_{\widecheck{\Omega}_k} \widetilde{\U}^{(t)\top}_k \left(\Proj_{\U_k} + \Proj_{\U_{k\perp}}\right) \\
& \quad \quad \quad \quad  \times_{i \in \overline{\Omega}_k} \widetilde{\U}^{(t)\top}_{i'}  \left(\Proj_{\U_{i'}} + \Proj_{\U_{i'\perp}}\right) \Big) \Big)_{\max(r_k)}\Big\|_q.
	\end{split}
\end{equation}
 Then 
\begin{equation}\label{ineq:th2-long-argument}
\begin{split}
\eqref{eq: equiva-expre} \overset{(a)}\leq & \left\|\left(\mathcal{M}_{\bar{k}}\left(\bcZ \times_{i \in \underline{\Omega}_k} \widetilde{\U}^{(t+1)\top}_{i'} \Proj_{\U_{i'}} \times_{\widecheck{\Omega}_k} \widetilde{\U}^{(t)\top}_k \Proj_{\U_k}  \times_{i \in \overline{\Omega}_k} \widetilde{\U}^{(t)\top}_{i'}  \Proj_{\U_{i'}}  \right) \right)_{\max(r_k)}\right\|_q  \\
& + \sum_{S_1 \in \cS_1^{(-\bar{k})}} \Big\|\Big(\mathcal{M}_{\bar{k}}\Big(\bcZ \times_{i \in S_{11} } \widetilde{\U}^{(t+1)\top}_{i'} \Proj_{\U_{i'\perp}} \times_{i \in S_{11}^c } \widetilde{\U}^{(t+1)\top}_{i'} \Proj_{\U_{i'}} \\ 
&\quad \quad \quad \quad \times_{i \in S_{12}} \widetilde{\U}^{(t)\top}_{i'} \Proj_{\U_{i'\perp}} \times_{i \in S^c_{12} } \widetilde{\U}^{(t)\top}_{i'} \Proj_{\U_{i'}}  \Big) \Big)_{\max(r_k)}\Big\|_q  +   \cdots  \\
&  + \sum_{S_j \in \cS_j^{(-\bar{k})}} \Big\|\Big(\mathcal{M}_{\bar{k}}\Big(\bcZ \times_{i \in S_{j1} } \widetilde{\U}^{(t+1)\top}_{i'} \Proj_{\U_{i'\perp}} \times_{i \in S_{j1}^c } \widetilde{\U}^{(t+1)\top}_{i'} \Proj_{\U_{i'}} \\
& \quad \quad \quad \quad \times_{i \in S_{j2} } \widetilde{\U}^{(t)\top}_{i'} \Proj_{\U_{i'\perp}} \times_{i \in S^c_{j2}} \widetilde{\U}^{(t)\top}_{i'} \Proj_{\U_{i'}}  \Big) \Big)_{\max(r_k)}\Big\|_q + \cdots\\
& + \left\|\left(\mathcal{M}_{\bar{k}}\left(\bcZ \times_{i \in \underline{\Omega}_k} \widetilde{\U}^{(t+1)\top}_{i'} \Proj_{\U_{i'\perp}} \times_{\widecheck{\Omega}_k} \widetilde{\U}^{(t)\top}_k \Proj_{\U_{k\perp}}  \times_{i \in \overline{\Omega}_k} \widetilde{\U}^{(t)\top}_{i'}  \Proj_{\U_{i'\perp}} \right) \right)_{\max(r_k)}\right\|_q\\
=& \left\|\left(\mathcal{M}_{\bar{k}}\left(\bcZ \times_{i \in \underline{\Omega}_k} \widetilde{\U}^{(t+1)\top}_{i'} \Proj_{\U_{i'}} \times_{\widecheck{\Omega}_k} \widetilde{\U}^{(t)\top}_k \Proj_{\U_{k}}  \times_{i \in \overline{\Omega}_k} \widetilde{\U}^{(t)\top}_{i'}  \Proj_{\U_{i'}}  \right) \right)_{\max(r_k)}\right\|_q  \\
& +\sum_{j=1}^{d-1} \sum_{S_j \in \cS_j^{(-\bar{k})}}  \Big\|\Big(\mathcal{M}_{\bar{k}}\Big(\bcZ \times_{i \in S_{j1} } \widetilde{\U}^{(t+1)\top}_{i'} \Proj_{\U_{i'\perp}} \times_{i \in S_{j1}^c } \widetilde{\U}^{(t+1)\top}_{i'} \Proj_{\U_{i'}} \\ 
& \quad \quad \quad \quad  \times_{i \in S_{j2}} \widetilde{\U}^{(t)\top}_{i'} \Proj_{\U_{i'\perp}} \times_{i \in S^c_{j2}} \widetilde{\U}^{(t)\top}_{i'} \Proj_{\U_{i'}}  \Big) \Big)_{\max(r_k)}\Big\|_q.
\end{split}
\end{equation}
Here $(a)$ is due to the triangle inequality for truncated Schatten-$q$ norm given in Lemma 4 of \cite{luo2020schatten}. The right hand side of \eqref{ineq:th2-long-argument} can be divided into the sum of $d$ groups and the value of $j^{th}$ group is denoted as $G_j$ where $$G_0:=\left\|\left(\mathcal{M}_{\bar{k}}\left(\bcZ \times_{i \in \underline{\Omega}_k} \widetilde{\U}^{(t+1)\top}_{i'} \Proj_{\U_{i'}} \times_{\widecheck{\Omega}_k} \widetilde{\U}^{(t)\top}_k \Proj_{\U_{k}}  \times_{i \in \overline{\Omega}_k} \widetilde{\U}^{(t)\top}_{i'}  \Proj_{\U_{i'}}  \right) \right)_{\max(r_k)}\right\|_q$$ and for $1 \leq j \leq d-1$, define $G_j$ to be,
\begin{equation*}
\begin{split}
G_j = \sum_{S_j \in \cS_j^{(-\bar{k})}} \Big\|\Big(\mathcal{M}_{\bar{k}}\Big(\bcZ \times_{i \in S_{j1} } \widetilde{\U}^{(t+1)\top}_{i'} \Proj_{\U_{i'\perp}} \times_{i \in S_{j1}^c } \widetilde{\U}^{(t+1)\top}_{i'} \Proj_{\U_{i'}} \\
\times_{i \in S_{j2}} \widetilde{\U}^{(t)\top}_{i'} \Proj_{\U_{i'\perp}} \times_{i \in S^c_{j2}} \widetilde{\U}^{(t)\top}_{i'} \Proj_{\U_{i'}}  \Big) \Big)_{\max(r_k)}\Big\|_q.
\end{split}
\end{equation*}
Next we are going to upper bound $G_j \,( 0 \leq j \leq d-1)$.

\begin{itemize}[leftmargin=*]
	\item Upper Bound of $G_0$.
	\begin{equation*}
	\begin{split}
	G_0 = &  \left\|\left(\mathcal{M}_{\bar{k}}\left(\bcZ \times_{i \in \underline{\Omega}_k} \widetilde{\U}^{(t+1)\top}_{i'} \Proj_{\U_{i'}} \times_{\widecheck{\Omega}_k} \widetilde{\U}^{(t)\top}_k \Proj_{\U_{k}}  \times_{i \in \overline{\Omega}_k} \widetilde{\U}^{(t)\top}_{i'}  \Proj_{\U_{i'}}  \right) \right)_{\max(r_k)}\right\|_q\\
	=  &\Big\|\Big( \mathcal{M}_{\bar{k}}\Big( \left(\bcZ \times_{i \neq \bar{k}}\U_{i'}^\top \right)\times_{i \in \underline{\Omega}_k} \widetilde{\U}^{(t+1)\top}_{i'} \U_{i'} \times_{\widecheck{\Omega}_k} \widetilde{\U}^{(t)\top}_k \U_{k} \\ 
	& \quad \quad \quad \quad \quad \times_{i \in \overline{\Omega}_k} \widetilde{\U}^{(t)\top}_{i'} \U_{i'}  \Big) \Big)_{\max(r_k)}\Big\|_q\\
	\leq  & \left\|\left(  \mathcal{M}_{\bar{k}}\left(\bcZ \times_{i \neq \bar{k}}\U_{i'}^\top  \right) \right)_{\max(r_k)}\right\|_q = \tau_{1k}.
	\end{split}
	\end{equation*}
	Here the inequality is due to the fact $\|\widetilde{\U}^{(t)\top}_{i'} \U_{i'} \| \leq 1$ for any $t$ and the following fact: for any matrix $\A$ s.t. $\|\A\| \leq 1$, \begin{equation} \label{ineq: matricization shrink}
	\sigma_k\left(\mathcal{M}_i(\bcZ \times_j \A)\right) = \sigma_k \left( \mathcal{M}_i(\bcZ) \cdot \A^\top  \right) \leq \sigma_k \left( \mathcal{M}_i(\bcZ) \right) \| \A\| \leq \sigma_k \left( \mathcal{M}_i(\bcZ) \right),
	\end{equation}
	where the first inequality is due to Lemma 7 of \cite{luo2020schatten}.
	\item Upper Bound of $G_j\, (1 \leq j \leq d-1)$.
	\begin{equation}\label{ineq: upper bound of G_k prepare}
	\begin{split}
	G_j & \leq {d-1 \choose j} \times \max_{S_j \in \cS_j^{(-\bar{k})}}\\
	& \Big\|\Big(\mathcal{M}_{\bar{k}}\Big(\bcZ \times_{i \in S_{j1} } \widetilde{\U}^{(t+1)\top}_{i'} \Proj_{\U_{i'\perp}} \times_{i \in S_{j1}^c } \widetilde{\U}^{(t+1)\top}_{i'} \Proj_{\U_{i'}} \times_{i \in S_{j2}} \widetilde{\U}^{(t)\top}_{i'} \Proj_{\U_{i'\perp}} \\
	& \quad \quad \quad \times_{i \in S^c_{j2}} \widetilde{\U}^{(t)\top}_{i'} \Proj_{\U_{i'}}  \Big) \Big)_{\max(r_k)}\Big\|_q\\
	& \overset{(a)} \leq {d-1 \choose j} \max_{S_j \in \cS_{j}^{(-\bar{k})}} \left(\prod_{i \in S_{j1}} \left\| \widetilde{\U}^{(t+1)\top}_{i'} \U_{i'\perp} \right\|_q \prod_{i \in S_{j2}} \left\| \widetilde{\U}^{(t)\top}_{i'} \U_{i'\perp} \right\|_q \right) \times \\
	& \left\|\left(\mathcal{M}_{\bar{k}} \left(\bcZ \times_{i \in S_{j1} } \frac{\widetilde{\U}^{(t+1)\top}_{i'} \Proj_{\U_{i'\perp}}}{\left\|\widetilde{\U}^{(t+1)\top}_{i'} \Proj_{\U_{i'\perp}}\right\|_q} \times_{i \in S_{j2}} \frac{\widetilde{\U}^{(t)\top}_{i'} \Proj_{\U_{i'\perp}}}{\left\|\widetilde{\U}^{(t)\top}_{i'} \Proj_{\U_{i'\perp}}\right\|_q} \times_{i \in S_j^c} \U_{i'}^\top   \right) \right)_{\max(r_k)}\right\|_q,
	\end{split}
	\end{equation} where (a) is due to the fact that $\| \tilde{\U}_i^{(t)\top} \U_{i'} \| \leq 1$ for any $t$. By simplifying \eqref{ineq: upper bound of G_k prepare}, we get
	\begin{equation} \label{ineq: upper bound of G_k}
	\begin{split}
	G_j & \leq {d-1 \choose j} \max_{S_j \in \cS_{j}^{(-\bar{k})}} \left(\prod_{i \in S_{j1}} \left\| \widetilde{\U}^{(t+1)\top}_{i'} \U_{i'\perp} \right\|_q \prod_{i \in S_{j2}} \left\| \widetilde{\U}^{(t)\top}_{i'} \U_{i'\perp} \right\|_q \right) \\
	& \times \Big\|\Big(\mathcal{M}_{\bar{k}} \Big(\bcZ \times_{i \in S_{j1} } \frac{\widetilde{\U}^{(t+1)\top}_{i'} \U_{i'\perp}}{\left\|\widetilde{\U}^{(t+1)\top}_{i'} \U_{i'\perp}\right\|_q} \U_{i'\perp}^\top \\
	& \quad \quad \quad \quad \quad  \times_{i \in S_{j2}} \frac{\widetilde{\U}^{(t)\top}_{i'} \U_{i'\perp}}{\left\|\widetilde{\U}^{(t)\top}_{i'} \U_{i'\perp}\right\|_q} \U_{i'\perp}^\top \times_{i \in S_j^c} \U_{i'}^\top   \Big) \Big)_{\max(r_k)}\Big\|_q\\
	& \overset{(b)} \leq {d-1 \choose j} \max_{S_j \in \cS_{j}^{(-\bar{k})}} \left(\prod_{i \in S_{j1}} \left\| \widetilde{\U}^{(t+1)\top}_{i'} \U_{i'\perp} \right\|_q \prod_{i \in S_{j2}} \left\| \widetilde{\U}^{(t)\top}_{i'} \U_{i'\perp} \right\|_q \right) \times \tau_{j+1} \\
	&  \leq  {d-1 \choose j} \max_{S_j \in \cS_{j}^{(-\bar{k})}} \left( e_t \right)^{|S_{j2}|} \left( e_{t+1} \right)^{|S_{j1}|} \tau_{j+1}
	\end{split}
	\end{equation}
	where (b) is due to the fact that $\left\| \frac{\widetilde{\U}^{(t+1)\top}_{i'} \U_{i'\perp}}{\left\|\widetilde{\U}^{(t+1)\top}_{i'} \U_{i'\perp}\right\|_q}\right\|_q \leq 1$ and the definition of $\tau_{j+1}$.
\end{itemize}
So in summary, plug \eqref{ineq: upper bound of G_k} into \eqref{ineq:th2-long-argument}, now we have the following bound
\begin{equation}\label{ineq: upper bnd of SVD error Z}
\begin{split}
&\left\|\left(\mathcal{M}_{\bar{k}}\left(\bcZ \times_{i \in \underline{\Omega}_k} \widetilde{\U}^{(t+1)\top}_{i'} \times_{\widecheck{\Omega}_k} \widetilde{\U}^{(t)\top}_k  \times_{i \in \overline{\Omega}_k} \widetilde{\U}^{(t)\top}_{i'}  \right) \right)_{\max(r_k)}\right\|_q\\
\leq & \tau_{1k} + \sum_{j=1}^{d-1} {d-1 \choose j} \max_{S_j \in \cS_{j}^{(-\bar{k})}} \left( e_t \right)^{|S_{j2}|} \left( e_{t+1} \right)^{|S_{j1}|} \tau_{j+1}.
\end{split}
\end{equation}
{\noindent \bf Step 2.} In this step we want to use induction to prove the following claims,
\begin{equation}\label{ineq:claim-e_t-induction}
e_{t} \leq 2^{\frac{d+3}{2}}\frac{\tau_{1}}{\lambda} + e_0/2^t \quad \text{ and } \quad  \tilde{e}_{t} \leq 2^{\frac{d+3}{2}}\frac{\tau_{1}}{\lambda} + \tilde{e}_0/2^t;\quad   t=0,1,\ldots. 
\end{equation}

Since the proof of two statements in \eqref{ineq:claim-e_t-induction} are similar, we mainly focus on the proof of the first statement. Claim \eqref{ineq:claim-e_t-induction} clearly holds if $t=0$. Assume the first claim of \eqref{ineq:claim-e_t-induction} holds for $t$ and next we show it also holds for $t+1$. 

Let's first show the upper bound of $e_{t}$ can be used to upper bound $e_{t+1,1}$. Notice \begin{equation}\label{ineq: recursive e_t1 bound}
\begin{split}
2\left\| \left( \mathcal{M}_{\bar{1}}\left(\bcZ \times_{i \neq \bar{1}} \widetilde{\U}_{i'}^{(t)\top} \right) \right)_{\max(r_1)}  \right\|_q \overset{(a)}\geq & \left\| \widetilde{\U}_{1\perp}^{(t+1)\top} \U_{1} \right\|_q  \sigma_{r_1} \left( \mathcal{M}_{\bar{1}}\left(\bcT \times_{i \neq \bar{1}}\widetilde{\U}_{i'}^{(t)\top} \right)  \right)  \\
\overset{(b)}\geq & \left(e_{t+1,1}\right) \sigma_{r_1} \left( \T_1 \left(  \otimes_{i\neq \bar{1}} \U_{i'} \right) \cdot \left( \otimes_{i\neq \bar{1}} \U_{i'}^\top \widetilde{\U}_{i'}^{(t)} \right)  \right) \\
\overset{(c)} \geq & \left(e_{t+1,1}\right)  \lambda \left( \prod_{i\neq \bar{1}} \sigma_{\min} \left(  \U_{i'}^\top \widetilde{\U}_{i'}^{(t)} \right)  \right) \\
\overset{(d)}\geq  & \left(e_{t+1,1}\right)  \lambda (1 - \tilde{e}_t^2)^{\frac{d-1}{2}}.
\end{split}
\end{equation}
Here (a) is due to Theorem 5 in \cite{luo2020schatten} and the fact the left singular space of $\mathcal{M}_{\bar{1}}\left(\bcT \times_{i \neq \bar{1}} \widetilde{\U}_{i'}^{(t)\top} \right)$ is $\U_1$ and its rank is less than $r_1$, (b) is due to the fact that the right singular space of $\T_1$ is $\otimes_{i\neq \bar{1}}\U_{i'}$, equation \eqref{eq: matricization relationship}, (c) is by Lemma 7 of \cite{luo2020schatten} and (d) is due to the fact $\sigma_{\min} \left(\U_{i'}^\top \widetilde{\U}_i^{(t)}  \right) = \sqrt{1 - \|\U_{i'\perp}^\top \widetilde{\U}_i^{(t)}  \|^2} \geq \sqrt{ 1 - \tilde{e}_{t}^2  }$ by Lemma 1 of \cite{cai2018rate}.


Plugging the upper bound of $\left\| \left( \mathcal{M}_{\bar{1}}\left(\bcZ \times_{i \neq \bar{1}} \widetilde{\U}_{i'}^{(t)\top} \right) \right)_{\max(r_1)}  \right\|_q$ in \eqref{ineq: upper bnd of SVD error Z} into \eqref{ineq: recursive e_t1 bound}, we get
\begin{equation} \label{ineq: recursive e_t1 bound eq2}
e_{t+1,1} \leq 2 \frac{\tau_{11} + \sum_{j=1}^{d-1} {d-1 \choose j} \left( e_t \right)^{j} \tau_{j+1}  }{ \lambda (1 - \tilde{e}_t^2)^{\frac{d-1}{2}} }.
\end{equation}

Next we show under condition \eqref{ineq: signal strength condition general d}, $e_{t+1,1} \leq 2^{\frac{d+3}{2}}\frac{\tau_{1}}{\lambda} + e_0/2^t$, i.e., the upper bound of $e_t$ can be used to upper bound $e_{t+1,1}$. First, under \eqref{ineq: signal strength condition general d}, we have
\begin{equation} \label{cond: lambda condition}
\lambda \geq  2^{\frac{d+6}{2}} \tau_1  \lor \left(\sum_{j=1}^{d-1} 2^{\frac{d+4-j}{2}} {d-1 \choose j} \tau_{j+1}\right),
\end{equation} due to the fact that $\xi \geq \tau_j$ for $j = 1, \ldots, d$. And under the assumption of $\lambda$, we have $\tilde{e}_t \leq \sqrt{2}/2$. 
So by \eqref{ineq: recursive e_t1 bound eq2}, we have
\begin{equation}\label{ineq: recursive varify}
\begin{split}
\left(\frac{1}{2}\right)^{\frac{d-1}{2}} e_{t+1,1} & \leq 2\frac{\tau_{11} + \sum_{j=1}^{d-1} {d-1 \choose j} \left( e_t \right)^{j} \tau_{j+1}  }{ \lambda } \\
& \leq  2 \left( \frac{\tau_{11}}{\lambda} + \frac{\sum_{j=1}^{d-1} {d-1 \choose j} (\frac{\sqrt{2}}{2})^{j-1} \tau_{j+1}  }{ \lambda } e_t \right)\\
& \overset{ \eqref{cond: lambda condition} }\leq 2 \left( \frac{\tau_{11}}{\lambda} + 2^{-\frac{d+3}{2}} e_t \right)
\end{split}
\end{equation}

Since \eqref{ineq:claim-e_t-induction} holds for $t$ and plug in the upper bound of $e_t$ in \eqref{ineq: recursive varify}, multiply $2^{\frac{d-1}{2}}$ at both side of \eqref{ineq: recursive varify}, we get
\begin{equation*}
e_{t+1,1} \leq 2^{\frac{d+3}{2}} \frac{\tau_1}{\lambda} + \frac{e_0}{2^{t+1}} \leq 2^{\frac{d+3}{2}} \frac{\tau_1}{\lambda} + \frac{e_0}{2^{t}}.
\end{equation*}
So the upper bound of $e_{t}$ can be used to bound $e_{t+1,1}$ and by doing similar analysis for all modes, we conclude that the upper bound for $e_t$ also holds for $e_{t+1}$, i.e., we have 
\begin{equation}\label{ineq: e_t+1 naive bound}
e_{t+1} \leq 2^{\frac{d+3}{2}} \frac{\tau_1}{\lambda} + \frac{e_0}{2^{t}}. 
\end{equation}


Now we can show the first statement in \eqref{ineq:claim-e_t-induction} also holds for $t+1$ given it holds for $t$. With \eqref{ineq: e_t+1 naive bound}, when we do the similar analysis of \eqref{ineq: recursive e_t1 bound} for other modes, we can use the same upper bound of $e_t$ to bound $e_{t+1}$. So repeat \eqref{ineq: recursive e_t1 bound}, for any $1 \leq k \leq m$, we have 
\begin{equation} \label{ineq: recursive e_tk bound}
e_{t+1, k} \left(1-\tilde{e}_t^2\right)^{\frac{d-1}{2}} \leq 2 \frac{ \tau_{1k} + \sum_{j=1}^{d-1} {d-1 \choose j} \left(e_t \right)^j \tau_{j+1} }{\lambda}.
\end{equation}
Since it is true for every $k$, we obtain the following recursive inequality for $e_{t+1}$,
\begin{equation}\label{ineq: recursive e_t bound}
e_{t+1} \left( \frac{1}{2} \right)^{\frac{d-1}{2}} \leq 2 \frac{ \tau_1 + \sum_{j=1}^{d-1} {d-1 \choose j} \left(e_t \right)^j \tau_{j+1} }{\lambda}.
\end{equation}
By applying the same argument as in \eqref{ineq: recursive varify}, we can show that  
\begin{equation*}
e_{t+1} \leq 2^{\frac{d+3}{2}} \frac{\tau_1}{\lambda} + \frac{e_0}{2^{t+1}},
\end{equation*} given \eqref{ineq:claim-e_t-induction} holds for $e_t$. Similarly we can show $\tilde{e}_{t+1} \leq 2^{\frac{d+3}{2}} \frac{\tau_1}{\lambda} + \frac{\tilde{e}_0}{2^{t+1}}$

Thus, \eqref{ineq:claim-e_t-induction} holds for all $t$. Given $t_{\max} \geq \log(e_0 \lambda/\tau_{1}) \vee 1$, we get the upper bounds for $e_t$ and $\tilde{e}_t$:
\begin{equation} \label{ineq: upper bound of e_t}
e_{t} \leq \left(2^{\frac{d+3}{2}} + 1\right) \frac{\tau_1}{\lambda}, \quad \tilde{e}_{t} \leq \left(2^{\frac{d+3}{2}} + 1\right) \frac{\tau_1}{\lambda}
\end{equation}

Plug \eqref{ineq: recursive e_t bound} and \eqref{cond: lambda condition} into the upper bound \eqref{ineq: recursive e_tk bound}, we get
\begin{equation*} 
\begin{split}
e_{t+1,k} & \leq 2\left( 1- \frac{\left(2^{\frac{d+3}{2}} + 1\right)^2 \tau_1^2}{\lambda^2} \right)^{- \frac{d-1}{2} } \left(\frac{\tau_{1k}}{\lambda} +  \sum_{j=1}^{d-1}\frac{ {d-1 \choose j}  \left(2^{\frac{d+3}{2}}+1\right)^j \tau_1^j \tau_{j+1} }{\lambda^{j+1}} \right),
\end{split}
\end{equation*} for $k =1, \ldots, m$. 

Finally the perturbation of signal subspaces follows by observing that $e_{t,k} := \left\|(\widetilde{\U}_{k\perp}^{(t)})^\top \U_k \right\|_q = \left\|\sin\Theta\left(\widetilde{\U}^{(t)}_k, \U_k\right)\right\|_q$ due to Lemma 6 of \cite{luo2020schatten}.

{\noindent \bf Step 3.} In this step, we are going to give an upper bound for the tensor reconstruction error for $\| \widetilde{\bcT}\times_{\Omega_1} P_{\widehat{\U}_1}\times \cdots \times_{\Omega_m} P_{\widehat{\U}_m} - \bcT\|_{\tHS}$. 

 First, notice the following decomposition
\begin{equation*}
\begin{split}
\bcT = & \bcT \times_1 \left(\Proj_{\widehat{\U}_{1'}}+\Proj_{\widehat{\U}_{1'\perp}}\right) \times_2 \left(P_{\widehat{\U}_{2'}}+P_{\widehat{\U}_{2'\perp}}\right) \times \cdots \times_d \left(P_{\widehat{\U}_{d'}}+P_{\widehat{\U}_{d'\perp}}\right)\\
= & \bcT \times_1 \Proj_{\widehat{\U}_{1'}} \times \cdots \times_d P_{\widehat{\U}_{d'}} + \bcT \times_1 P_{\widehat{\U}_{1'\perp}} \times_2 P_{\widehat{\U}_{2'}} \times \cdots \times_d P_{\widehat{\U}_{d'}}\\
& + \bcT \times_1 \I_{p_{1'}} \times_2 P_{\widehat{\U}_{2'\perp}} \times \cdots \times_d P_{\widehat{\U}_{d'}} + \bcT \times_1 \I_{p_{1'}} \times_2 \I_{p_{2'}} \times \cdots \times_d P_{\widehat{\U}_{d'\perp}}\\
& + \ldots + \bcT \times_{i \leq d-1} \I_{p_{i'}} \times_d P_{\widehat{\U}_{d'\perp}}\\
= & \bcT \times_{\Omega_1} \Proj_{\widehat{\U}_{1}} \times \cdots \times_{\Omega_m} P_{\widehat{\U}_m} + \sum_{k=1}^d \bcT \times_{i < k} \I_{p_{i'}} \times_{k} P_{\widehat{\U}_{k'\perp}} \times_{i> k} P_{\widehat{\U}_{i'}}.
\end{split}
\end{equation*}
Thus,
\begin{equation}\label{ineq: hatT l2 bound}
\begin{split}
& \left\|\widetilde{\bcT} \times_{\Omega_1} P_{\widehat{\U}_1} \times \cdots \times_{\Omega_m} P_{\widehat{\U}_m} - \bcT \right\|_{\tHS} \\
\leq & \left\|(\widetilde{\bcT}-\bcT) \times_{\Omega_1} P_{\widehat{\U}_1} \times \cdots \times_{\Omega_m} P_{\widehat{\U}_m}\right\|_{\tHS} + \sum_{k=1}^d \left\|\bcT \times_{i < k} \I_{p_{i'}} \times_{k} P_{\widehat{\U}_{k'\perp}} \times_{i> k} P_{\widehat{\U}_{i'}}\right\|_{\tHS}\\
\leq & \left\| \bcZ \times_{\Omega_1} P_{\widehat{\U}_1} \times \cdots \times_{\Omega_m} P_{\widehat{\U}_m}\right\|_{\tHS} + \sum_{k=1}^d \left\|\widehat{\U}_{k'\perp}^\top \mathcal{M}_{k}(\bcT)\right\|_F.
\end{split}
\end{equation}
Notice that for $k_1, k_2 \in \Omega_k$, $\mathcal{M}_{k_1}(\bcT) = \mathcal{M}_{k_2}(\bcT) P_1$ where $P_1$ is a permutation matrix. So $\sum_{k=1}^d \left\|\widehat{\U}_{k'\perp}^\top \mathcal{M}_{k}(\bcT)\right\|_F = \sum_{k=1}^m |\Omega_k| \left\|\widehat{\U}_{k\perp}^\top \T_k \right\|_F $.

\begin{equation}\label{ineq: hatU_kT_k bound}
\begin{split}
&\left\|\widehat{\U}_{k\perp}^\top \T_k\right\|_F \left( 1- \frac{\left(2^{\frac{d+3}{2}} + 2\right)^2 \tau_1^2}{\lambda^2} \right)^{ \frac{d-1}{2} } \\ 
\overset{(a)}\leq & \left\|\widehat{\U}_{k\perp}^\top \T_k\right\|_F \sigma_{\min} \left( \left( \otimes_{i\in \underline{\Omega}_k} \U_{i'}^\top\widetilde{\U}_i^{(t_{\max})}  \otimes_{i \in \widecheck{\Omega}_k \bigcup \overline{\Omega}_k} \U_{i'}^\top\widetilde{\U}_{i'}^{(t_{\max} - 1)}  \right) \right)\\
\overset{(b)}= &\left\|\widehat{\U}_{k\perp}^\top \T_k \otimes_{i\neq \bar{k}} \U_{i'} \right\|_F \sigma_{\min} \left( \left( \otimes_{i\in \underline{\Omega}_k} \U_{i'}^\top\widetilde{\U}_i^{(t_{\max})}  \otimes_{i \in \widecheck{\Omega}_k \bigcup \overline{\Omega}_k} \U_{i'}^\top\widetilde{\U}_{i'}^{(t_{\max} - 1)}  \right) \right)\\
\overset{(c)}\leq & \left\|\widehat{\U}_{k\perp}^\top \T_k \left(\otimes_{i\neq \bar{k}} \U_{i'} \right) \left( \otimes_{i\neq \bar{k}} \U_{i'} \right)^\top \left(\otimes_{i\in \underline{\Omega}_k } \widetilde{\U}_{i'}^{(t_{\max})} \otimes_{i\in \widecheck{\Omega}_k \bigcup \overline{\Omega}_k  } \widetilde{\U}_{i'}^{(t_{\max} - 1)} \right)   \right\|_F\\
\overset{(d)}\leq & 2\left\|  \left( \mathcal{M}_{\bar{k}}\left(  \bcZ \times_{i \in \underline{\Omega}_k} \widetilde{\U}_{i'}^{(t_{\max})\top} \times_{i \in \widecheck{\Omega}_k \bigcup \overline{\Omega}_k } \widetilde{\U}_{i'}^{(t_{\max}-1)\top} \right)  \right)_{\max(r_k)} \right\|_F \overset{(e)} \leq 2\xi,
\end{split}
\end{equation}
where $\xi := \underset{\|\bcY \|_{\tHS} \leq 1,  \rank(\bcY) \leq (r_{1'}, \ldots, r_{d'})}\sup \langle \bcZ, \bcY \rangle$. (a) is due to the fact $\sigma_{\min} \left(\U_{i'}^\top \widetilde{\U}_i^{(t)}  \right) = \sqrt{1 - \|\U_{i'\perp}^\top \widetilde{\U}_i^{(t)}  \|^2} \geq \sqrt{ 1 - \tilde{e}_{t}^2  }$ and $\max(\tilde{e}_{t_{\max}},\tilde{e}_{t_{\max}-1} ) \leq \frac{\left(2^{\frac{d+3}{2}} + 2\right) \tau_1}{\lambda}$; (b) is due to the fact that the right singular space of $\T_k$ is $\otimes_{i \neq \bar{k}} \U_{i'}$; (c) is due to Lemma 7 of \cite{luo2020schatten}, equation \eqref{eq: matricization relationship} and properties of Kronecker product, (d) is due to Theorem 2 in \cite{luo2020schatten} and $$\rank\left(\widehat{\U}_{k\perp}^\top \T_k \left(\otimes_{i\in \underline{\Omega}_k } \widetilde{\U}_{i'}^{(t_{\max})} \otimes_{i\in \widecheck{\Omega}_k \bigcup \overline{\Omega}_k  } \widetilde{\U}_{i'}^{(t_{\max} - 1)} \right)\right) \leq r_k;$$ and the last inequality (e) is due to the fact that
\begin{equation*}
\begin{split}
&\left\|  \left( \mathcal{M}_{\bar{k}}\left(  \bcZ \times_{i \in \underline{\Omega}_k} \widetilde{\U}_{i'}^{(t_{\max})\top} \times_{i \in \widecheck{\Omega}_k \bigcup \overline{\Omega}_k } \widetilde{\U}_{i'}^{(t_{\max}-1)\top} \right)  \right)_{\max(r_k)} \right\|_F\\
= &\sup_{\|\X\|_F \leq 1, \rank(\X) \leq r_k} \left\langle \mathcal{M}_{\bar{k}}\left(  \bcZ \times_{i \in \underline{\Omega}_k} \widetilde{\U}_{i'}^{(t_{\max})\top} \times_{i \in \widecheck{\Omega}_k \bigcup \overline{\Omega}_k } \widetilde{\U}_{i'}^{(t_{\max}-1)\top} \right), \X \right\rangle \leq \xi,
\end{split}
\end{equation*}
for $k =1, \ldots, d$ and the equality is due to Lemma 2 of \cite{luo2020schatten}.

Combined with the fact $\left\| \bcZ \times_{\Omega_1} P_{\widehat{\U}_1} \times \cdots \times_{\Omega_m} P_{\widehat{\U}_m}\right\|_{\tHS} \leq \xi$ and \eqref{ineq: hatU_kT_k bound} and plug them into \eqref{ineq: hatT l2 bound}, we have
\begin{equation}
\begin{split}
&\left\|\widetilde{\bcT} \times_{\Omega_1} P_{\widehat{\U}_1} \times \cdots \times_{\Omega_m} P_{\widehat{\U}_m} - \bcT \right\|_{\tHS} \\
\leq & \left\| \bcZ \times_{\Omega_1} P_{\widehat{\U}_1} \times \cdots \times_{\Omega_m} P_{\widehat{\U}_m}\right\|_{\tHS}+ \sum_{k=1}^d \left\|\widehat{\U}_{k'\perp}^\top \mathcal{M}_k(\bcT)  \right\|_F\\
= & \left\|  \bcZ\times_{\Omega_1} \widehat{\U}_1^\top \times \cdots \times_{\Omega_m} \widehat{\U}_m^\top \right\|_{\tHS} + \sum_{k=1}^m |\Omega_k| \left\|\widehat{\U}_{k\perp}^\top \mathcal{M}_k(\bcT) \right\|_F\\
\leq &  \left(  1+ 2d \left( 1- \frac{\left(2^{\frac{d+3}{2}} + 2\right)^2 \tau_1^2}{\lambda^2} \right)^{- \frac{d-1}{2} }  \right)\xi. 
\end{split}
\end{equation}
Finally, notice that \eqref{ineq: t step tensor reconstruction order d} follows exact the same proof as Step 3 except that in \eqref{ineq: hatU_kT_k bound}, we need to replace $\widetilde{\U}_i^{(t_{\max})},\widetilde{\U}_i^{(t_{\max}-1)}$ with $\widetilde{\U}_i^{(t)},\widetilde{\U}_i^{(t-1)}$ and use $\max(\tilde{e}_t, \tilde{e}_{t-1}) \leq 2^{\frac{d+3}{2}} \frac{\tau_1}{\lambda} + \frac{\tilde{e}_0}{2^{t-1}}$. Therefore, we have finished the proof of this theorem. \quad $\blacksquare$

\subsection{Proof of Theorem \ref{th: tensor denoising}}
In this setting, we consider applying Corollary \ref{th:HOOI d diff}. To apply Corollary \ref{th:HOOI d diff}, we only need to compute $\tau_{1i}$, $\xi$. 

Consider the case $q = \infty$ and let's first bound $\tau_{1k}$. Notice 
\begin{equation*}
\mathcal{M}_k(\bcZ \times_{i \neq k} \U_i^\top ) = \mathcal{M}_k(\bcZ) \otimes_{i \neq k} \U_i.
\end{equation*} Each row of $\mathcal{M}_k(\bcZ) \otimes_{i \neq k} \U_i$ is independent multivariate Gaussian with covariance matrix $\sigma \I_{r_{-k}}$ where $r_{-k} = \prod_{i \neq k} r_i$. By random matrix theory \citep{vershynin2010introduction}, we have
\begin{equation*}
\bbP \left( \| \mathcal{M}_k(\bcZ \times_{i \neq k} \U_i^\top )  \|/\sigma \leq \sqrt{p_k} + \sqrt{r_{-k}} \right) \geq 1 - c\exp(-(p_k + r_{-k})/2).
\end{equation*} 
Since $r_{\max} \leq p_{\min}^{\frac{1}{d-1}}$, we have $\tau_{1k} \leq C \sigma \sqrt{p_k}$ with probability at least $1 - c \exp(-(p_k + r_{-k})/2)$.

Similarly by Lemma 5 of \cite{zhang2018tensor}, we have 
\begin{equation*}
\xi \leq C \sigma \sqrt{\sum_{i=1}^d p_i r_i}
\end{equation*} w.p. at least $1 - \exp(-c p_{\min})$.

As we mentioned in Remark \ref{rem: size of xi constant}, $\tau_j \leq \xi$. So the results follows by plugging the bound of $\tau_j, \xi$ into Corollary \ref{th:HOOI d diff} and noticing that 
$$\frac{\sum_{j=1}^{d-1} {d-1 \choose j} \left( \left(2^{\frac{d+3}{2}}+1\right) \tau_1/\lambda \right)^{j} \tau_{j+1} }{\lambda} \leq C \frac{\tau_1 \max_j \tau_j}{\lambda^2},$$
under the assumption of the signal to ratio $\lambda/\sigma$. 
\quad $\blacksquare$

\subsection{Proof of Lemma \ref{lem: why HOOI work in tensor block}}
The proof of this Lemma is straight forward. 
\begin{equation*}
\begin{split}
\bcT &= \bcB \times_1 \bPi_1 \times \cdots \times_d \bPi_d\\
& = \bcB \times_1 \bPi_1 (\bPi_1^\top \bPi_1)^{-\frac{1}{2}} (\bPi_1^\top \bPi_1)^{\frac{1}{2}} \times \cdots \times_d \bPi_d (\bPi_d^\top \bPi_d)^{-\frac{1}{2}} (\bPi_d^\top \bPi_d)^{\frac{1}{2}}\\
& = \left( \bcB \times_1  (\bPi_1^\top \bPi_1)^{\frac{1}{2}}\times \cdots \times_d (\bPi_d^\top \bPi_d)^{\frac{1}{2}}\right) \times_1 \bPi_1 (\bPi_1^\top \bPi_1)^{-\frac{1}{2}} \times \cdots \times_d \bPi_d (\bPi_d^\top \bPi_d)^{-\frac{1}{2}}\\
& = \bcS \times_1 \bPi_1 (\bPi_1^\top \bPi_1)^{-\frac{1}{2}} \V_1 \times \cdots \times_d \bPi_d (\bPi_d^\top \bPi_d)^{-\frac{1}{2}} \V_d,
\end{split}
\end{equation*} where the last inequality comes from the assumption about the decomposition of $\bcB \times_1 (\bPi_1^\top \bPi_1)^{\frac{1}{2}} \times \cdots \times_d (\bPi_d^\top \bPi_d)^{\frac{1}{2}}$. 
\quad $\blacksquare$

\subsection{Proof of Theorem \ref{th: tensor block model}}
First by the same argument in the proof of Theorem \ref{th: tensor denoising}, with probability at least $1 - \exp(-c p_{\min})$, $\tau_{1k} \leq C \sigma \sqrt{p_k}$ and $\xi \leq C\sigma \sqrt{\sum_{i=1}^d p_i r_i}$. 

Under the assumption \eqref{assu: cocluster size assumption}, we have $\lambda(\bcS):= \min_i \sigma_{r_i}(\mathcal{M}_i(\bcS)) \geq C \lambda \sqrt{\frac{\prod_{i=1}^d p_i}{\prod_{i=1}^d r_i}} \geq C \sqrt{\frac{p_{\max}^2 r_{\max}}{p_{\min}}}$ by the definition of $\bcS$.

Notice that under the signal strength in Theorem \ref{th: tensor block model}, the first order perturbation error in $\| \sin \Theta(\widehat{\U}_i, \U_i)  \|$ dominates as we discussed in Remark \ref{rem: unilateral bound}. So from Corollary \ref{th:HOOI d diff} we have 
\begin{equation*}
	\| \sin \Theta(\widehat{\U}_i, \U_i)  \| \leq C \frac{\sqrt{p_k}}{\lambda/\sigma} \sqrt{\frac{\prod_{i=1}^d r_i}{\prod_{i=1}^d p_i}},
\end{equation*}
and the bound for $\|\widehat{\bcT} - \bcT\|_{\tHS}$.

For the cocluster membership recovery, first we have
\begin{equation} \label{ineq: singular subspace Fro bound}
\| \sin \Theta(\widehat{\U}_i, \U_i)  \|_F \leq \sqrt{r_i} \| \sin \Theta(\widehat{\U}_i, \U_i)  \| \leq C\frac{ \sqrt{p_i r_i} }{\lambda/\sigma} \sqrt{\frac{\prod_{i=1}^d r_i}{\prod_{i=1}^d p_i}}.
\end{equation}

The following proof is relative standard for proving misclassification error based on singular subspace estimation. Without loss of generality, we focus on the mode-$1$ clustering analysis. Take $\delta_k = \sqrt{\frac{1}{p_1^{(k)}} + \frac{1}{\max \{p_1^{(l)}: l\neq k \} } }$, and let $S_k = \{i \in G_1^{(k)}: \| (\widehat{\bPi}_1 \widehat{\X}_1)_{[i,:]} - (\U_1)_{[i,:]} \| \geq \frac{\delta_k}{2} \}$ where $\U_1$ is defined in Lemma \ref{lem: why HOOI work in tensor block}. Here we can imagine $S_k$ as the set of nodes in cluster $k$ of mode $1$ that may not be correctly clustered. By Lemma 5.3 of \cite{lei2015consistency} and \eqref{ineq: singular subspace Fro bound}, we have 
\begin{equation*}
	\sum_{k=1}^{r_1} |S_k| \delta_k^2 \leq C \frac{ p_i r_i }{(\lambda/\sigma)^2} \frac{\prod_{i=1}^d r_i}{\prod_{i=1}^d p_i}.
\end{equation*}
Moreover, under the signal strength condition in Theorem \ref{th: tensor block model}, we have that there exists a permutation matrix such that after permutation, the nodes in $G_1 = \bigcup_{k=1}^{r_1} (G_1^{(k)} \setminus S_k) $ can be perfectly recovered based on Lemma 5.3 of \cite{lei2015consistency}. So to analyze the misclassification error, we just need to consider about $S_k$.
\begin{equation*}
	 \widetilde{err}(\widehat{\bPi}_1, \bPi_1) \leq \max_{1 \leq k \leq r_1} \frac{|S_k|}{p_1^{(k)}} \leq \sum_{1 \leq k \leq r_1} \frac{|S_k|}{p_1^{(k)}} \overset{(a)}\leq \sum_{k=1}^{r_1} |S_k| \delta_k^2 \leq C \frac{ p_1 r_1 }{(\lambda/\sigma)^2} \frac{\prod_{i=1}^d r_i}{\prod_{i=1}^d p_i},
\end{equation*} where (a) is because of the choice of $\delta_k$ and 
\begin{equation*}
	 err(\widehat{\bPi}_1, \bPi_1) \leq \sum_{1 \leq k \leq r_1} \frac{|S_k|}{p_1}  \leq C \frac{ p_{1,\max}' r_1 }{(\lambda/\sigma)^2} \frac{\prod_{i=1}^d r_i}{\prod_{i=1}^d p_i}.
\end{equation*} 
The proof for other modes are similar and this finishes the proof of this theorem.
\quad $\blacksquare$
\end{sloppypar}

\vskip 0.2in
\bibliography{reference}

\end{document}